\documentclass{article}
\usepackage[a4paper,top=3cm,bottom=3.5cm,left=2.5cm,right=2.5cm,marginparwidth=1.75cm]{geometry}
\usepackage[hidelinks]{hyperref}
\usepackage{amsmath, amssymb, amsthm, bbm, graphicx, cleveref, caption, tikz, enumerate, subcaption, todonotes,comment}
\usetikzlibrary{decorations.pathreplacing,arrows.meta,calc,patterns,angles,quotes,decorations.markings}
\definecolor{mycolour}{rgb}{0,0.6,0.6}

\newcommand{\ddp}[2]{\frac{\partial#1}{\partial#2}}

\newcommand{\D}{\partial D}
\renewcommand{\S}{\mathcal{S}}
\renewcommand{\Cap}{\mathrm{Cap}}
\newcommand{\K}{\mathcal{K}}

\renewcommand*{\Re}{\operatorname{Re}}
\renewcommand*{\Im}{\operatorname{Im}}
\newcommand{\A}{\mathcal{A}}
\newcommand{\de}{\: \mathrm{d}}
\newcommand{\R}{\mathbb{R}}
\renewcommand{\i}{\mathrm{i}\mkern1mu}
\newcommand{\CC}{\mathbb{C}}
\graphicspath{{figures/}}
\newcommand{\C}{\mathcal{C}}
\newcommand{\F}{\mathcal{F}}

\newcommand{\Dc}{\mathcal{D}}
\renewcommand{\L}{\mathcal{L}}
\newcommand{\M}{\mathcal{M}}
\newcommand{\N}{\mathbb{N}}

\renewcommand{\P}{\mathcal{P}}
\newcommand{\Rc}{\mathcal{R}}
\renewcommand{\S}{\mathcal{S}}
\newcommand{\T}{\mathcal{T}}

\newcommand{\Z}{\mathbb{Z}}
\newcommand{\p}{\partial}
\renewcommand{\epsilon}{\varepsilon}
\newcommand{\dx}{\: \mathrm{d}}

\renewcommand{\b}[1]{\textbf{#1}}

\newcommand{\ie}{\textit{i.e.}}
\newcommand{\nm}{\noalign{\smallskip}}
\newcommand{\ds}{\displaystyle}
\newcommand{\uin}{u_\mathrm{in}}
\newcommand{\Ka}[2]{(\mathcal{K}_D^{-#1,#2})^*}
\renewcommand{\k}{\mathbf{k}}
\newcommand{\w}{\mathbf{w}}

\newtheorem{thm}{Theorem}
\newtheorem{prop}[thm]{Proposition}
\newtheorem{lemma}[thm]{Lemma}
\newtheorem{defn}[thm]{Definition}
\newtheorem{cor}[thm]{Corollary}
\theoremstyle{definition}
\newtheorem{remark}[thm]{Remark}

\numberwithin{equation}{section}
\numberwithin{thm}{section}

\title{Functional analytic methods for discrete approximations of subwavelength resonator systems}

\author{Habib Ammari\thanks{\footnotesize Department of Mathematics, ETH Z\"urich, R\"amistrasse 101, CH-8092 Z\"urich, Switzerland (habib.ammari@math.ethz.ch, bryn.davies@sam.math.ethz.ch, erik.orvehed.hiltunen@sam.math.ethz.ch).}\and Bryn Davies\footnotemark[1]  \and Erik Orvehed Hiltunen\footnotemark[1]}

\date{}

\begin{document}
	\maketitle

\begin{abstract}
	We survey functional analytic methods for studying subwavelength resonator systems. In particular, rigorous discrete approximations of Helmholtz scattering problems are derived in an asymptotic subwavelength regime. This is achieved by re-framing the Helmholtz equation as a non-linear eigenvalue problem in terms of integral operators. In the subwavelength limit, resonant states are described by the eigenstates of the generalised capacitance matrix, which appears by perturbing the elements of the kernel of the limiting operator. Using this formulation, we are able to describe subwavelength resonance and related phenomena. In particular, we demonstrate large-scale effective parameters with exotic values. We also show that these systems can exhibit localised and guided waves on very small length scales. Using the concept of topologically protected edge modes, such localisation can be made robust against structural imperfections.
\end{abstract}
\vspace{0.5cm}
\noindent{\textbf{Mathematics Subject Classification (MSC2000):} 35J05, 35C20, 35P20, 74J20.

\vspace{0.2cm}

	\vspace{0.5cm}

\tableofcontents

\newpage

\section{Introduction}\label{sec:intro}
\subsection{Wave manipulation at subwavelength scales}

A widespread ambition in wave physics is to be able to manipulate waves at scales that are much smaller than their wavelengths. On the other hand, an intuitive physical paradigm is that the propagation of a wave is not significantly affected by small objects or inhomogeneities. In particular, if an object is much smaller than the incident wavelength, then it will typically have a negligible scattering effect. This simplified phenomenon is closely related to Abbe’s famous diffraction limit, which describes how the resolution of imaging systems depends on the operating wavelength. In order to overcome this limit, and be able to manipulate waves at subwavelength scales, there is widespread interest in settings where small objects exhibit \emph{subwavelength resonance} and strongly scatter waves with comparatively large wavelengths.

The first high-profile example of subwavelength resonance came in the setting of acoustics when Marcel Minnaert observed the resonance of small air bubbles in water \cite{minnaert1933musical}. The very large contrast between the material parameters of air and water is understood to be the crucial mechanism here \cite{ammari2018minnaert, meklachi2018asymptotic}. This phenomenon has since been observed in a variety of other settings, such as Helmholtz resonators \cite{ammari2015superresolution},  plasmonic nanoparticles \cite{ammari2017plasmonicscalar,ammari2016plasmonicMaxwell} and high-contrast dielectric particles \cite{ammari2019dielectric,ammari2020mathematical}. In general, bounded material inclusions whose parameters differ greatly from the background medium and which experience subwavelength resonance will be referred to as \emph{subwavelength resonators} in this work.

The value of subwavelength resonators is that they can be used as the building blocks for large, complex structures which can exhibit a variety of exotic and useful properties. These micro-structured materials are examples of \emph{metamaterials}: materials with a repeating micro-structure that exhibit properties surpassing those of the individual building blocks \cite{kadic20193d}. The widespread interest in metamaterials began with the realisation that they could be designed to have effectively negative material parameters \cite{smith2004metamaterials} and, as a result, could be used to design perfect lenses \cite{veselago1968electrodynamics, pendry2000negative} as well as cloaking and shielding devices \cite{milton2006cloaking, alu2005achieving}. Moreover, due to the subwavelength nature of the resonance, these structures enable wave control on very small length scales. Most notably, waves can be confined or guided using very small devices \cite{ammari2018subwavelength,lemoult2016soda,yves2017crytalline}.

More recently, the study of micro-structured resonant media has focussed on designing structures whose properties are robust with respect to imperfections in their construction. This is important for realising the applications of this theory, since small errors will be introduced during the manufacturing process and many of the properties of these micro-structured media are very sensitive. Developments in this area have been based on studying the topological properties of periodic structures to create so-called \emph{topologically protected} modes \cite{ammari2020topological}. These concepts have previously been widely studied in a variety of settings, most notably in quantum mechanics for the Schrödinger operator \cite{fefferman2014topologically,fefferman2018honeycomb,drouot2019bulk}.

\subsection{Analysis of scattering problems}

There is a large body of work dedicated to studying the scattering of waves by a collection of objects \cite{martin2006multiple}. A popular simplifying assumption is to consider scattering by circular or spherical inclusions. In the case of a single inclusion, characterisations of the scattered field can be obtained through the use of expansions in terms of Bessel functions or spherical harmonics \cite{capdeboscq2012scattered, hansen2007asymptotically}. Likewise, in the case of two spheres a bispherical coordinate system can be used to give explicit representations of solutions \cite{ammari2020close}.

In order to study scatterers with a more general class of shapes, integral equation methods are commonly used \cite{colton1983integral}. Boundary integral formulations can be used to reduce the dimension of the scattering problem, by rephrasing it as a problem posed on the boundaries of the scatterers \cite{ammari2009layer, ammari2018minnaert}. Similarly, approaches that use Lippmann--Schwinger representations to express solutions in terms of volume integrals have been used for both scalar models \cite{meklachi2018asymptotic, ammari2019dielectric} and for the Maxwell equations \cite{costabel2010volume, costabel2015volume}. The fundamental idea here is that by representing solutions using appropriate integral operators a scattering problem can be equivalently phrased as a non-linear eigenvalue problem. With this formulation, scattering resonances can be characterised as the poles of meromorphic operator-valued functions \cite{dyatlov2019mathematical}.  In some settings, this can be paired with a scattering matrix \cite{ammari2020exceptional} or transfer matrix \cite{lin2021mathematical} formulation to give a concise description of the response of the system.

Given the multi-scale nature of subwavelength metamaterials, asymptotic techniques are often used to understand their properties. In particular, a common approach for studying subwavelength problems is to assume that the resonator is asymptotically small while the other material parameters are fixed \cite{meklachi2018asymptotic, ammari2019dielectric}. This is convenient because it can be implemented easily via a change of variables to give a concise description of a structure that is significantly (in an asymptotic sense) smaller than the operating wavelength. Related to this, homogenization techniques are often used to describe effective properties of micro-structured media. However, standard homogenization techniques do not apply here since these phenomena are based on local resonance of the small repeating units \cite{ammari2017effective,ammari2019bloch}.

One downside to modelling subwavelength resonators as being asymptotically small is that it can tend to simplify the otherwise exotic behaviour as it reduces the underlying mechanism to just a rescaling of the model. Conversely, in this work we instead fix the resonators' size and position and consider an asymptotic limit in the material contrast parameter (which, in the case of acoustic waves, describes the ratio of the density inside and outside the resonators). This has the fundamental difference that the limiting problem is not trivial and has a spectrum of eigenvalues that can be understood. The asymptotic perturbation theory of Gohberg and Sigal \cite{gohberg1971operator,ammari2018mathematical} can then be used to prove the existence of subwavelength resonant frequencies, which are defined as resonant frequencies which satisfy a given asymptotic condition. This approach reveals the fundamental differences between a system's subwavelength resonant modes and the higher-frequency resonances.

\begin{figure}
	\centering
	\begin{subfigure}[t]{0.4\linewidth}
		\includegraphics[width=\linewidth]{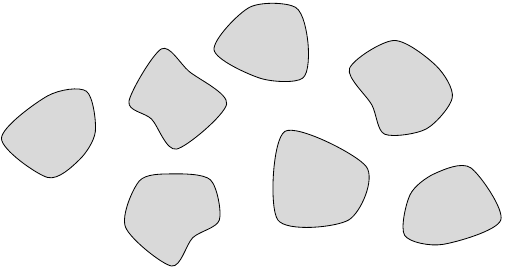}
		\caption{A system of finitely many resonators.}
	\end{subfigure}
	\hspace{0.5cm}
	\begin{subfigure}[t]{0.5\linewidth}
		\includegraphics[width=\linewidth]{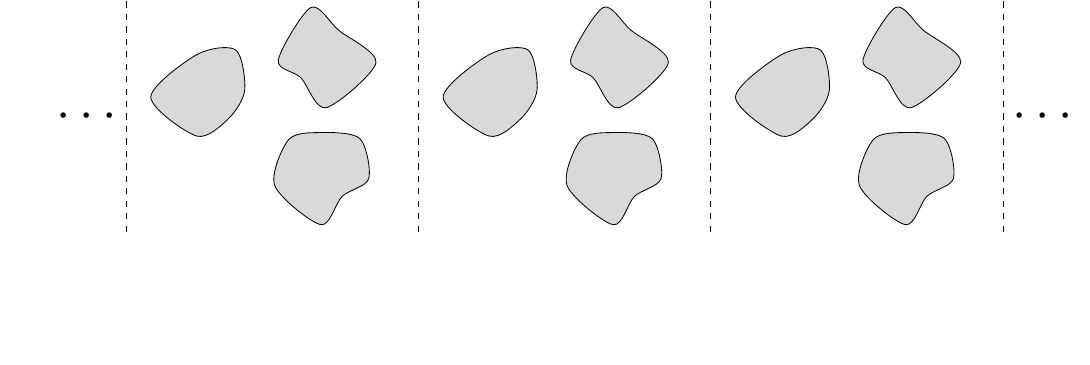}
		\caption{A periodic system of resonators.}
	\end{subfigure}
	\caption{The functional analytic method developed here is useful for studying scattering by a system of material inclusions, which act as subwavelength resonators in an appropriate high-contrast regime. We are able to derive concise asymptotic results in terms of the capacitance matrix for the case of either finitely many resonators or a periodically repeating array of finitely many resonators.}\label{fig:finiteperiodic}
\end{figure}

\subsection{Functional analytic approach} \label{sec:approach}

We represent the subwavelength resonators as material inclusions $D_i$ in $\R^d$ for $i\in\mathcal{I}$, where $\mathcal{I}\subset\N$ is some index set. The number of connected components that make up $D$ should either be finite or the geometry should be periodic, such that it is given by an array of finitely many resonators that repeats indefinitely, as illustrated in \Cref{fig:finiteperiodic}. 
We suppose that $d\in\{2,3\}$ and the material inclusions 
$D_i$ for $i\in\mathcal{I}$ are disjoint, connected sets with boundaries in $C^{1,s}$ for some $0<s<1$. We will study scalar wave settings where we throughout use $\omega$ to denote the frequency of the waves. We use $v_i$ to denote the wave speed in resonator $D_i$, then $k_i=\omega/v_i$ is the wave number in $D_i$. Similarly, the wave speed and wave number in the background medium are denoted by $v$ and $k$. We assume that $v>0$. We are interested in understanding solutions to Helmholtz resonance problems of the form
\begin{equation} \label{eq:finite_scattering}
	\left\{
	\begin{array} {ll}
		\ds \Delta {u}+ k^2 {u}  = 0 & \text{in } \R^d \setminus \overline{D}, \\
		\nm
		\ds \Delta {u}+ k_i^2 {u}  = 0 & \text{in } D_i, \text{ for } i=1,\dots,N, \\
		\nm
		\ds  {u}|_{+} -{u}|_{-}  =0  & \text{on } \partial D, \\
		\nm
		\ds  \delta_i \frac{\partial {u}}{\partial \nu} \bigg|_{+} - \frac{\partial {u}}{\partial \nu} \bigg|_{-} =0 & \text{on } \partial D_i \text{ for } i \in \mathcal{I}, \\
		\nm
		\multicolumn{2}{l}{\ds u(x) \ \text{satisfies an outgoing radiation condition},}
	\end{array}
	\right.
\end{equation}
where $\ds D=\cup_{i\in\mathcal{I}} D_i$ and the outgoing radiation condition depends on whether $D$ is a finite or a periodically infinite system of resonators.  Such Helmholtz equations, which can be used to model acoustic and polarised electromagnetic waves, represent the simplest model for wave propagation that still exhibits the rich phenomena associated to subwavelength physics.

We wish to characterise solutions to \eqref{eq:finite_scattering} in terms of the system's subwavelength resonant modes. The parameters $\delta_i$ in \eqref{eq:finite_scattering} are of crucial importance, and can be interpreted as the material contrast. They are allowed to be complex, to account for sources of energy loss or gain; see \Cref{sec:exceptional}. In order to achieve subwavelength resonance we will assume that these parameters are small, corresponding to a large contrast between the materials. So that we can perform concise asymptotics in terms of the material contrast, we will introduce the real-valued parameter $\delta:=|\delta_1|$ and assume that $\delta_i=O(\delta)$ as $\delta\to0$ for all $i\in\mathcal{I}$. We will then make a definition of a resonant mode being \emph{subwavelength} as an asymptotic property in terms of $\delta$.

\begin{defn}[Subwavelength resonant frequency] \label{defn:resonance}
	Given $\delta>0$, a subwavelength resonant frequency $\omega=\omega(\delta)\in\CC$ is defined to be such that

		$\,$(i) there exists a non-trivial solution to \eqref{eq:finite_scattering}, known as an associated resonant mode;

		(ii) $\omega$ depends continuously on $\delta$ and satisfies $\omega\to0$ as $\delta\to0$.
\end{defn}

The starting point for using functional analytic methods to understand resonance problems is to re-frame the problem \eqref{eq:finite_scattering} as an operator equation. For instance, we will show that finding a solution to \eqref{eq:finite_scattering} is equivalent to finding a (non-trivial) function $\Phi\in L^2(\p D)$ such that an integral equation of the form
\begin{equation} \label{eq:charvalproblem}
\A(\omega,\delta)[\Phi]=0,
\end{equation}
is satisfied; see \eqref{aomegadelta} and \eqref{aomegadeltalpha}. In the example studied in this work, $\A(\omega,\delta)$ will be an operator $L^2(\D)\times L^2(\D)\to H^1(\D)\times L^2(\D)$ and $\Phi$ an element of $L^2(\D)\times L^2(\D)$. Here, $H^1$ is the usual Sobolev space of square-integrable functions whose weak derivative is square integrable.

If $X$ and $Y$ are two Banach spaces, then we write $\L(X,Y)$ to denote the space of bounded linear operators from $X$ into $Y$. In this work, we are interested in the case that $X$ and $Y$ are themselves spaces of functions and we have the following definition to describe the zeros of an operator-valued function that maps into $\L(X,Y)$.

\begin{defn}[Characteristic value] \label{defn:char_value}
	A point $z_0\in\CC$ is said to be a characteristic value of $\T:\CC\to\L(X,Y)$, which is an operator-valued function of a complex variable, if there exists some $\phi\in X$ such that $\phi(z_0)\neq0$ and $\T(z_0)\phi(z_0)=0$.%\todo{holomorphic?}
%	\begin{enumerate}[(i)]
%		\item $\phi$ is holomorphic at $z_0$ and $\phi(z_0)\neq0$,
%		\item $\T(z)\phi(z)$ is holomorphic at $z_0$ and $\T(z_0)\phi(z_0)=0$.
%	\end{enumerate}
\end{defn}

Comparing Definitions~\ref{defn:resonance} and \ref{defn:char_value}, we see that finding a subwavelength resonant frequency of the system is equivalent to finding, for a given $\delta$, a characteristic value $\omega$ of $\A(\omega,\delta)$ which is such that $\omega(\delta)\to0$ as $\delta\to0$. Our approach to finding such solutions is to consider perturbations of the elements of the kernel of $\A(0,0)$. We will see that this space has dimension equal to the number of distinct resonators in the structure. Once we understand $\ker(\A(0,0))$, we can characterise characteristic values of $\A(\omega,\delta)$ for small $\omega$ and $\delta$ as perturbations of this space. This analysis is based on the asymptotic perturbation theory of Gohberg and Sigal 
(\ie, the generalised Rouch\'e theorem and argument principle  to operator-valued functions) \cite{gohberg1971operator, ammari2018mathematical} and allows us both to prove the existence of subwavelength resonant modes (satisfying \Cref{defn:resonance}) and to derive asymptotic formulas for their values.

This functional analytic approach has been used to describe subwavelength resonance in a variety of different physical settings. For instance, it was used to characterise a system of subwavelength Helmholtz resonators in \cite{ammari2015superresolution}, plasmonic particles in \cite{ammari2017plasmonicscalar} and high-contrast dielectric resonators in \cite{ammari2019dielectric,ammari2020mathematical}. In this work, we will explore its use to study scattering by a high contrast material inclusion, such as an bubble in water for the case of acoustic waves. This approach was first developed in this setting by \cite{ammari2018minnaert} but, as we shall see, has since been developed to cover a variety of different settings and applications.

\subsection{Capacitance coefficients}

In the high-contrast Helmholtz setting that we will consider here, the functional analytic method described above will yield an approximation in terms of \emph{capacitance coefficients}. Capacitance coefficients have a long history in electrostatics, where they govern the relationship between the distributions of potential and charge in a system of conductors. In particular, Maxwell introduced the matrix of capacitance coefficients $C\in\R^{N\times N}$ to be such that if $V\in\R^N$ is the vector of potentials on a system of $N$ conductors then $Q=CV$ is the vector of charges on the conductors \cite{maxwell1873treatise,diaz2011positivity}.

Capacitance coefficients appear in the setting of subwavelength Helmholtz problems when we describe the principal part of the meromorphic operator $(\A(\omega,0))^{-1}$, which is a finite-rank operator governing the perturbation of the kernel of $\A(0,0)$. Using a pole-pencil decomposition, we are able to project the problem onto $\ker(\A(0,0))$, which gives a finite-dimensional characterisation in terms of the \emph{generalised capacitance matrix}. These ideas are elaborated in \Cref{app:abstract}, and shows that the capacitance coefficients appear naturally from the functional analytic approach described in \Cref{sec:approach}.

%Capacitance coefficients are defined using the Laplace Green's function
%In the case when $D$ is a single conductor, the function $\psi_1$ represents the charge distribution which gives a unit potential on $D$, and $\Cap_D$ is known as the self capacitance.
%Thanks to their usefulness, understanding the behaviour of capacitance coefficients has been a long-standing problem.

In this article, we will survey how the generalised capacitance matrix offers a rigorous and intuitive discrete approximation to subwavelength Helmholtz scattering and resonance problems. This gives leading-order asymptotic expressions for both resonant modes and scattered solutions in terms of the eigenvalues and eigenvectors of the generalised capacitance matrix, which are accompanied by precise error bounds. We will see that a wide variety of different applications and phenomena can be studied using the capacitance approximation, demonstrating the power of reducing a differential problem to a matrix approximation in this way.

\section{Finite systems} \label{sec:finite}

\begin{figure}
	\begin{center}
		\includegraphics[width=0.6\linewidth]{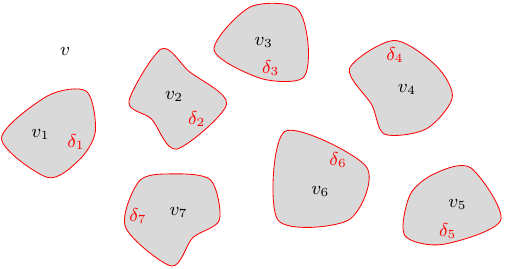}
	\end{center}
	\caption{A finite collection of $N$ resonators, with wave speeds $v_i$ for $i=1,\dots,N$, in a surrounding medium with wave speed $v$. The contrast between the $i$\textsuperscript{th} resonator and the background is given by $\delta_i$, where a small value of $\delta_i$ describes a large contrast.} \label{fig:finite}
\end{figure}

In this section, we apply the functional analytic method set out in \Cref{sec:approach} to a system of finitely many resonators. An example of the setting is sketched in \Cref{fig:finite}.
In particular,  we consider a Helmholtz resonance problem \eqref{eq:finite_scattering} for the finite domain $D=D_1\cup\dots\cup D_N$, where $N\in\N$ and the outgoing radiation condition (known as 
the Sommerfeld radiation condition) says that
\begin{equation} \label{eq:SRC}
	\lim_{|x|\to\infty} |x|^{\frac{d-1}{2}}\left(\ddp{}{|x|}-\i k\right)u=0, \quad \text{uniformly in all directions } x/|x|,
\end{equation}
and guarantees that energy is radiated outwards by the scattered solution. We assume that all contrast parameters are small while the wave speeds have order $1$. In other words, the parameter $\delta>0$ is such that
\begin{equation}
	\delta_i=O(\delta), \ v,v_i=O(1) \qquad \text{as } \delta\to0, \text{ for } i=1,\dots,N.
\end{equation}
In order to concisely represent the different $\delta_i$, we introduce the function $\widetilde\delta(x)\in L^2(\p D)$ as
\begin{equation}\label{eq:dtilde}
\widetilde{\delta}(x) = \delta_i \ \text{ for } x\in \p D_i.
\end{equation}

Due to the Sommerfeld radiation condition \eqref{eq:SRC}, the problem \eqref{eq:finite_scattering} has only a trivial solution $u=0$ for any $\omega$ real. This follows from combining Rellich's lemma (see, for instance, \cite[Section 2.8]{ammari2018mathematical}) together with the unique continuation principle.

\subsection{Main results of the capacitance formulation}\label{sec:finite_main}
The main tool that allows us to reveal the resonant properties of the system $D=D_1\cup\dots\cup D_N$ is the Helmholtz single layer potential. This is an operator $\S_D^\omega:L^2(\D)\to H_{\mathrm{loc}}^1(\R^d)$ which is defined as
\begin{equation}
\S_D^\omega[\varphi](x)=\int_{\D} G^\omega(x-y)\varphi(y)\de\sigma(y), \quad x\in\R^d, \ \varphi\in L^2(\D),
\end{equation}
where $H_{\mathrm{loc}}^1(\R^d)$ is the set of functions in $H^1(A)$ for all compact subsets $A\subset \R^d$, while $G^\omega$ is the Helmholtz Green's function, given by
\begin{equation} \label{eq:G}
G^\omega(x)=
\left\{
\begin{array}{l l}
-\frac{\i}{4} H_0^{(1)}(\omega|x|), & d=2,\\
-\frac{1}{4\pi|x|}e^{\i \omega |x|}, & d=3,
\end{array}\right.
\qquad x\neq0, \ \Re(\omega)>0.
\end{equation}
Here, $H_0^{(1)}$ is the Hankel function of the first kind and order zero. In the case $\omega=0$, $G^0$ is the Laplace Green's function given by
\begin{equation}
	G^0(x)=
	\left\{
	\begin{array}{l l}
		\frac{1}{2\pi}\ln|x|, & d=2,\\
		-\frac{1}{4\pi|x|}, & d=3,
	\end{array}\right.
	\qquad x\neq0.
\end{equation}
The single layer potential is useful because it allows us to seek solutions to \eqref{eq:finite_scattering} of the form
\begin{equation} \label{eq:representation}
u(x) = \begin{cases}
\S_D^k[\phi](x) & x\in \R^d\setminus\overline{D}, \\[0.3em]
\S_{D}^{k_i}[\psi](x) & x\in D_i,
\end{cases}
\end{equation}
where $\psi,\phi\in L^2(\D)$ are density functions that need to be found. The value of this representation is that a solution of the form \eqref{eq:representation} necessarily satisfies the Helmholtz equations and the radiation condition in problem \eqref{eq:finite_scattering}. Also, the different wave numbers $k_i$ inside $D_i$ have been taken into account by using different single layer potentials $\S_{D}^{k_i}$. We can collectively represent these single layer potentials through a single operator $\widetilde{\S}_D^{\omega}$, defined piecewise for $i=1,\dots,N$ as
\begin{equation} \label{eq:Stilde}
	\widetilde{\S}_D^{\omega}[\varphi](x) = \S_D^{k_i}[\varphi](x), \quad x\in D_i, \ \varphi\in L^2(\D).
\end{equation}
We emphasize that $\widetilde{\S}_D^{0} = \S_D^{0}$. It now remains only to find $\psi$ and $\phi$ such that the transmission conditions on $\D$ are fulfilled. This can be achieved through the introduction of an additional integral operator, the Neumann--Poincar\'e operator associated to $D$. This is an operator $\K_D^{\omega,*}$ on $L^2(\D)$ which is defined as
\begin{equation}
\K_D^{\omega,*}[\varphi](x)=\int_{\D} \ddp{}{\nu_x}G^\omega(x-y)\varphi(y)\de\sigma(y), \quad x\in\D, \ \varphi\in L^2(\D).
\end{equation}
We define $\widetilde{\K}_D^{\omega,*}$ in the same spirit as $\widetilde{\S}_D^\omega$, namely
\begin{equation} \label{eq:Ktilde}
	\widetilde{\K}_D^{\omega}[\varphi](x) = \K_D^{k_i}[\varphi](x), \quad x\in\D_i, \ \varphi\in L^2(\D).
\end{equation}
We are now able to describe how $\S_D^\omega$ and its normal derivative behave on $\D$. In particular, it holds that for any $\varphi\in L^2(\D)$ (see, for instance, \cite{polarization,nedelec,taylor})
\begin{equation} \label{eq:jump1}
\S_D^\omega[\varphi]\big|_+=\S_D^\omega[\varphi]\big|_-,
\end{equation}
and
\begin{equation} \label{eq:jump2}
\ddp{}{\nu}\S_D^\omega[\varphi]\Big|_\pm=\left(\pm \frac{1}{2}I+\K_D^{\omega,*}\right),
\end{equation}
where the subscripts $+$ and $-$ denote taking the limit from outside and inside the boundary $\D$, respectively. With the so-called \emph{jump conditions} \eqref{eq:jump1} and \eqref{eq:jump2} in hand, we can derive the following lemma, which characterises the resonance problem \eqref{eq:finite_scattering} as a boundary integral equation.

\begin{lemma} \label{lem:BIE_finite}
	In the regime $\omega\rightarrow 0$, the Helmholtz problem \eqref{eq:finite_scattering} is equivalent to finding $\psi,\phi\in L^2(\D)$ such that
	\begin{equation} \label{eq:BIE_finite}
	\A(\omega,\delta)
	\begin{pmatrix} \psi \\ \phi \end{pmatrix}
	= \begin{pmatrix} 0 \\ 0 \end{pmatrix},
	\end{equation}
	where the operator $\A(\omega,\delta):L^2(\D)\times L^2(\D)\to H^1(\D)\times L^2(\D)$ is defined as
	\begin{equation}\label{aomegadelta}
	\A(\omega,\delta)= \begin{pmatrix}
	\widetilde{\S}_D^\omega & -\S_D^k \\
	-\frac{1}{2}I+\widetilde{\K}_D^{\omega,*} & -\widetilde\delta \left(\frac{1}{2}I+\K_D^{k,*}\right)
	\end{pmatrix},
	\end{equation}
where, as in \eqref{eq:dtilde}, $\widetilde{\delta}(x) = \delta_i$ for $x\in \p D_i$.
\end{lemma}
The Helmholtz problem \eqref{eq:finite_scattering} is equivalent to the characteristic value problem \eqref{eq:BIE_finite} provided that $\omega/v_i$ is not a Dirichlet eigenvalue of $D_i$ for any $i$. This condition is naturally satisfied in the regime $\omega\rightarrow 0$.

The approach outlined in \Cref{sec:approach} can now be applied to \Cref{lem:BIE_finite} to prove the existence of subwavelength resonances, as defined in \Cref{defn:resonance}, and derive their asymptotic behaviour as $\delta \to 0$. The idea here is to study the kernel of $\A(0,0)$, where
\begin{equation} \label{eq:A00}
\A(0,0)= \begin{pmatrix}
{\S}_D^0 & -\S_D^0 \\
-\frac{1}{2}I+{\K}_D^{0,*} & 0
\end{pmatrix},
\end{equation}
and then understand how $\ker\A(0,0)$ is perturbed when $\delta$ and $\omega$ are non-zero. The following lemma describes the two operators that appear in $\A(0,0)$, as given in \eqref{eq:A00}.

\begin{lemma} \label{lem:A00properties}
	Consider a system of $N$ subwavelength resonators $D=D_1\cup\dots D_N$ in $\mathbb{R}^3$. Then, it holds that

	(i) the Laplace single layer potential $\mathcal{S}_D^0:L^2(\partial D)\to H^1(\partial D)$ is invertible,

	(ii) $\ker(-\frac{1}{2}I+\mathcal{K}_D^{0,*})=\mathrm{span}\{\psi_1,\psi_2,\dots,\psi_N\}$ where $\psi_i:=(\mathcal{S}_D^0)^{-1}[\chi_{\partial D_i}]$, and $\chi_{\p D_i}$ denotes the characteristic function of $\p D_i$, for $i=1, \dots, N$.
\end{lemma}

From \Cref{lem:A00properties} we can see that $\A(0,0)$ has an $N$-dimensional kernel. Therefore, $\omega=0$ is a characteristic value of $\A(\omega,0)$. Due to symmetry, the multiplicity of $\omega=0$ is, in fact, $2N$. When the material parameters are real, it is easy to see that $\overline{\A(\omega,\delta)}=\A(-\overline{\omega},\delta)$, from which we can see that the resonant frequencies will be symmetric with respect to the imaginary axis, in the sense described in \Cref{lem:symmetric_res} (\emph{cf.} the analysis of \cite{dyatlov2019mathematical}).

\begin{lemma} \label{lem:symmetric_res}
	The set of resonant frequencies is symmetric in the imaginary axis. In particular, if $\delta_i, v_i\in \R$ for all $i=1,\dots,N$, and if $\omega$ is such that \eqref{eq:BIE_finite} is satisfied for some non-zero $\psi,\phi\in L^2(\D)$, then it will also hold that
	\begin{equation*}
	\A(-\overline{\omega},\delta)
	\begin{pmatrix} \overline{\psi} \\ \overline{\phi} \end{pmatrix}
	= \begin{pmatrix} 0 \\ 0 \end{pmatrix}.
	\end{equation*}
\end{lemma}
With \Cref{lem:symmetric_res} in mind, we will subsequently state results only for the resonant frequencies with non-negative real parts. We can now show the following two theorems, using the approach described in \Cref{sec:approach}.

\begin{thm}
	Consider a system of $N$ subwavelength resonators in $\mathbb{R}^d$ for $d\in\{2,3\}$. For sufficiently small $\delta>0$, there exist $N$ subwavelength resonant frequencies $\omega_1(\delta),\dots,\omega_N(\delta)$ with non-negative real parts.
\end{thm}

\begin{defn}[Capacitance matrix] \label{defn:CM}
	For a system of $N\in\N$ resonators $D_1,\dots,D_N$ in $\R^3$ we can define the capacitance matrix $C=(C_{ij})\in\R^{N\times N}$ to be the square matrix given by
	\begin{equation*}
	C_{ij}=-\int_{\D_i} (\S_D^0)^{-1}[\chi_{\p D_j}] \de\sigma,\quad i,j=1,\dots,N.
	\end{equation*}
\end{defn}
Due to the different material parameters inside each resonator we introduce the generalised capacitance matrix, which is the main quantity we use in order to describe the subwavelength resonators.
\begin{defn}[Generalised capacitance matrix] \label{defn:GCM}
	For a system of $N\in\N$ resonators $D_1,\dots,D_N$ in $\R^3$ we can define the generalised capacitance matrix, denoted by $\C=(\C_{ij})\in\CC^{N\times N}$, to be the square matrix given by
	\begin{equation}\label{eq:GCM}
	\C_{ij}=\frac{\delta_i v_i^2}{|D_i|} C_{ij}, \quad i,j=1,\dots,N.
	\end{equation}
\end{defn}

\begin{thm} \label{thm:res}
	Let $d=3$. Consider a system of $N$ subwavelength resonators in $\mathbb{R}^3$. As $\delta\to0$, the $N$ subwavelength resonant frequencies satisfy the asymptotic formula
	\begin{equation*}
	\omega_n = \sqrt{\lambda_n}+O(\delta), \quad n=1,\dots,N,
	\end{equation*}
	where $\{\lambda_n: n=1,\dots,N\}$ are the eigenvalues of the generalised capacitance matrix $\mathcal{C}\in\mathbb{C}^{N\times N}$, which satisfy $\lambda_n=O(\delta)$ as $\delta\to0$.
\end{thm}

\begin{remark}
	The assumption that the dimension $d=3$ in Definitions~\ref{defn:CM} and \ref{defn:GCM} is important as the Laplace single layer potential $\S_D^0$ is known to be invertible in this case. As we will see in \Cref{sec:twodim}, this is not generally the case when $d=2$ meaning that the corresponding version of \Cref{thm:res} is slightly less elegant.
\end{remark}

\begin{cor} \label{prop:eigenvector}
	Let $d=3$. Let $\textbf{v}_n$ be the normalised eigenvector of $\mathcal{C}$ associated to the eigenvalue $\lambda_n$. Then the normalised resonant mode $u_n$ associated to the resonant frequency $\omega_n$ is given, as $\delta\to0$, by
	\begin{equation*}
	u_n(x)=\begin{cases}
	\textbf{v}_n\cdot\textbf{S}_D^{k}(x)+O(\delta^{1/2}), \quad x\in\mathbb{R}^3\setminus \overline{D}, \\
	\textbf{v}_n\cdot\textbf{S}_D^{k_i}(x)+O(\delta^{1/2}), \quad x\in D_i,
	\end{cases}
	\end{equation*}
	where $\textbf{S}_D^k:\mathbb{R}^3\to\mathbb{C}^N$ is the vector-valued function given by
	\begin{equation*}
	\textbf{S}_D^k (x)=\begin{pmatrix}
	\mathcal{S}_D^k[\psi_1](x) \\[-0.4em]
	\vdots \\[-0.3em]
	\mathcal{S}_D^k[\psi_N](x)
	\end{pmatrix},  \quad x\in\mathbb{R}^3\setminus \partial D,
	\end{equation*}
	with $\psi_i:=(\mathcal{S}_D^0)^{-1}[\chi_{\partial D_i}]$.
\end{cor}

\begin{remark}
	The capacitance matrix is defined solely in terms of the kernel of the integral operators $\A(0,0)$ and $\A^*(0,0)$. In fact, the adjoint $\K_D^0$ of the Neumann--Poincar\'e operator satisfies $$\ker(-\frac{1}{2}I+\K_D^0) = \mathrm{span}\{\chi_{\p D_1},\chi_{\p D_2},\dots,\chi_{\p D_N}\}.$$ Then,
	$$C_{ij} = -\left\langle \chi_{\p D_i}, \psi_j\right\rangle,$$
	where, as in \Cref{lem:A00properties}, $\{\psi_1,\dots,\psi_N\}$ is a basis for $\ker(-\frac{1}{2}I+\K_D^{0,*})$. In \Cref{app:abstract}, we use these ideas to define the generalised capacitance matrix purely in terms of the integral operator $\A$, thus providing a general method to study subwavelength resonance systems.
\end{remark}

\subsection{Properties of the capacitance matrix}

Through \Cref{thm:res} and \Cref{prop:eigenvector}, we have reduced the resonance problem \eqref{eq:finite_scattering} to a matrix eigenproblem for the generalised capacitance matrix $\C$. We now wish to understand the properties of $\C$.

\begin{lemma} \label{lem:capacitance_altdefn}
	For $i=1,\dots,N$, let $V_i$ be defined as the solution to the exterior boundary value problem
	\begin{equation*}
	\left\{
	\begin{array} {ll}
	\ds \Delta V_i  = 0 & \text{in } \R^3 \setminus \overline{D}, \\
	\nm
	\ds  V_i=\delta_{ij}  & \text{on } \partial D_j, \text{ for } j=1,\dots,N, \\
	\nm
	\ds V_i(x)=O\left(|x|^{-1}\right) & \text{as } |x|\to\infty,
	\end{array}
	\right.
	\end{equation*}
	where $\delta_{ij}$ is the Kronecker delta. Then, the capacitance coefficients, defined in \Cref{defn:CM}, are given by
	\begin{equation*}
	C_{ij}=\int_{\R^3\setminus D} {\nabla V_i}\cdot \nabla V_j \de x, \quad\text{for } i,j=1,\dots,N.
	\end{equation*}
\end{lemma}

From \Cref{lem:capacitance_altdefn} emerges a slightly different explanation for why the capacitance matrix approximation works. In the limiting case, when $\delta_i=0$ for all $i$, the Helmholtz problem \eqref{eq:finite_scattering} is reduced to a Neumann boundary value problem in the interior of $D$ and a Dirichlet boundary value problem in the exterior of $D$. When $\omega=0$ (\ie,  $k=k_i=0$), the interior Neumann problem is solved by constant functions, meaning that the solution of the exterior Dirichlet boundary value problem is a linear combination of the functions $V_1,\dots,V_N$.

The projection onto this finite dimensional space yields a leading-order approximation of the solution as $\omega,\delta\to0$, in the form of an eigenvalue problem for the generalised capacitance matrix. \Cref{lem:capacitance_altdefn} is also useful as it allows us to immediately see, among other things, the symmetry of the capacitance matrix.

\begin{lemma}
	The capacitance matrix $C$ is symmetric and positive definite.
\end{lemma}

The symmetry and positive definiteness of the capacitance matrix $C$ is useful for understanding the properties of the generalised capacitance matrix $\C$ which is the product of $C$ with a diagonal matrix containing the weights $\delta_iv_i^2/|D_i|$. In the case that $\delta_iv_i^2>0$, this diagonal matrix is positive definite so, for example, we have the following lemma, which can be proved using the fact that $C$ is always Hermitian.

\begin{lemma} \label{lem:eigenbasis}
	If $\delta_iv_i^2$, for $i=1,\dots,N$, are real-valued positive numbers, then the generalised capacitance matrix $\C\in\CC^{N\times N}$ has $N$ linearly independent eigenvectors.
\end{lemma}

\begin{remark}
	We will see, in \Cref{sec:exceptional}, that non-zero imaginary parts of $\delta_iv_i^2$ can be used to model damping and amplification in the system. In this case, we can create \emph{exceptional points} where eigenvalues and eigenvectors coincide and $\C$ is not diagonalisable.
\end{remark}

\begin{figure}
	\begin{center}
		\includegraphics[width=0.5\linewidth]{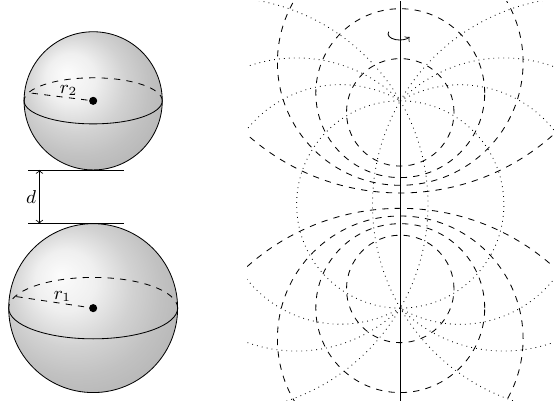}
	\end{center}
	\caption{A system of two spherical resonators can be described using bispherical coordinates. Such a coordinate system is convenient since the boundaries of the spheres lie on level sets and the capacitance coefficients can be calculated explicitly.} \label{fig:twoshperes}
\end{figure}

%\subsection{Formulas for capacitance coefficients}

Explicit formulas for capacitance coefficients are generally beyond reach. However, by making some additional assumptions, we can make the capacitance matrix easier to understand. For example, spherical resonators can be easily described using appropriate radial coordinate systems. In this way, we can see that if $D$ is a single sphere of radius $r$ then its capacitance is given by $\Cap_D := -\int_{\p D} (\mathcal{S}_{D}^0)^{-1}[\chi_{\p D}] \dx \sigma= 4\pi r$. Similarly, if we consider the case of two spherical resonators, as depicted in \Cref{fig:twoshperes}, then we can use a bispherical coordinate system to derive explicit formulas for the capacitance coefficients \cite{ammari2020close}.

\begin{lemma} \label{lemmablowup}
	Suppose that $D=D_1\cup D_2$ consists of two spheres of radius $r_1$ and $r_2$ separated by a distance $d >0$. Define the quantities $\alpha$, $\xi_1$ and $\xi_2$ as
	\begin{equation*}
	\alpha=\frac{\sqrt{d(2r_1+d)(2r_2+d)(2r_1+2r_2+d)}}{2(r_1+r_2+d)}
	\quad\text{and}\quad
	\xi_i=\sinh^{-1}\left(\frac{\alpha}{r_i}\right), \text{ for } i=1,2.
	\end{equation*}
	Then, it holds that
	\begin{align*}
	C_{11}=8\pi\alpha\sum_{n=0}^{\infty} \frac{e^{(2n+1)\xi_2}}{e^{(2n+1)(\xi_1+\xi_2)}-1}, \qquad
	C_{22}=8\pi\alpha\sum_{n=0}^{\infty} \frac{e^{(2n+1)\xi_1}}{e^{(2n+1)(\xi_1+\xi_2)}-1},
	\end{align*}
	\begin{equation*}
	C_{12}=C_{21}=-8\pi\alpha\sum_{n=0}^{\infty}\frac{1}{e^{(2n+1)(\xi_1+\xi_2)}-1}.
	\end{equation*}
\end{lemma}

\begin{remark} \label{rmk:blowup}
Based on \Cref{lemmablowup}, an asymptotic  analysis  of  the  behaviour  of  two  subwavelength  spherical resonators that  are  close  to  touching (\ie, as $d\rightarrow 0$) can be conducted.  In  \cite{ammari2020close}, it is shown that  the  two  subwavelength  resonant frequencies  associated to the two sphere system have  different  asymptotic  behaviours  and   estimates  for  the rate at which the gradient of each eigenmode blows up are derived.
In the acoustic setting,  the gradient of the
acoustic pressure describes the forces that the resonators exert on one another in the presence of sound waves. When the subwavelength resonant frequencies are excited,  enhancement of the forces in the gap region between the two spherical resonators is observed. This field enhancement
phenomena is due to subwavelength resonances and is similar to the one observed in electromagnetics for nearly-touching high contrast dielectric resonators \cite{dielectricblowup}, where subwavelength resonances occur \cite{ammari2020mathematical}. The results obtained in \cite{ammari2020close} could be generalised to shapes that are strictly convex in a region of the close-to-touching points. This  relies  on  using  spheres  with  the  same  curvature  to  approximate  the structure,  as  has  been  done in  the setting of antiplane  elasticity in \cite{ammari2013}. 
\end{remark}

In the case of larger systems of resonators we cannot hope to find such concise representations for the capacitance coefficients. However, a very useful property is that if we multiply elements in some domain $B\subset\R^3$ by some factor $a\in\R$, then a scaling argument can be used to see that $\textrm{Cap}_{a B} = a \textrm{Cap}_B$. With this in mind, we can obtain explicit expressions in the case when the resonators are small compared to the distance between them. The following lemma follows from appropriate scaling arguments, where we fix the resonators and scale the distances between them \cite{ammari2020high}.

\begin{lemma} \label{lem:dilutefinite}
	For $j=1,\dots,N$, let $B_j$ be fixed, bounded subsets of $\R^3$ with boundary in $C^{1,s}$ for some $0<s<1$. Then, consider a dilute system of $N$ subwavelength resonators with size of order $\epsilon$, given by
		\begin{equation*}
		D=\bigcup_{j=1}^N \left(B_j + \epsilon^{-1} z_j\right),
		\end{equation*}
		where $0<\epsilon\ll1$ and $z_j\in\R^3$ are fixed vectors that describe the relative position of each resonator. In the limit as $\epsilon\rightarrow 0$ the asymptotic behaviour of the capacitance matrix is given by
		$$
		C_{ij} =
		\begin{cases}
		\displaystyle \mathrm{Cap}_{B_i} + O(\epsilon^2), &\quad i=j,\\[0.3em]
		\displaystyle -\frac{\epsilon\mathrm{Cap}_{B_i}\mathrm{Cap}_{B_j}}{4\pi|z_i - z_j|} + O(\epsilon^2), &\quad i\neq j.\\
		\end{cases}
		$$
\end{lemma}

\begin{remark}
In \cite{fepponmodal}, several other important properties of the
        capacitance matrix,  
        such as a Perron-Frobenius type theorem, spectral bounds and
        properties on the coefficients in the case of symmetries, are established.
\end{remark}
\begin{remark}\label{rmk:tightbind}
	It is interesting to compare and contrast the capacitance formulation to the tight-binding approximation that is commonly employed in quantum-mechanical settings \cite{fefferman2018honeycomb,wallace1947band}. Both these formulations construct matrix eigenvalue problems as discrete approximations to continuous differential problems. A crucial difference, however, is that the generalised capacitance matrix accounts for strong interactions between the resonators. The analogy with the tight-binding model is closer when the resonators are dilute. In this case, we see from \Cref{lem:dilutefinite} that $C_{ij}$ only depends on the $i$\textsuperscript{th} and $j$\textsuperscript{th} resonator and is unaffected by remaining resonators. Moreover, in the dilute regime, the eigenmodes of the system can be approximated by a linear combination of the eigenmodes of the individual resonators. This property, which is a key assumption in the tight-binding approximation, does not hold in the case of non-dilute subwavelength resonators.
\end{remark}

\subsection{Modal decompositions}
The solution to the resonance problem, given in \Cref{thm:res} and \Cref{prop:eigenvector}, can be used to understand the scattering behaviour of $D$. That is, we can use an expansion in terms of the resonant modes $u_1,\dots,u_N$ (\textit{i.e.} a modal decomposition) to express the scattered field when $D$ is illuminated by some incident wave $\uin$. We therefore consider the problem
\begin{equation} \label{eq:scattering_eq}
	\left\{
	\begin{array} {ll}
		\ds \Delta {u}+ k^2 {u}  = 0 & \text{in } \R^3 \setminus \overline{D}, \\
		\nm
		\ds \Delta {u}+ k_i^2 {u}  = 0 & \text{in } D_i, \text{ for } i=1,\dots,N, \\
		\nm
		\ds  {u}|_{+} -{u}|_{-}  =0  & \text{on } \partial D, \\
		\nm
		\ds  \delta_i \frac{\partial {u}}{\partial \nu} \bigg|_{+} - \frac{\partial {u}}{\partial \nu} \bigg|_{-} =0 & \text{on } \partial D_i \text{ for } i=1,\dots,N, \\
		\nm
		\multicolumn{2}{l}{\ds u-\uin \ \text{satisfies the Sommerfeld radiation condition}.}
	\end{array}
	\right.
\end{equation}
Here, the frequency $\omega$ of the incident field is real, $u$ is the total field while $u-\uin$ is the scattered field. We assume that the incident field satisfies $\Delta \uin + k^2\uin = 0$ in $\R^d$ and $\nabla \uin\big|_D = O(\omega)$.  The next result, from \cite{ammari2020biomimetic} (see also \cite{fepponmodal}), shows the modal decomposition approximation of the scattered field.
\begin{thm} \label{lem:modal_res}
	Let $V$ be the matrix of eigenvectors of $\C$. If $\omega=O(\sqrt{\delta})$ as $\delta \to 0$ and $|\omega -\omega_i| > K\sqrt{\delta}$ for $i=1,\dots,N$, for some constant $K>0$, then the solution to the scattering problem \eqref{eq:scattering_eq} can be written, uniformly for $x$ in compact subsets of $\mathbb{R}^3$, as
	\begin{equation*}
		u(x)-\uin(x) = \sum_{n=1}^N a_n u_n(x) - \S_D^k\left[\left(\S_D^k\right)^{-1}[\uin]\right](x) + O(\sqrt{\delta}),
	\end{equation*}
	for coefficients $a_n=a_n(\omega)$ which satisfy the problem
	\begin{equation*}
	V\begin{pmatrix}
	\omega^2-\omega_1^2 & & \\ & \ddots & \\ & & \omega^2-\omega_N^2
	\end{pmatrix}
	\begin{pmatrix} a_1 \\ \vdots \\ a_N \end{pmatrix}
	=
	\begin{pmatrix}
	\frac{\delta_1v_1^2}{|D_1|} \int_{\D_1} (\S_D^0)^{-1}[\uin]\de\sigma \\
	\vdots \\
	\frac{\delta_Nv_N^2}{|D_N|} \int_{\D_N} (\S_D^0)^{-1}[\uin]\de\sigma
	\end{pmatrix}+O(\delta^{3/2}).
	\end{equation*}
\end{thm}

\begin{remark}
	The term $\S_D^k\left[\left(\S_D^k\right)^{-1}[\uin]\right](x)$ in \Cref{lem:modal_res} is equal to $\uin(x)$ if $x\in D$ but not for $x$ outside of the resonators.
\end{remark}

\subsection{Higher-order approximations}

The arguments used to derive the asymptotic formula in \Cref{thm:res} can be continued to higher orders. For details, see \cite{ammari2020biomimetic}. For simplicity, we assume that the material parameters on each resonator are the same.

\begin{thm} \label{res_impart}
	Let $d=3$. Consider a system of $N$ subwavelength resonators in $\mathbb{R}^3$. Suppose that the material parameters are the same on each resonator, \textit{i.e.} $v_1=v_2=\dots=v_N$ and $\delta_1=\delta_2=\dots=\delta_N$. As $\delta\to0$, the $N$ subwavelength resonant frequencies satisfy the asymptotic formula
	\begin{equation*}
	\omega_n = \sqrt{\lambda_n}-\i \tau_n+O(\delta^{3/2}), \quad n=1,\dots,N,
	\end{equation*}
	where $\lambda_n$ for $n=1,\dots,N$ are the eigenvalues of the generalised capacitance matrix $\mathcal{C}$ and $ \tau_n$ are given by
	\begin{equation*}
	\tau_n=\delta_1 \frac{v_1^2}{8\pi v} \frac{\textbf{v}_n^\top C J C\textbf{v}_n}{\| \textbf{v}_n\|_D^2},
	\end{equation*}
	with $C$ being the capacitance matrix, $J$ the $N\times N$ matrix of ones, $\textbf{v}_n$ the eigenvector associated to $\lambda_n$ and we use the norm $\| x\|_D:=\big(\sum_{i=1}^N |D_i| x_i^2\big)^{1/2}$. Further, for each $n=1,\dots,N$, it holds that $\sqrt{\lambda_n}=O(\delta^{1/2})$ and $\tau_n=O(\delta)$ as $\delta\to0$.
\end{thm}

\begin{remark}
	If the material parameters $v_1,\dots,v_N$ and $\delta_1,\dots,\delta_N$ are real, then $\lambda_n$ and $\tau_n$ from \Cref{res_impart} are all non-negative real numbers. This follows from the fact that the capacitance matrix $C$ is symmetric and positive definite. Thus, in this case the $O(\delta^{1/2})$-term is the leading-order approximation of the real part while the imaginary part appears at $O(\delta)$.
\end{remark}

\begin{remark}
	Due to the loss of energy (\emph{e.g.} to the far field), the resonant frequencies will have negative imaginary parts when the material parameters are real. In many cases it will hold that $\tau_n=0$ for some $n$, meaning that the imaginary parts exhibit higher-order behaviour in $\delta$. For example, the imaginary part of the second (dipole) frequency for a pair of identical resonators with real parameters is known to be $O(\delta^{2})$ \cite{ammari2019double}.
\end{remark}

\subsection{Two-dimensional models} \label{sec:twodim}
Throughout \Cref{sec:finite} we have mainly considered the problem of a resonator array in $\R^3$. This was convenient for two reasons. Firstly, for small frequencies the Laplace single layer potential $\S_D^0$ approximates the Helmholtz single layer potential $\S_D^\omega$ at leading order, in the sense that $\S_D^\omega=\S_D^0+O(\omega)$ in the operator norm as $\omega\to0$. On top of this, the fact that $\S_D^0$ is invertible in three dimensions was central to our definition of the capacitance matrix. If we consider a Helmholtz problem in two dimensions, however, we do not have either of these helpful properties. In two dimensions, $\S_D^0$ is not generally injective and the low-frequency expansion of $\S_D^\omega$ is given by
\begin{equation} \label{eq:expansion2d}
\S_D^\omega = \frac{1}{2\pi}\log\left(\frac{1}{2}\omega e^{\gamma-\i\frac{\pi}{2}}\right) I_{\D} + \S_D^0 + O(\omega^2\log\omega), \quad\text{as }\omega\to0,
\end{equation}
where $I_{\D}$ is the map defined as $$I_{\D}[\varphi]=\int_{\D}\varphi\de\sigma$$ for $\varphi\in L^2(\D)$ and $$\gamma=\lim_{n\to\infty}(\sum_{k=1}^n \frac{1}{k}-\log n)\approx0.577\dots$$ is the Euler--Mascheroni constant. That is, the leading-order term in the expansion of $\S_D^\omega$ has a $\log\omega$ singularity as $\omega\to0$.

The invertibility of the Laplace single layer potential in two dimensions can be readily fixed. Let $L^2_0(\p D)$ be the mean-zero subspace of $L^2(\D)$ defined as
\begin{equation*}
L^2_0(\p D) = \left\{f \in L^2(\p D) \ : \ \int_{\p D} f \dx \sigma = 0 \right\}.
\end{equation*}
Then, we have that the Laplace single layer potential is well-behaved on $L^2_0(\p D)$ \cite{ammari2019cochlea}.

\begin{lemma}
	Let $d=2$. The Laplace single layer potential $\S_D^0$ is invertible from $L^2_0(\p D)$ onto its image.
\end{lemma}

With this in hand, we can show that an invertible version of $\S_D^0$ can be defined by adding a term proportional to the integral operator $I_{\D}[\varphi]=\int_{\D}\varphi\de\sigma$. For example, we have the following lemma \cite{ammari2019cochlea}.
\begin{lemma}
	Let $d=2$ and $I_{\D}$ be the integral map given by $I_{\D}[\varphi]=\int_{\D}\varphi\de\sigma$. For any $\omega\in\CC\setminus\{z\in\CC:z=\i y \text{ for some } y\geq0\}$, the operator $\widehat{\S}_D^\omega:L^2(\D)\to L^2(\D)$, defined as
	\begin{equation*}
	\widehat{\S}_D^\omega=\frac{1}{2\pi}\log\left(\frac{1}{2}\omega e^{\gamma-\i\frac{\pi}{2}}\right) I_{\D} + \S_D^0,
	\end{equation*}
	is invertible.
\end{lemma}

Notice, finally, that $\widehat{\S}_D^\omega$ is nothing other than the leading-order approximation of $\S_D^\omega$ from \eqref{eq:expansion2d}. This means that, up to some technical modifications, we can repeat the argument used to derive \Cref{thm:res} to obtain the following analogous result.
\begin{thm}
	Let $d=2$. A system of $N$ subwavelength resonators in $\mathbb{R}^2$ has $N$ subwavelength resonant frequencies. Further to this, if for any $\omega$ and $\delta$ we define the $N\times N$-matrix $\A^{(2)}_{\omega,\delta}$ as
	\begin{equation} \label{adimension2}
	\begin{split}
	(\A^{(2)}_{\omega,\delta})_{ij} &=\omega^2\log\omega
	+\left(\left(1+\frac{c_1}{b_1}-\log v_i\right)
	-\frac{\mathcal{S}_D[\psi_j]|_{\partial D_i}}{4b_1( \int_{\partial D}\psi_j )} \right)\omega^2 \\
	&\hspace{12em}-\frac{v_i^2}{4b_1|D_i|}\left(\frac{\int_{\partial D_i}\psi_j }{\left(\int_{\partial D}\psi_j\right)}
	+\frac{\log (v/v_i)}{2\pi} \int_{\partial D_i} (\hat{\mathcal{S}}_D^k)^{-1}[\chi_{\partial D}]\right)\delta_i,
	\end{split}
	\end{equation}
	where $b_1=-\frac{1}{8\pi}$ and $c_1=b_1(\gamma -\log 2-1-\i\frac{\pi}{2})$, then the subwavelength resonant frequencies are such that the determinant of $\A^{(2)}_{\omega,\delta}$ vanishes, at leading order:
	\begin{equation*}	\det(\A^{(2)}_{\omega,\delta})=O(\omega^4\log\omega+\delta\omega^2\log\omega), \quad\text{as } \omega,\delta\to0.
	\end{equation*}
\end{thm}

%\begin{equation*}
%B_\delta^{(i)}(\omega)[\phi]:= \omega^2\log\omega
%+\left(\left(1+\frac{c_1}{b_1}-\log v_i\right)
%-\frac{\mathcal{S}_D[\phi]|_{\partial D_i}}{4b_1( \int_{\partial D}\phi )} \right)\omega^2 -\frac{v_i^2}{4b_1|D_i|}\left(\frac{\int_{\partial D_i}\phi }{\left(\int_{\partial D}\phi\right)}
%+\frac{\log (v/v_i)}{2\pi} \int_{\partial D_i} (\hat{\mathcal{S}}_D^k)^{-1}[\chi_{\partial D}]\right)\delta_i,
%\end{equation*}

\begin{remark} \label{rmklogcapacity}
    Notice that $\mathcal{S}_D[\psi_j]$ is constant in $D_j$ since $\psi_j \in \ker(-\frac{1}{2}I+\mathcal{K}_D^{0,*})$. Then 
    $\ds - \frac{\mathcal{S}_D[\psi_j]|_{\partial D_j}}{\int_{\partial D}\psi_j}$ is nothing else than $1/(2\pi) \times$ the logarithm of the capacity of $D_j$ (see, for instance, \cite[p. 39]{polarization}). Furthermore, the matrix with entries 
$$\ds - \frac{\mathcal{S}_D[\psi_j]|_{\partial D_i}}{\int_{\partial D}\psi_j}$$ in \eqref{adimension2} can be considered as the two-dimensional analogue of the capacitance matrix introduced in \Cref{defn:CM}.
\end{remark}

\subsection{Numerical approaches}

\begin{figure}
	\begin{center}
		\includegraphics[width=0.65\linewidth]{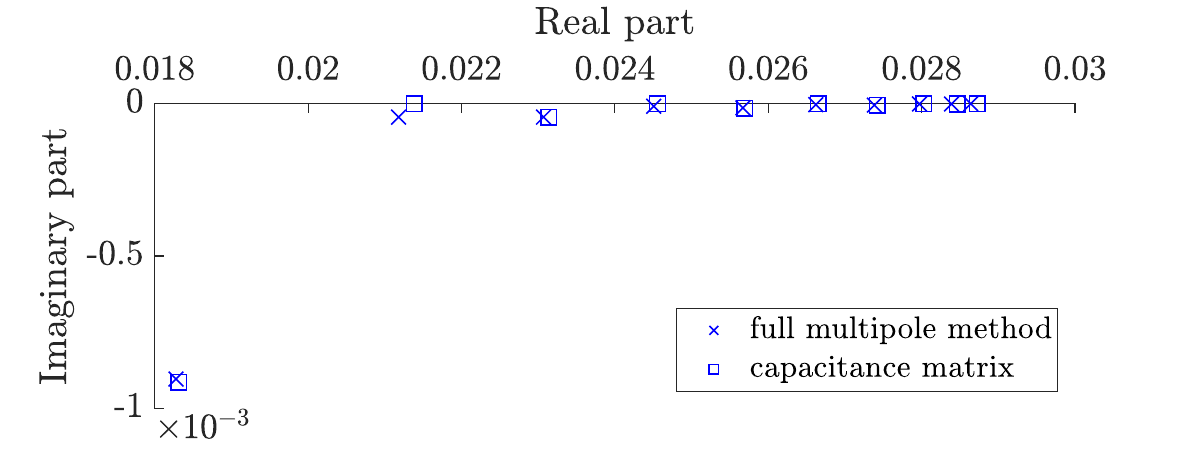}
	\end{center}
	\caption{The subwavelength resonant frequencies of a system of ten spherical resonators. We compare the values computed using the multipole expansion method to discretised the full boundary integral equation and the values computed using the capacitance matrix. The computations using the full multipole method took $41$ seconds while the approximations from the capacitance matrix took just $0.02$ seconds, on the same computer. Each resonator has unit radius and we use $\delta=1/5000$.} \label{fig:frequencies}
\end{figure}

The resonant frequencies of \eqref{eq:finite_scattering} can be computed numerically in a variety of ways. For example, we can make use of the boundary integral formulation \eqref{eq:BIE_finite} to derive a discrete version of the problem. This can be achieved, for example, by discretizing the boundary using boundary elements or a multipole expansion. In order to find the resonant frequencies for a given value of $\delta$, one needs to use a numerical root finding algorithm to find $\omega$ such that the boundary integral equation \eqref{eq:BIE_finite} is satisfied. The discrete version of the boundary integral operator (which depends non-linearly on $\omega$) will need to be recomputed at each step in this iterative algorithm.

Muller's method can be used as a numerical root finding algorithm, following the methodology of  \cite[Section 1.6]{ammari2009layer}. 
Open source codes for its implementation are provided in \cite{ ammari2018mathematical,muller}. We recall that Muller's method allows to find (complex) roots of holomorphic functions by using quadratic interpolants. 
In \Cref{fig:frequencies}, we apply Muller's method to obtain the zeros of the eigenvalues of the discrete matrix approximation of the operator 
$\A(\omega,\delta)$ in \eqref{eq:charvalproblem} derived from using the multiple expansion. 
Since at $\omega$ for which one of the these eigenvalues is zero, 
$\A(\omega, \delta)$ is not invertible, computing such roots $\omega$ yields then the desired subwavelength resonant frequencies.

The generalised capacitance matrix can be used to obtain accurate numerical approximations with a significant reduction in computational power. Provided that $\delta$ is sufficiently small, \Cref{thm:res} and \Cref{res_impart} give an approximation of the resonant frequencies that is sufficiently accurate for many purposes. In \Cref{fig:frequencies}, we show the resonant frequencies of a system of ten spherical resonators computed using both the full multipole method and using the approximation with the capacitance matrix. The values derived from the eigenvalues of the capacitance matrix give a good approximation and required just $0.02$ seconds of computation time, compared to the $41$ seconds required for the full multipole computations, on the same computer. If greater precision than that of the capacitance matrix approximation is required, then the values derived from the capacitance matrix can be used as initial values for root finding algorithms to reduce computational time. In particular, when Muller's method is used, in order to obtain all the roots, we can initialize it with the subwavelength resonant frequencies predicted by the capacitance matrix asymptotic analysis.

\section{Periodic systems}
Next, we will investigate the case when the resonators are repeated periodically as illustrated in \Cref{fig:periodic}, often referred to as a \emph{metamaterial}.  Conceptually, there are three different cases depending on the dimensions of periodicity of the lattice, $d_l$, and of the ambient space $d$:
\begin{itemize}
	\item $d-d_l = 0$, which we will refer to as a \emph{crystal}. In this fully periodic case (either $d=d_l =2$ or $d=d_l = 3$), the metamaterial has no boundary to the surrounding space;
	\item $d-d_l = 1$, which we will refer to as a \emph{screen}. In this case, the structure consists of a thin sheet of resonators;
	\item $d-d_l = 2$, which we will refer to as a \emph{chain}. There is one example of this case, namely $d=3,d_l = 1$.
\end{itemize}
We assume that $D$, as defined in \Cref{sec:finite}, is repeated in a periodic lattice $\Lambda$. We let $P_l: \R^d \to \R^{d_l}$ be the projection onto the first $d_l$ coordinates, and $P_\perp: \R^d \to \R^{d-d_l}$ be the projection onto the last $d-d_l$ coordinates.

We let $l_1,\dots, l_{d_l}\in \R^d$ denote lattice vectors generating the lattice $\Lambda$, in other words such that
$$\Lambda := \left\{ m_1 l_1+\dots+m_{d_l} l_{d_l} ~|~ m_i \in \Z \right\}. $$
For simplicity, we assume that $P_\perp l_i = 0$, which means that the lattice is aligned with the first $d_l$ coordinate axes.

For a point $x\in \R^d$, we will throughout use the notation $x = (x_l,x_0)$, where $x_l \in \R^{d_l}$ is the vector along the first $d_l$ dimensions and $x_0 \in \R^{d-d_l}$.
Denote by $Y \subset \R^{d}$ a fundamental domain of the given lattice. Explicitly, we take
$$ Y := \left\{ c_1 l_1+\dots+c_{d_l} l_{d_l} ~|~ 0 \le c_1,\dots,c_{d_l} \le 1 \right\}. $$
The dual lattice of $\Lambda$, denoted $\Lambda^*$, is generated by $\alpha_1,\dots,\alpha_{d_l}$ satisfying $ \alpha_i\cdot l_j = 2\pi \delta_{ij}$ and $P_\perp \alpha_i = 0$,  for $i,j = 1,\dots,d_l.$ The \emph{Brillouin zone} $Y^*$ is defined as $Y^*:= \big(\R^{d_l}\times\{\mathbf{0}\}\big) / \Lambda^*$, where $\mathbf{0}$ is the zero-vector in $\R^{d-d_l}$. We let $Y_l:=Y \cap \{x_0=0\}$ and remark that $Y^*$ can be written as $Y^*=Y^*_l\times\{\mathbf{0}\}$, where  $Y^*_l$ has the topology of a torus in $d_l$ dimensions.

The periodically repeated $i$\textsuperscript{th} resonator $\Dc_i$ and the full periodic structure $\Dc$ are given, respectively, by
$$\Dc_i = \bigcup_{m\in \Lambda} D_i + m, \qquad \Dc = \bigcup_{i=1}^N \Dc_i.$$
With this notation in hand, we rewrite \eqref{eq:finite_scattering}  as follows:
\begin{equation} \label{eq:scattering_periodic}
	\left\{
	\begin{array} {ll}
		\ds \Delta {u}+ k^2 {u}  = 0 & \text{in } \R^d \setminus \Dc, \\[0.3em]
		\ds \Delta {u}+ k_i^2 {u}  = 0 & \text{in } \Dc_i, \ i=1,\dots,N, \\
		\nm
		\ds  {u}|_{+} -{u}|_{-}  = 0  & \text{on } \partial \Dc, \\
		\nm
		\ds  \delta_i \frac{\partial {u}}{\partial \nu} \bigg|_{+} - \frac{\partial {u}}{\partial \nu} \bigg|_{-} = 0 & \text{on } \partial \Dc_i, \ i=1,\dots,N, \\
		\nm
		\ds u(x_l,x_0) & \text{satisfies the outgoing radiation condition as }  |x_0| \rightarrow \infty.
	\end{array}
	\right.
\end{equation}
Notice that while the set of material inclusions is compact in \Cref{sec:finite}, it is not the case in the infinite periodic setting. Consequently, the radiation condition in \eqref{eq:scattering_periodic} is no longer the Sommerfeld radiation condition. Moreover, it depends on the dimensionality (see, for instance, \cite{gerard1989resonance,ammari2002mathematical,ammari2018mathematical}).
Note also that we only impose a radiation condition for $x$ away from the periodic structure and, in the fully-periodic case, we do not assume any radiation condition. As we shall see, since there is no radiation condition in the first $d_l$ coordinate dimensions, the spectrum $\sigma$ of \eqref{eq:scattering_periodic} is in general continuous. In order to effectively study this equation, we will use the \emph{Floquet-Bloch} theory, which is outlined below.

\begin{figure}
	\begin{center}
		\includegraphics[width=\linewidth]{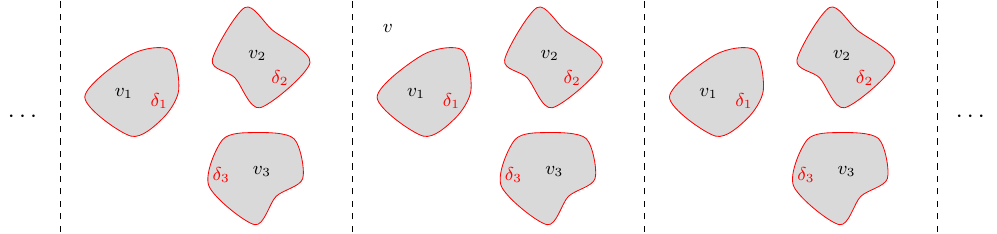}
	\end{center}
	\caption{A periodic array of material inclusions. Here, three material inclusions (resonators) are sketched with periodicity in one dimension. Each interior has a different wave speed $v_1$, $v_2$, $v_3$ and the surrounding medium has a wave speed $v$. The contrast between the $i$\textsuperscript{th} resonator and the background is given by $\delta_i$, where a small value of $\delta_i$ describes a large contrast.} \label{fig:periodic}
\end{figure}

\subsection{Floquet-Bloch theory}
A function $f(x)\in L^2(\R^d)$ is said to be $\alpha$-quasiperiodic, with quasiperiodicity $\alpha\in Y^*$, if $e^{-\i \alpha\cdot x}f(x)$ is $\Lambda$-periodic.  Given a function $f\in L^2(\R^d)$, the Floquet transform of $f$ is defined as
\begin{equation}\label{eq:floquet}
	\F[f](x,\alpha) := \sum_{m\in \Lambda} f(x-m) e^{\i \alpha \cdot m}, \quad x,\alpha\in \R^d.
\end{equation}
$\F[f]$ is always $\alpha$-quasiperiodic in $x$ and periodic in $\alpha$. The Floquet transform is an invertible map $\F:L^2(\R^d) \rightarrow L^2(Y\times Y^*)$, with inverse given by (see, for instance, \cite{ammari2018mathematical, kuchment1993floquet})
\begin{equation*}
	\F^{-1}[g](x) = \frac{1}{|Y_l^*|}\int_{Y^*} g(x,\alpha) \dx \alpha, \quad x\in \R^d,
\end{equation*}
where $g(x,\alpha)$ is extended quasiperiodically for $x$ outside of the unit cell $Y$.

If we apply the Floquet transform to \eqref{eq:scattering_periodic} we obtain, where $u^\alpha(x) := \F[u](x,\alpha)$,
\begin{equation} \label{eq:scattering_quasi}
	\left\{
	\begin{array} {ll}
		\ds \Delta {u^\alpha}+ k^2 {u^\alpha}  = 0 & \text{in } \R^d \setminus \Dc, \\[0.3em]
		\ds \Delta {u^\alpha}+ k_i^2 {u^\alpha}  = 0 & \text{in } \Dc_i, \ i=1,\dots,N, \\
		\nm
		\ds  {u^\alpha}|_{+} -{u^\alpha}|_{-}  = 0  & \text{on } \partial \Dc, \\
		\nm
		\ds  \delta_i \frac{\partial {u^\alpha}}{\partial \nu} \bigg|_{+} - \frac{\partial {u^\alpha}}{\partial \nu} \bigg|_{-} = 0 & \text{on } \partial \Dc_i, \ i=1,\dots,N,\\
		\nm
		\ds u^\alpha(x_l,x_0) & \text{is $\alpha$-quasiperiodic in }  x_l, \\
		\nm
		\ds u^\alpha(x_l,x_0) & \text{satisfies $\alpha$-quasiperiodic radiation condition as }  |x_0| \rightarrow \infty.
	\end{array}
	\right.
\end{equation}
The $\alpha$-quasiperiodic radiation condition depends on the dimensionality, and we refer (for example) to \cite{ammari2018mathematical,bao1995mathematical,bonnet1994guided,botten2013electromagnetic} for its explicit form. The spectrum $\sigma$ of the original problem \eqref{eq:scattering_periodic} is parametrised by the spectra $\sigma(\alpha), \ \alpha \in Y^*,$ of the problem \eqref{eq:scattering_quasi}, which in turn are known to consist of discrete values $\omega = \omega_i^\alpha$:
$$\sigma = \bigcup_{\alpha\in Y^*} \sigma(\alpha), \quad \sigma(\alpha) = \bigcup_{i=1}^\infty \omega_i^\alpha.$$
\begin{defn}[Band function]
	The resonant frequencies $\omega_i^\alpha$, seen as functions of $\alpha$, are called band functions. The collection of band functions is called the band structure.
\end{defn}
\begin{defn}[Band gap]
	A band gap of $\Dc$ is a connected component of $\CC\setminus \sigma$. If the spectrum $\sigma$ is real, we define a band gap of $\Dc$ as a connected component of $\R \setminus \sigma$, which consists of intervals in $\R$.
\end{defn}
As in \Cref{sec:finite}, we will focus on the subwavelength part of the spectrum, which are the resonant frequencies $\omega_i^\alpha$ which tend to $0$ as $\delta\to 0$. The results reviewed in this section are valid for $d\in\{2,3\}$ and $0<d_l\leq d.$

In the case $k \neq |\alpha + q|$ for all $q\in \Lambda^*$, we can define the quasiperiodic Green's function $G^{\alpha,k}(x)$ as the Floquet transform of $G^k(x)$ in the first $d_l$ coordinate dimensions, \ie{},
\begin{equation}\label{eq:xrep}
	G^{\alpha,k}(x) := \sum_{m \in \Lambda} G^k(x-m)e^{\i \alpha \cdot m}.
\end{equation}
Here $G^k$ is the usual Helmholtz Green's function, defined in \eqref{eq:G} for $d=2$ or $d=3$. The series in \eqref{eq:xrep} converges uniformly for $x$ and $y$ in compact sets of $\R^d$, $x\neq y$,  and $k \neq |\alpha + q|$ for all $q\in \Lambda^*$. Shortly speaking, using this Green's function we can define analogous quantities and get analogous results as in the finite case. The quasiperiodic single layer potential  $\mathcal{S}_D^{\alpha,k}$ is then defined as
$$\mathcal{S}_D^{\alpha,k}[\varphi](x) := \int_{\partial D} G^{\alpha,k} (x-y) \varphi(y) \dx\sigma(y),\quad x\in \mathbb{R}^d.$$
\begin{lemma} \label{lem:inv}
	The quasiperiodic single layer potential
 	$\S_D^{\alpha,k} : L^2(\p D) \rightarrow H^1(\p D)$ is invertible if $k$ is small enough and $k \neq  |\alpha + q|$ for all $q \in \Lambda^*$.
\end{lemma}
The quasiperiodic single layer potential satisfies many conceptually similar properties as the ``regular'' single layer potential (which is a consequence of the fact that the singularity of corresponding Green's functions are the same, \textit{i.e.}, $G^{\alpha,k} - G^k$ is a smooth function of $x$ around the origin). For example, the quasiperiodic single layer potential satisfies the jump relations on $\p D$:
\begin{equation} \label{eq:jump1_quasi}
	\S_D^{\alpha,k}[\varphi]\big|_+ = \S_D^{\alpha,k}[\varphi]\big|_- \quad \mbox{on}~ \p D,
\end{equation}
and
\begin{equation} \label{eq:jump2_quasi}
	\frac{\p}{\p\nu} \S_D^{\alpha,k}[\varphi] \Big|_{\pm} = \left( \pm \frac{1}{2} I +\Ka{\alpha}{k}\right)[\varphi]\quad \mbox{on}~ \p D,
\end{equation}
where $\Ka{\alpha}{k}$ is the quasiperiodic Neumann--Poincar\'e operator, given by
$$ \Ka{\alpha}{k}[\varphi](x):= \int_{\p D} \frac{\p}{\p\nu_x} G^{\alpha,k}(x-y) \varphi(y) \dx\sigma(y).$$

Above, we assumed that $k \neq |\alpha + q|$ for all $q\in \Lambda^*$. When $k$ is small and lies in the subwavelength regime, this condition separates into two cases:
\begin{itemize}
	\item $k < \inf_{q\in \Lambda^*} |\alpha + q|$. Waves in this regime are exponentially decaying away from the structure. Such waves, which vanish in the far-field, are known as \emph{evanescent} waves;
	\item $|\alpha| < k < \inf_{q\in \Lambda^*\setminus \{0\}} |\alpha + q|$. Waves in this regime are propagating far away from the structure, and this regime is known as the \emph{first radiation continuum}.
\end{itemize}
If we briefly assume that $v_i$ and $\delta_i$ are real we can interpret the two regimes as follows. When $k < \inf_{q\in \Lambda^*} |\alpha + q|$, the problem \eqref{eq:scattering_quasi} can be viewed as the spectral problem for a self-adjoint operator, and the resonant frequencies are real. When $|\alpha| < k < \inf_{q\in \Lambda^*\setminus \{0\}} |\alpha + q|$, due to the radiation condition this equation no longer corresponds to a self-adjoint operator. Therefore, the resonators attain a small but non-zero imaginary part corresponding to the coupling with the far-field. The transition between these two regimes occurs when $k = |\alpha|$ for some $q$ (known as a \emph{Rayleigh-Wood anomaly}), which are the points where the spectrum becomes real and the modes become localised to the structure.

In the two regimes mentioned above, we have the following integral representation	 (analogously to \Cref{lem:BIE_finite}).
\begin{lemma} \label{lem:BIE_periodic}
	Let $d\in \{2,3\}$ and $0<d_l\leq d$. Assume that $k \neq  |\alpha + q|$ for all $q \in \Lambda^*$. Then the Helmholtz problem \eqref{eq:scattering_periodic} is equivalent to finding $\psi,\phi\in L^2(\D)$ such that
	\begin{equation}
		\A^\alpha(\omega,\delta)
		\begin{pmatrix} \psi \\ \phi \end{pmatrix}
		= \begin{pmatrix} 0 \\ 0 \end{pmatrix},
	\end{equation}
	where the operator $\A^\alpha(\omega,\delta):L^2(\D)\times L^2(\D)\to H^1(\D)\times L^2(\D)$ is defined as
	\begin{equation}\label{aomegadeltalpha}
		\A^\alpha(\omega,\delta)= \begin{pmatrix}
			\widetilde{\S}_D^\omega & -\S_D^{\alpha,k}\\
			-\frac{1}{2}I+\widetilde{\K}_D^{\omega,*} & -\widetilde\delta \left(\frac{1}{2}I+\Ka{\alpha}{k}\right)
		\end{pmatrix},
	\end{equation}
	with $\widetilde{\delta}$, $\widetilde{\S}_D^\omega$ and $\widetilde{\K}_D^{\omega,*}$ defined in \eqref{eq:dtilde}, \eqref{eq:Stilde} and \eqref{eq:Ktilde}, respectively.
\end{lemma}

\begin{figure}
	\begin{center}
		\includegraphics[width=0.6\linewidth]{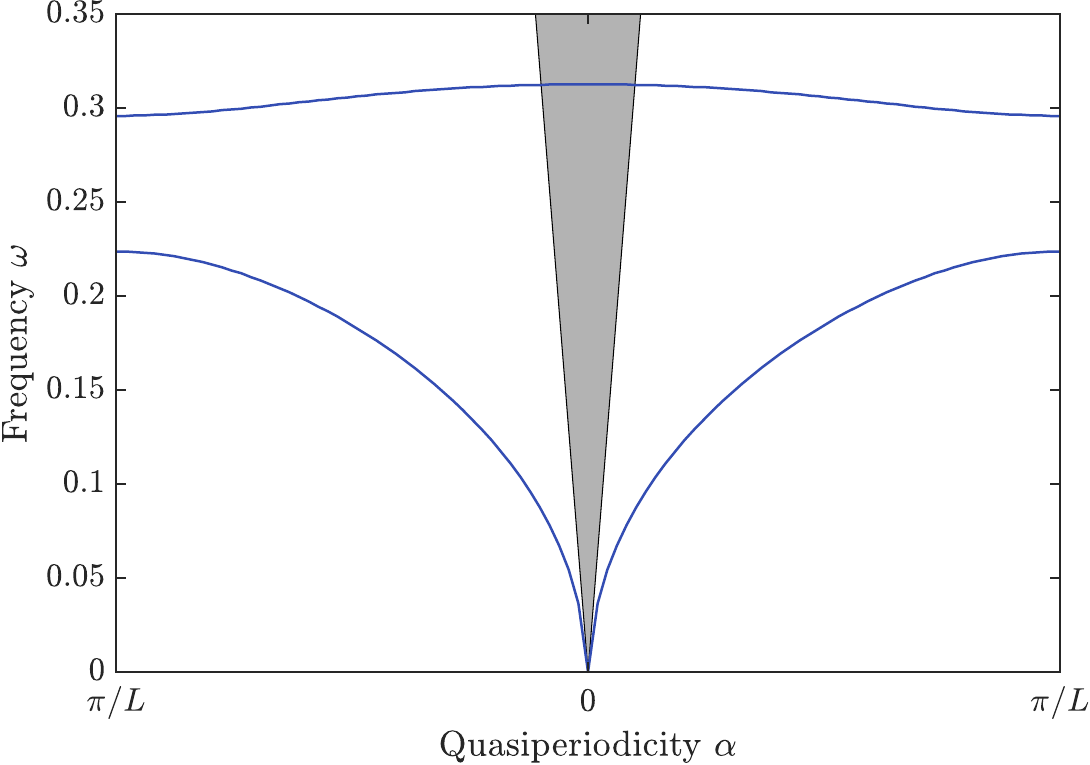}
	\end{center}
	\caption{Example of the subwavelength band structure of a resonator array with two resonators in the unit cell. The shaded region is the first radiation continuum, defined by $|\alpha| < \omega/v < \inf_{q\in \Lambda^*\setminus \{0\} }|\alpha+q|$, while the unshaded region correspond to evanescent modes. Here we see an example of band gap opening: there is an uncovered interval between the first and the second bands. Moreover, between the subwavelength bands and the higher (non-subwavelength) bands there will always be a band gap.} \label{fig:band5000}
\end{figure}

\subsection{Evanescent-mode resonances}
If we assume $|\alpha| > c > 0$ for some constant $c$ independent of $\omega$ and $\delta$, the quasiperiodic Helmholtz single layer potential is well approximated by the corresponding Laplace single layer potential in the sense that $\S_D^{\alpha,\omega}$ converges to $ \S_D^{\alpha,0}$ in the operator norm as $\omega \rightarrow 0$. Because of this result, which holds even in the case $d=2$ and  $0<d_l \leq 2$ (see, for instance, \cite[Section 2.12]{ammari2018mathematical}), we can use analogous methods as those outlined in \Cref{sec:finite}, and obtain similar results.
\begin{lemma} \label{lem:Aa00properties}
	Let $d\in \{2,3\}$ and $0<d_l\leq d$. Assume $|\alpha| > c > 0$ for some constant $c$ independent of $\omega$ and $\delta$, and consider a fundamental cell containing $N$ subwavelength resonators $D=D_1\cup\dots D_N$ in $Y$. Then, it holds that $\ker(-\frac{1}{2}I+\Ka{\alpha}{0})=\mathrm{span}\{\psi_1^\alpha,\psi_2^\alpha,\dots,\psi_N^\alpha\}$, where $\psi_i^\alpha:=(\mathcal{S}_D^{\alpha,0})^{-1}[\chi_{\partial D_i}]$.
\end{lemma}

\begin{thm}
Let $d\in \{2,3\}$ and $0<d_l\leq d$. Assume $|\alpha| > c > 0$ for some constant $c$ independent of $\omega$ and $\delta$, and consider a fundamental cell containing $N$ subwavelength resonators $D=D_1\cup\dots D_N$ in $Y$. For sufficiently small $\delta>0$, there exist $N$ subwavelength resonant frequencies $\omega_1^\alpha(\delta),\dots,\omega_N^\alpha(\delta)$ with non-negative real parts.
\end{thm}
\begin{remark}
	The above theorem describes only the \emph{subwavelength} part of the band structure. For small enough $\delta$, there will be a band gap between the first $N$ bands (which are in the subwavelength regime) and the higher bands (which are not close to $0$ for small $\delta$; see \Cref{fig:band5000}).
\end{remark}

\begin{defn}[Quasiperiodic capacitance matrix] \label{defn:QCM}
	Assume $\alpha \neq 0$. For a system of $N\in\N$ resonators $D_1,\dots,D_N$ in $Y$ we define the quasiperiodic capacitance matrix $C^\alpha=(C^\alpha_{ij})\in\CC^{N\times N}$ to be the square matrix given by
	\begin{equation*}
		C^\alpha_{ij}=-\int_{\D_i} (\S_D^{\alpha,0})^{-1}[\chi_{\p D_j}] \de\sigma,\quad i,j=1,\dots,N.
	\end{equation*}
\end{defn}
\begin{lemma}
	The quasiperiodic capacitance matrix is a Hermitian matrix.
\end{lemma}
\begin{defn}[Generalised quasiperiodic capacitance matrix] \label{defn:GQCM}
	Assume $\alpha \neq 0$. For a system of $N\in\N$ resonators $D_1,\dots,D_N$ in $Y$ we can define the generalised quasiperiodic capacitance matrix, denoted by $\C^\alpha=(\C^\alpha_{ij})\in\CC^{N\times N}$, to be the square matrix given by
	\begin{equation*}
		\C^\alpha_{ij}=\frac{\delta_i v_i^2}{|D_i|} C^\alpha_{ij}, \quad i,j=1,\dots,N.
	\end{equation*}
\end{defn}

The next result characterises the first $N$ resonances of the periodic structure, and shows that they are in the subwavelength regime.
\begin{thm} \label{thm:res_periodic}
	Let $d\in \{2,3\}$ and $0<d_l\leq d$. Consider a system of $N$ subwavelength resonators in $Y$, and assume $|\alpha| > c > 0$ for some constant $c$ independent of $\omega$ and $\delta$. As $\delta\to0$, the $N$ subwavelength resonant frequencies satisfy the asymptotic formula
	\begin{equation*}
		\omega_n^\alpha = \sqrt{\lambda_n^\alpha}+O(\delta^{3/2}), \quad n=1,\dots,N,
	\end{equation*}
	where $\{\lambda_n^\alpha: n=1,\dots,N\}$ are the eigenvalues of the generalised quasiperiodic capacitance matrix $\mathcal{C}^\alpha\in\mathbb{C}^{N\times N}$, which satisfy $\lambda_n^\alpha=O(\delta)$ as $\delta\to0$.
\end{thm}
\begin{remark}
	The error term 	$O(\delta^{3/2})$ has higher order compared to the error term $O(\delta)$ in \Cref{thm:res}. This is a consequence of the fact that the $O(\omega)$-term in the expansion of $\S_D^{\alpha,\omega}$ vanishes.
\end{remark}

\begin{cor} \label{prop:eigenvector_quasi}
	Let $d\in \{2,3\}$, $0<d_l\leq d$ and assume $|\alpha| > c > 0$ for some constant $c$ independent of $\omega$ and $\delta$. Let $\textbf{v}_n^\alpha$ be the eigenvector of $\mathcal{C}^\alpha$ associated to the eigenvalue $\lambda_n^\alpha$. Then the resonant mode $u_n^\alpha$ associated to the resonant frequency $\omega_n^\alpha$ is given, as $\delta\to0$, by
	\begin{equation*}
		u_n^\alpha(x)=\begin{cases}
			\textbf{v}_n^\alpha\cdot\textbf{S}_D^{\alpha,k}(x)+O(\delta^{1/2}), \quad x\in\mathbb{R}^d\setminus \overline{\Dc}, \\
			\textbf{v}_n^\alpha\cdot\textbf{S}_D^{\alpha,k_i}(x)+O(\delta^{1/2}), \quad x\in \Dc_i,
		\end{cases}
	\end{equation*}
	where $\textbf{S}_D^{\alpha,k}:\mathbb{R}^d\to\mathbb{C}^N$ is the vector-valued function given by
	\begin{equation*}
		\textbf{S}_D^{\alpha,k} (x)=\begin{pmatrix}
			\mathcal{S}_D^{\alpha,k}[\psi_1^\alpha](x) \\[-0.4em]
			\vdots \\[-0.3em]
			\mathcal{S}_D^{\alpha,k}[\psi_N^\alpha](x)
		\end{pmatrix},  \quad x\in\mathbb{R}^d\setminus \partial \Dc,
	\end{equation*}
	with $\psi_i^\alpha:=(\mathcal{S}_D^{\alpha,0})^{-1}[\chi_{\partial D_i}]$.
\end{cor}
\Cref{fig:efuncS} shows the resonant mode of a square crystal in two dimensions ($d=d_l=2$) for $\alpha$ close to the corner of the Brillouin zone. Here, we can observe the two-scale behaviour associated to subwavelength metamaterials: the resonant modes are oscillating on the small scale, with amplitudes which satisfy large-scale oscillations (for more details on this, we refer to \cite{ammari2019bloch}).
\begin{figure}[tbh]
	\centering
	\vspace{-10pt}
	\begin{subfigure}[b]{0.43\linewidth}
		\centering
		\includegraphics[scale=0.4]{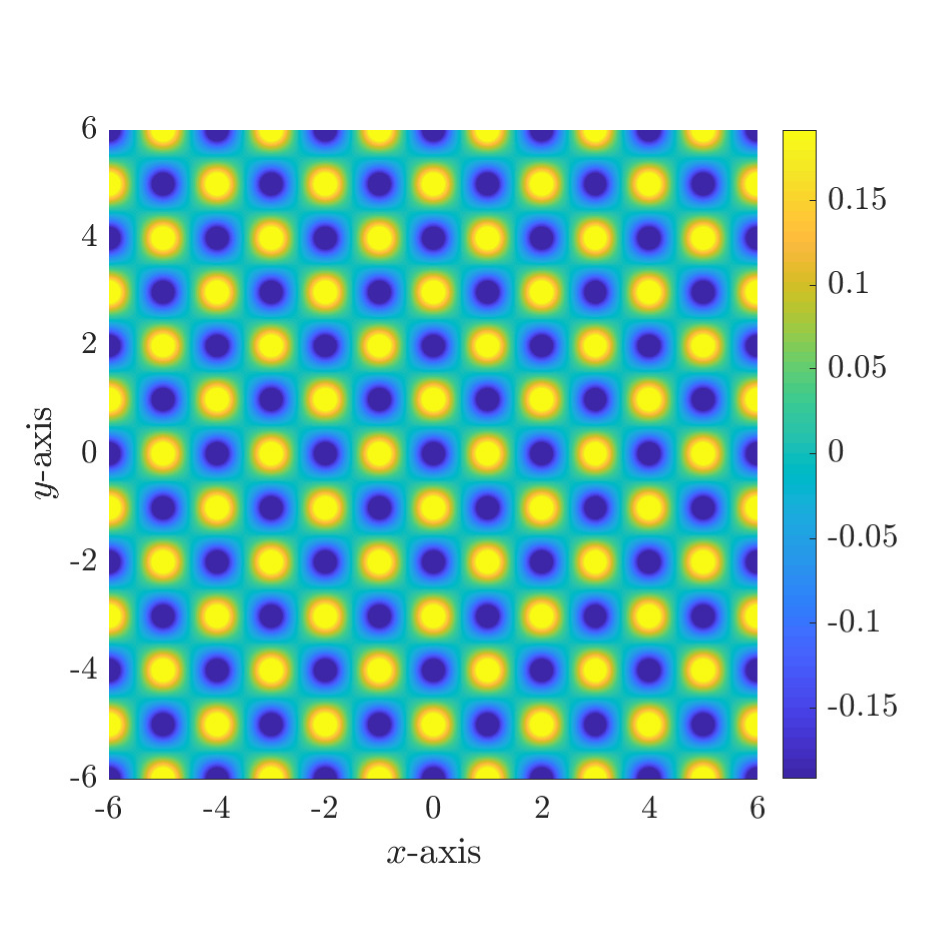}
		\vspace{-15pt}
		\caption{Small-scale behaviour of the resonant mode.}
		\label{fig:S}
	\end{subfigure}
	\hspace{3pt}
	\begin{subfigure}[b]{0.52\linewidth}
		\includegraphics[scale=0.45]{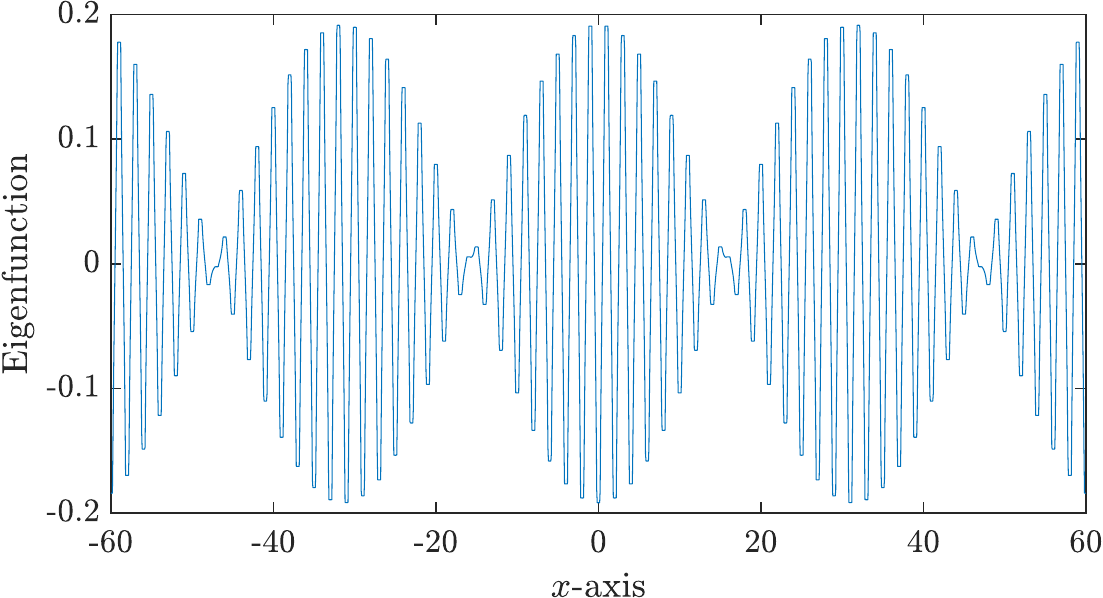}
		\caption{One-dimensional plot over many unit cells.}
		\label{fig:s1D}
	\end{subfigure}
	\caption{Plot of the resonant mode of a square crystal for $\alpha$ close to the corner of the Brillouin zone. We see that the resonant modes have a distinct two-scale behaviour: rapidly oscillating on the small scale, and a large scale envelope which satisfies a homogenized equation.} \label{fig:efuncS}
\end{figure}
\begin{remark}
 As shown in \Cref{sec:twodim}, the discrete formulations for approximating the subwavelength resonant frequencies for finite systems in two and three dimensions are slightly different. In contrast with the finite system setting,  for infinitely periodic systems,
 exactly the same capacitance matrix formulation  holds in both  two and three dimensions. This is due to the fact that 
 the single layer potential associated with the Helmholtz equation in free space has a  logarithmic singularity in $\omega$ near zero (see \eqref{eq:expansion2d}) while the $\alpha$-quasiperiodic single layer potentials are smooth functions of $\omega$ around zero when $\alpha \neq 0$. 
\end{remark}

\subsection{Resonances in the first radiation continuum}
Here we study the regime $|\alpha| < k < \inf_{q\in \Lambda^*\setminus \{0\} }|\alpha+q|$. Since the analysis depends on the dimensionalities, we will exemplify it in the case  of a metascreen, $d-d_l = 1$.

Recall that $k = \omega/v$. In the current regime, we must have $\alpha \to 0$ as $\omega \to 0$. Therefore, we assume that
$$\alpha =  \omega \alpha_0 \in Y^*,$$
for some $\alpha_0$, independent of $\omega$ and such that $|\alpha_0| < 1/v$. In scattering problems, this limit corresponds to incident waves with a fixed direction of incidence (specified by $\alpha_0$) and a frequency $\omega$ in the subwavelength regime.

\subsubsection{Green's function and capacitance matrix formulation}\label{sec:periodicGreens}
In the current setting, the quasiperiodic Green's function admits the spectral representation
\begin{equation} \label{eq:specrep} G^{\alpha,k}(x) = \frac{e^{\i\alpha\cdot x}e^{\i k_0|x_0|}}{2\i k_0|Y_l|} - \sum_{q\in \Lambda^*\setminus \{0\}}\frac{e^{\i(\alpha+q)\cdot x}e^{-\sqrt{|\alpha+q|^2 - k^2}|x_0|}}{2|Y_l|\sqrt{|\alpha+q|^2 - k^2} },\end{equation}
where $x=(x_l,x_0)$ and $k_0 = \sqrt{k^2 - |\alpha|^2}$. The series in \eqref{eq:specrep} converges uniformly for $x$ in compact sets of $\R^d$, $x\neq 0$,  and $|\alpha| < k < \inf_{q\in \Lambda^*\setminus \{0\} }|\alpha+q|$. In the case when $k = \alpha = 0$, we define the \emph{periodic} Green's function $G^{0,0}$  as
\begin{equation} \label{eq:specrep_0}
	G^{0,0}(x) = \frac{|x_0|}{2|Y_l|} - \sum_{q\in \Lambda^*\setminus \{0\}}\frac{e^{\i q \cdot x}e^{-|q||x_0|}}{2|Y_l||q|}.
\end{equation}
Here, $G^{0,0}$ is periodic in the $x_l$-variable. When $\omega \rightarrow 0$, we then have the asymptotic expansion
\begin{equation} \label{eq:Gexp}
	G^{\omega\alpha_0,k}(x) = \frac{1}{2\i k_0 |Y_l|} + G^{0,0}(x) + \frac{\alpha \cdot x}{2k_0|Y_l|} + O(\omega).
\end{equation}
In particular, the Green's function has a singularity when $\omega\to 0$. In fact, this will make the analysis conceptually similar to the case in \Cref{sec:twodim}. We define the operator $\widehat{\S}_D^{\alpha,k}: L^2(\p D) \rightarrow H^1(\p D)$ as
\begin{equation}
	\label{eq:Shat}\widehat{\S}_D^{\alpha,k}[\varphi](x) = {\S}_D^{0,0}[\varphi](x) - \frac{\i - \alpha\cdot x}{2 k_0 |Y_l|} \int_{\p D}\varphi \dx \sigma - \int_{\p D}\frac{\alpha \cdot y}{2k_0|Y_l|} \varphi(y)\dx \sigma (y).
\end{equation}
We then have the asymptotic expansion  $\S_D^{\omega\alpha_0,\omega} = \widehat{\S}_D^{\omega\alpha_0,\omega} + O(\omega)$, with respect to the operator norm, as $\omega \rightarrow 0$ \cite{ammari2021bound}.
\begin{lemma} \label{lem:ker}
	$\S_D^{0,0}$ is invertible from the mean-zero space $L^2_0(\p D)$ onto its image.
\end{lemma}

\begin{lemma} \label{lem:holo}
	For any $\alpha_0 \in Y^*$ with $0<|\alpha_0| < 1$,  $\left({\S}_D^{\omega\alpha_0,\omega}\right)^{-1}$ is a holomorphic operator-valued function of $\omega$ in a neighbourhood of $\omega = 0$.
\end{lemma}
Crucially, \Cref{lem:holo} shows that inverse $({\S}_D^{\omega\alpha_0,\omega})^{-1}$ does not have the $\omega^{-1}$-singularity around $\omega=0$. This will allow us to define capacitance coefficients in the current setting. Since $({\S}_D^{\omega\alpha_0,\omega})^{-1}$ is a holomorphic function of $\omega$ we have $\left({\S}_D^{\omega\alpha_0,\omega}\right)^{-1} = \S_0^{\alpha_0} +  O(\omega)$ as $\omega \to 0$, for some operator $\S_0^{\alpha_0}$ which is independent of $\omega$.

\begin{defn}[Periodic capacitance matrix] \label{defn:PCM}
	For $\alpha_0$ with $|\alpha_0| < 1$,  and for a system of $N\in\N$ resonators $D_1,\dots,D_N$ in $Y$ we can define the periodic capacitance matrix $C^0=(C^0_{ij})\in\R^{N\times N}$ to be the square matrix given by
	\begin{equation*}
		C_{ij}^{0}=-\int_{\D_i} \S_0^{\alpha_0}[\chi_{\p D_i}] \de\sigma,\quad i,j=1,\dots,N.
	\end{equation*}
\end{defn}
Similarly to \Cref{lem:capacitance_altdefn}, the periodic capacitance matrix can alternatively be written
\begin{equation}
    \label{eqc0}
C_{ij}^0 = \int_{Y\setminus D} \nabla V_i^0\cdot  \nabla V_j^0 \de x,
\end{equation}
where $V_i^0 = \S_D^{0,0}[\psi_i^0]$ and $\psi_i^0 =\S_0^{\alpha_0}[\chi_{\p D_i}]$. Notice that $V_i^0$ is the unique solution to following problem:
$$
\left\{ \begin{array}{ll}
		\ds \Delta V_i^0  = 0 & \text{in } Y \setminus D, \\
		\nm
		\ds V_i^0   = \delta_{ij}  & \text{on } \partial D_j, \\
		\nm
		\ds V_i^0(x_l,x_0) & \text{is periodic in }  x_l, \\
		\nm
		\ds V_i^0(x_l,x_0) = O(|x_0|^{-1}) & \text{as }  |x_0| \rightarrow \infty.
\end{array}
\right.$$

Because of \eqref{eqc0}, $C^0$ share many of the properties of the capacitance matrices in previous settings.
\begin{lemma} \label{lem:cap0}
	The periodic capacitance matrix $C^{0}$ is a real, symmetric, positive semi-definite matrix with one vanishing eigenvalue. Moreover, $C^0$ is independent of $\alpha_0$.
\end{lemma}
The name \emph{periodic} comes from the fact that an $\alpha$-quasiperiodic function is, in the case $\alpha=0$, a periodic function. Since $\alpha = \omega \alpha_0$ and $\omega$ is small, we are working close to the periodic case. In fact, since $C^0$ is independent of $\alpha_0$, all relations involving $C^0$ are equal to the periodic case $\alpha_0 = 0$.
\begin{defn}[Generalised periodic capacitance matrix] \label{defn:GPCM}
	For a system of $N\in\N$ resonators $D_1,\dots,D_N$ in $Y$ we can define the generalised periodic capacitance matrix, denoted by $\C^0=(\C^0_{ij})\in\CC^{N\times N}$, to be the square matrix given by
	\begin{equation*}
		\C^0_{ij}=\frac{\delta_i v_i^2}{|D_i|} C^0_{ij}, \quad i,j=1,\dots,N.
	\end{equation*}
\end{defn}
Since $\S_D^{\omega\alpha_0,\omega}$ has a $\omega^{-1}$-singularity as $\omega\to 0$, we cannot apply the Gohberg-Sigal theory to study the characteristic value perturbation of $\A^\alpha(\omega,\delta)$, defined in \Cref{lem:BIE_periodic}, as $\delta \rightarrow 0$. Instead, we rephrase the system in the following result \cite{ammari2020exceptional}.
\begin{lemma} \label{lem:BIE_hat}
	Let $k= \omega/v$ and assume that $\alpha = \omega\alpha_0$ for some $\alpha_0$ independent of $\omega$ and $\delta$ such that $|\alpha_0|< 1/v$. Then the Helmholtz problem \eqref{eq:scattering_periodic} is equivalent to finding $\eta\in H^1(\D)$ such that
	\begin{equation}
		\widehat{\A}^\alpha(\omega,\delta)\eta
		= 0,
	\end{equation}
	where the operator $\widehat{\A}^\alpha(\omega,\delta):H^1(\D)\to L^2(\D)$ is defined as
	\begin{equation*}
		\widehat{\A}^\alpha(\omega,\delta)=
			\left(-\frac{1}{2}I+\widetilde{\K}_D^{\omega,*}\right)\left(\widetilde{\S}_D^\omega\right)^{-1} -\widetilde\delta \left(\frac{1}{2}I+\Ka{\alpha}{k}\right)\left(\S_D^{\omega\alpha_0,k}\right)^{-1}.
	\end{equation*}
\end{lemma}
Now we can apply the functional analytic approach, outlined in \Cref{sec:approach}, to derive the following result on the subwavelength resonant frequencies in the current setting \cite{ammari2020exceptional, ammari2021bound}.
\begin{thm} \label{thm:res0}
	Let $d \in\{2,3\}$ and $d-d_l = 1$. Let $k= \omega/v$ and assume that $\alpha = \omega\alpha_0$ for some $\alpha_0$ independent of $\omega$ and $\delta$ such that $|\alpha_0|< 1/v$.  As $\delta\to0$, there are $N$ subwavelength resonant frequencies which satisfy the asymptotic formula
	\begin{equation*}
		\omega_n^\alpha = \sqrt{\lambda_n^0}+O(\delta), \quad n=1,\dots,N,
	\end{equation*}
	where $\{\lambda_n^0: n=1,\dots,N\}$ are the eigenvalues of the generalised periodic capacitance matrix $\mathcal{C}^0\in\mathbb{C}^{N\times N}$, which satisfy $\lambda_n^0=O(\delta)$ as $\delta\to0$.
\end{thm}

\subsubsection{Higher-order approximations}
In the case $\delta_iv_i^2\in \R$ for all $i$, the leading-order parts of the resonances, given in \Cref{thm:res0}, are real. It is also interesting to compute corresponding imaginary parts, which will specify the bandwidth of the resonant behaviour. We have the asymptotic expansion of $\S_D^{\omega\alpha_0,\omega} = S_0^{\alpha_0} + \omega\S_{-1}^{\alpha_0} + O(\omega^2)$. Then we can define the ``higher-order'' matrix $C^{1,\alpha_0}$ as
$$C_{ij}^{1,\alpha_0} = -\int_{\D_i}\S_{-1}^{\alpha_0}[\chi_{\p D_j}]\de \sigma$$
and corresponding generalised matrix $\C^{1,\alpha_0}$ as
$$\C_{ij}^{1,\alpha_0} = \frac{\delta_iv_i^2}{|D_i|} C_{ij}^{1,\alpha_0}.$$
Therefore, we have the following result, which gives the resonances $\omega_i^\alpha$ to a higher order \cite{ammari2020exceptional}.
\begin{thm} \label{thm:higher_periodic}
	Let $d \in\{2,3\}$ and $d-d_l = 1$. Let $k= \omega/v$ and assume that $\alpha = \omega\alpha_0$ for some $\alpha_0$ independent of $\omega$ and $\delta$ such that $|\alpha_0|< 1/v$.  As $\delta\to0$, the $N$ subwavelength resonant frequencies satisfy $\omega_n^\alpha = \widehat{\omega}_n^\alpha + O(\delta^{3/2})$ where $\widehat{\omega}_n^\alpha$, for $n=1,\dots,N$, are the roots $\omega = \widehat{\omega}_n^\alpha$ of the equation
	$$\det\left(\C^{0} + \omega\C^{1,\alpha_0} -\omega^2I\right) = 0.$$
\end{thm}
We define
$\mathbf{c}_n \in \mathbb{C}^d$ by
$$\mathbf{c}_n= \int_{\p D}y\psi_n^0(y)\de \sigma(y), \quad n = 1,\dots,N.$$
While the capacitance coefficients can be thought of as total charge (or ``mass''), $\mathbf{c}_i$ is the centre of mass (up to rescaling). Briefly put, we can compute $\C^{1,\alpha_0}$ in terms of these coefficients, which allows us to compute explicit expressions of $\widehat{\omega}_n^\alpha$.

\subsubsection{Modal decomposition}
The ``higher-order'' term $\S_{-1}^{\alpha_0}$ in the expansion of the single layer potential enters the expression for the resonant modes in this case. From the arguments used to derive \Cref{thm:higher_periodic} we have the following result on the resonant modes \cite{ammari2020exceptional}.
\begin{thm} \label{prop:eigenvector_periodic}
	Let $d \in\{2,3\}$ and $d-d_l = 1$. Let $k= \omega/v$ and assume that $\alpha = \omega\alpha_0$ for some $\alpha_0$ independent of $\omega$ and $\delta$ such that $|\alpha_0|< 1/v$. Let $\textbf{v}_n^0$ be the eigenvector of $\mathcal{C}^0$ associated to the eigenvalue $\lambda_n^0$. Then the resonant mode $u_n^\alpha$ associated to the resonant frequency $\omega_n^\alpha$ is given, as $\delta\to0$, by
	\begin{equation*}
		u_n^\alpha(x)=\begin{cases}
			\textbf{v}_n^0\cdot\textbf{S}_D^{\alpha,k}(x)+O(\delta^{1/2}), \quad x\in\mathbb{R}^d\setminus \overline{\Dc}, \\
			\textbf{v}_n^0\cdot\textbf{S}_D^{\alpha,k_i}(x)+O(\delta^{1/2}), \quad x\in \Dc_i,
		\end{cases}
	\end{equation*}
	where $\textbf{S}_D^{\alpha,k}:\mathbb{R}^d\to\mathbb{C}^N$ is the vector-valued function given by
	\begin{equation*}
		\textbf{S}_D^{\alpha,k}(x)=\begin{pmatrix}
			\mathcal{S}_D^{\alpha,k}[\psi_1^0 + k \psi_1^{1,\alpha_0}](x) \\[-0.4em]
			\vdots \\[-0.3em]
			\mathcal{S}_D^{\alpha,k}[\psi_N^0 + k \psi_N^{1,\alpha_0}](x)
		\end{pmatrix},  \quad x\in\mathbb{R}^d\setminus \partial \Dc,
	\end{equation*}
	with $\psi_i^0:=\S_0^{\alpha_0}[\chi_{\partial D_i}]$ and $\psi_i^{1,\alpha_0}:=\S_{-1}^{\alpha_0}[\chi_{\p D_i}]$.
\end{thm}
With these resonant modes at hand, we can easily solve the associated scattering problem
\begin{equation} \label{eq:scattering_p}
	\left\{
	\begin{array} {ll}
		\ds \Delta {u^\alpha}+ k^2 {u^\alpha}  = 0 & \text{in } \R^d \setminus \Dc, \\[0.3em]
		\ds \Delta {u^\alpha}+ k_i^2 {u^\alpha}  = 0 & \text{in } \Dc_i, \ i=1,\dots,N, \\
		\nm
		\ds  {u^\alpha}|_{+} -{u^\alpha}|_{-}  = 0  & \text{on } \partial \Dc, \\
		\nm
		\ds  \delta_i \frac{\partial {u^\alpha}}{\partial \nu} \bigg|_{+} - \frac{\partial {u^\alpha}}{\partial \nu} \bigg|_{-} = 0 & \text{on } \partial \Dc_i, \ i=1,\dots,N,\\
		\nm
		\ds u^\alpha(x_l,x_0) & \text{is $\alpha$-quasiperiodic in }  x_l, \\
		\nm
		\ds u^\alpha-\uin & \text{satisfies $\alpha$-quasiperiodic radiation condition as }  |x_0| \rightarrow \infty.
	\end{array}
	\right.
\end{equation}
Here, we assume that $\uin$ is a plane wave with some wave vector $\k$; $\uin(x) = e^{\i \k\cdot x}$. Furthermore, $\alpha$ is now specified in terms of $\k$ as $\alpha = P_l\k$, where (as before) $P_l$ is the projection to the first $d_l$ coordinates. Then, for small $\omega$ we assume that $\k = \omega \w$ for some fixed $\w$.

\begin{thm} \label{lem:modal_per}
	Let $d \in\{2,3\}$ and $d-d_l = 1$. Assume that $\uin = e^{\i \k\cdot x}$, where $\k = \omega \w$. Also, assume that $|\omega -\omega_i^0| > K\sqrt{\delta}$ for $i=1,\dots,N$, for some constant $K>0$. Then, as $\delta \to 0$, the solution to the scattering problem \eqref{eq:scattering_p} can be written, uniformly for $x$ in compact subsets of $\mathbb{R}^d$, as
	\begin{equation}\label{eq:scattered}
		(u^\alpha - \uin)(x) = \sum_{n=1}^N a_n u_n^\alpha(x) - \S_D^{\alpha,k}\left(\S_D^{\alpha,k}\right)^{-1}[\uin](x) + O(\sqrt{\delta}),
	\end{equation}
where $V$ is the matrix of eigenvectors of $\C^0$ and the coefficients $a_n=a_n(\omega)$ satisfy 
\begin{equation*}
V\begin{pmatrix}
	\omega^2-(\omega_1^0)^2 & & \\ & \ddots & \\ & & \omega^2-(\omega_N^0)^2
\end{pmatrix}
\begin{pmatrix} a_1 \\ \vdots \\ a_N \end{pmatrix}
=
\begin{pmatrix}
	\frac{\delta_1v_1^2}{|D_1|} \int_{\D_1} \left({\S}_D^{\alpha,k}\right)^{-1}[\uin]\de\sigma \\
	\vdots \\
	\frac{\delta_Nv_N^2}{|D_N|} \int_{\D_N}\left({\S}_D^{\alpha,k}\right)^{-1}[\uin]\de\sigma
\end{pmatrix}+O(\delta^{3/2}).
\end{equation*}
\end{thm}
\subsubsection{Connection to quasiperiodic capacitance matrix}
In the previous sections, we have used the notation $C^0$ for the periodic capacitance matrix, and $C^\alpha$ for the quasiperiodic capacitance matrix, defined for $\alpha\in Y^*\setminus\{0\}$. This choice of suggestive notation is deliberate, and in this section we will show that these capacitance matrices can be combined into a continuous function of $\alpha \in Y^*$, in the sense that $C^0 = \lim_{\alpha \to 0} C^\alpha$.

For $\alpha\in Y^*$, $\alpha\neq 0$, the Green's function $G^{\alpha,0}$ satisfies
\begin{equation}
	G^{\alpha,0}(x) = - \frac{1}{2|Y_l|}\sum_{q\in \Lambda^*}\frac{e^{\i(\alpha+q)\cdot x}e^{-|\alpha+q||x_0|}}{|\alpha+q|},
\end{equation}
so that for small $\alpha$ we have
\begin{equation}
	G^{\alpha,0}(x) = - \frac{1}{2 |\alpha| |Y_l|} - \frac{\i\alpha\cdot x}{2 |\alpha| |Y_l|} + G^{0,0}(x) + O(|\alpha|).
\end{equation}
In this setting we define the operator $\widehat{\S}_D^{\alpha,0}: L^2(\p D) \rightarrow H^1(\p D)$ as
\begin{equation}	\label{eq:Shat0}
	\widehat{\S}_D^{\alpha,0}[\varphi](x) = {\S}_D^{0,0}[\varphi](x) - \frac{1 + \i\alpha\cdot x}{2 |\alpha| |Y_l|} \int_{\p D}\varphi \dx \sigma + \int_{\p D}\frac{\i\alpha \cdot y}{2|\alpha||Y_l|} \varphi(y)\dx \sigma (y).
\end{equation}
The structure of this operator is entirely analogous to $\widehat{\S}_D^{\alpha,k}$ defined in \Cref{sec:periodicGreens} (where, conceptually, the limit $\omega\to 0$ now corresponds to $|\alpha|\to0$). Therefore, we can apply the same method  to obtain the following result \cite{ammari2021bound}.
\begin{lemma}
	The periodic capacitance matrix $C^0$ and the quasiperiodic capacitance matrix $C^\alpha$, for $\alpha\in Y^*\setminus\{0\}$, satisfy
	$$C^0 = \lim_{\alpha\to 0}C^\alpha.$$
\end{lemma}

\section{Applications of the generalised capacitance matrix}
The generalised capacitance matrix can be applied to explain a variety of interesting physical phenomena. We particularly want to study extraordinary macroscopic properties (such as exotic effective parameters, rainbow trapping or unidirectional scattering) as well as robust localisation at subwavelength scales.

\subsection{Double-negative materials}\label{sec:doubleneg}

The field of metamaterials more or less began with the realisation that micro-structured media could be designed to have effectively negative material parameters. This behaviour can also be replicated in the setting studied in this work, as reported in \cite{ammari2019double,fepponhomogenization}. Consider a large number $N$ of identical resonator pairs, given by
\begin{equation} \label{eq:negmaterial}
	D^N=\bigcup_{1\leq j\leq N} \left(z_j^N+s R_{d_j^N} D\right),
\end{equation}
where $D=D_1\cup D_2$ is some fixed pair of resonators, $0<s\ll1$ is some characteristic size, $z_j^N\in\mathbb{R}^3$ are the positions of the resonator pairs and $R_{d_j^N}$ are rotations in $\mathbb{R}^3$ which orient the pair in the direction of the unit vector $d_j^N$. We will assume that $sN=\Lambda$ for some fixed $\Lambda>0$ and that there is some bounded domain $\Omega$ such that $\{z_j^N:1\leq j\leq N\}\subset\Omega$ for any $N\geq1$. We also want that the resonators are regularly distributed in the sense that there exists some number $\nu$ such that
\begin{equation} \label{reg1}
	\min_{i\neq j} |z_i^N-z_j^N|\geq \nu N^{-1/3}, \quad\text{for any }N\geq1.
\end{equation}
Additionally, we want that there exists some positive function $V\in C^1(\overline{\Omega})$ and a matrix-valued function $B\in C^1(\overline{\Omega})$ such that there are constants $C_1$ and $C_2$ which satisfy
\begin{equation}\label{reg2}
	\max_{1\leq j\leq N}\left|\frac{1}{N} \sum_{i\neq j} G^k(z_i^N-z_j^N)f(z_i^N)-\int_\Omega G^k(y-z_j^N) V(y)f(y) \de y\right|\leq C_1 N^{-\alpha/3} \|f\|_{C^{0,\alpha}(\Omega)},
\end{equation}
for all $f\in C^{0,\alpha}(\Omega)$ with $0<\alpha\leq1$ and all $N\geq1$, where $G^k$ is defined by \eqref{eq:G}.
Similarly, we have
\begin{equation} \label{reg3}
	\max_{1\leq j\leq N}\left|\frac{1}{N} \sum_{i\neq j} \left(f(z_i^N)\cdot d_i^N\right)\left( d_i^N\cdot \nabla G^k(z_i^N-z_j^N)\right)-\int_\Omega f(y) B G^k(y-z_j^N) \de y\right|\leq C_2 N^{-\alpha/3} \|f\|_{\left(C^{0,\alpha}(\Omega)\right)^3},
\end{equation}
for all $f\in \left(C^{0,\alpha}(\Omega)\right)^3$ with $0<\alpha\leq1$ and all $N\geq1$.
\begin{remark}
	The regularity assumptions \eqref{reg2} and \eqref{reg3} are challenging to comprehend in general. In the case that the positions $\{z_j^N:1\leq j\leq N\}$ are uniformly distributed in $\Omega$, then $V(x)$ is constant. Likewise, if the orientations are such that the average of $d_j^N(d_j^N)^\top$ in any neighbourhood of $\Omega$ converges to the identity, then $B(x)$ is equal to some positive function times the identity matrix.
\end{remark}

We must make some additional assumptions to achieve the desired double-negative behaviour. We assume that the material parameters are such that $v_1=v_2\in\mathbb{R}$ and $\delta_1=\delta_2=\mu^2s^2$ for some number $\mu>0$. This means that the generalized capacitance coefficients associated to a single resonator pair do not depend on $s$. As such, the leading-order term in the expansion of the subwavelength resonant frequencies, as $s\to0$, is fixed. Additionally, we will assume that the pair of resonators is symmetric in the sense that %\todo{think about examples}
\begin{equation} \label{symmetrydimer}
	P:=\int_{\D} y_1(\psi_1-\psi_2)\de\sigma(y)>0,
	\quad \text{while} \quad
	\int_{\D} y_2(\psi_1-\psi_2)\de\sigma(y)=\int_{\D} y_3(\psi_1-\psi_2)\de\sigma(y)=0,
\end{equation}
for $\psi_1$ and $\psi_2$ as defined in \Cref{lem:A00properties}. In this case, we also have that the subwavelength resonant frequencies for the resonator pair are given, as $s\to0$, by
\begin{align}
	\omega_1&= \sqrt{\lambda_1} -\i \tau_1 \mu^2s+O(s^2), \label{eq:dimer2a} \\  
	\omega_2&= \sqrt{\lambda_2} + \mu^3 \eta_1 s^2 -\i \mu^4 \eta_2 s^3 + O(s^4),  \label{eq:dimer2}
\end{align}
for real numbers $\tau_1$, $\eta_1$ and $\eta_2$.

We want the incident frequency to be close to the second resonant frequency, which corresponds to the dipole resonant mode. In particular, we assume that there is some $a< \mu^3\eta_1$ such that
\begin{equation} \label{eq:freq_assump}
	\omega=\sqrt{\lambda_2}+as^2.
\end{equation}
Under these assumptions, we can derive an effective medium theory for the system as $N\to\infty$. Let $u^N$ be the field scattered by an array of $N$ resonator pairs of the form \eqref{eq:negmaterial}. Suppose there is some macroscopic field $u\in C^{1,\alpha}(\overline{\Omega})$ which is such that $u^N$ converges to $u$ in $C^{1,\alpha}(\overline{\Omega})$ as $N\to\infty$. Then, we can show that $u$ must satisfy the equation \cite{ammari2019double}
\begin{equation} \label{eq:effectivemedium}
	\nabla\cdot M_1(x) \nabla u(x) +M_2(x)u(x)=0 \quad \text{in } \mathbb{R}^3,
\end{equation}
where
\begin{equation*}
	M_1=\begin{cases}
		I & \text{in } \mathbb{R}^3\setminus\Omega,\\
		I-\Lambda g^1 B & \text{in } \Omega,
	\end{cases}
\qquad\text{and}\qquad
	M_2=\begin{cases}
		k^2 & \text{in } \mathbb{R}^3\setminus\Omega,\\
		k^2-\Lambda g^0 V & \text{in } \Omega,
	\end{cases}
\end{equation*}
for positive constants $g^0$ and $g^1$ given by
\begin{equation*}
	g^0=\frac{2(C_{11}+C_{12})}{1-\lambda_1/\lambda_2}
	\qquad \text{and} \qquad
	g^1=\frac{\mu^2v_1^2}{2|D|\lambda_2(\mu^3\eta_1-a)}P^2.
\end{equation*}

\Cref{eq:effectivemedium} shows the desired double-negative behaviour. Since we assumed that $a< \mu^3\eta_1$, it holds that $g^1>0$. Thus, if $B(x)$ is a positive definite matrix for all $x\in\Omega$ and $\Lambda$ is chosen to be sufficiently large, then $M_1(x)$ will be a negative definite matrix for all $x\in\Omega$. Likewise, since $g^0>0$ and $V>0$ by assumption, we can choose $\Lambda$ to be sufficiently large that $M_2(x)<0$ for all $x\in\Omega$. When considering the effective parameters $M_1$ and $M_2$ as a function of the incident frequency $\omega \, (=k v)$, one can see that there is a range of frequencies in which both the effective material parameters are negative.
We refer to \cite{ammari2019double} for the details.
\begin{remark}
    For a system $D$ of two identical resonators $D_1$ and $D_2$ which are symmetric with respect to the origin and separated from each other, the symmetry assumption \eqref{symmetrydimer} is satisfied \cite{ammari2019double}. The asymptotic formulas \eqref{eq:dimer2a} and \eqref{eq:dimer2} show that the system $D$ features two slightly different subwavelength resonant frequencies. Such frequencies are called hybridised resonant frequencies. They correspond to two fundamentally different resonant modes. The first resonant mode is a monopole mode while the second mode is a dipole mode. The resonant frequency associated with the dipole mode is usually referred to as the anti-resonant frequency. For an appropriate volume fraction, when the incident frequency is close to the anti-resonant frequency, the collection of the dipole modes contribute to the effective parameter in the high-order term in the differential operator ($M_1$) while the collection of the monopole modes to the zero-th order term ($M_2$).   
\end{remark}

\subsection{Frequency separation with graded arrays}

By introducing a gradient to the properties of a metamaterial, it is possible to cause different frequencies to be localised at different positions in the structure. This frequency separation is useful for a variety of applications since it means an incoming signal can be separated into its different frequency components. This phenomenon is often known as \emph{rainbow trapping} and has been observed in a variety of electromagnetic \cite{tsakmakidis2007trapped, jang2011plasmonic}, acoustic \cite{ammari2019cochlea, zhu2013acoustic, jimenez2017rainbow} and other metamaterials \cite{bennetts2018graded}. Approximating the system using the generalised capacitance matrix, we can study the extent to which a system of \emph{coupled} subwavelength resonators is able to perform frequency separation in a controllable manner.

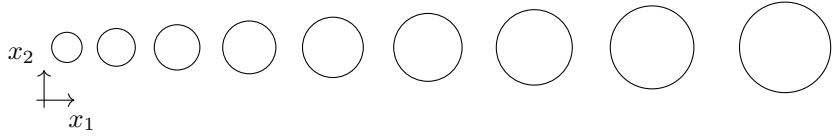
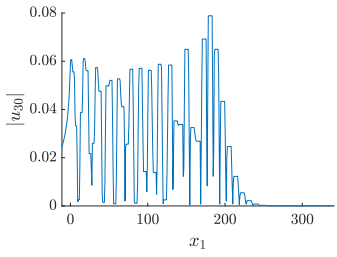
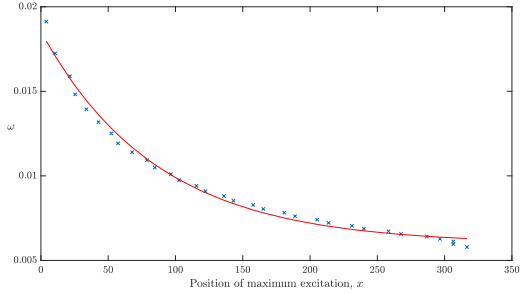
\begin{figure}
	\begin{center}
		\begin{subfigure}{0.75\linewidth}
			\begin{center}
				\begin{tikzpicture}
					\draw (0.3,0) circle (0.2);
					\draw (0.95,0) circle (0.25);
					\draw (1.75,0) circle (0.3);
					\draw (2.7,0) circle (0.35);
					\draw (3.8,0) circle (0.4);
					\draw (5.05,0) circle (0.45);
					\draw (6.45,0) circle (0.5);
					\draw (8,0) circle (0.55);
					\draw (9.75,0) circle (0.6);

					\draw [->] (-0.1,-0.7) -- (0.4,-0.7);
					\draw [->] (0,-0.8) -- (0,-0.3);

					\node at (0.5,-1) {$x_1$};
					\node at (-0.3,-0.1) {$x_2$};
				\end{tikzpicture}
			\end{center}
			\caption{A graded array of resonators. High frequencies will give a maximum response towards the left of the array while lower frequencies will be detected further to the right.} \label{fig:graded}
		\end{subfigure}

		\begin{subfigure}{0.3\linewidth}
			\includegraphics[width=\linewidth]{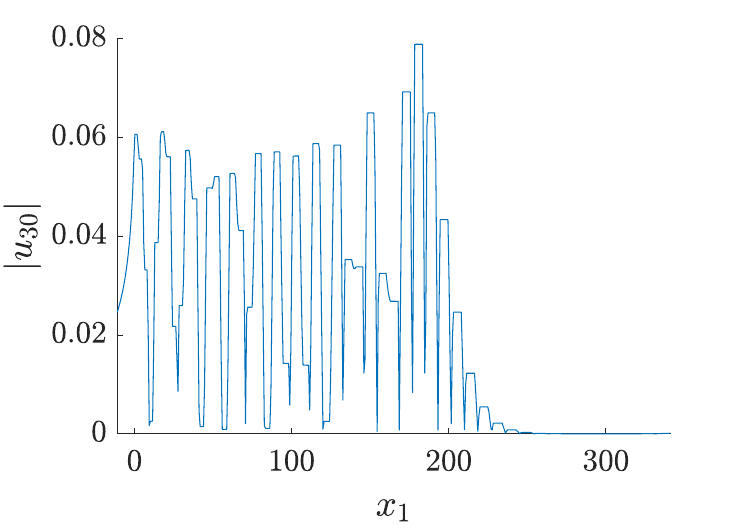}
			\caption{The modulus of the 30\textsuperscript{th} eigenmode for a graded array of 50 resonators, as an example. There is a clear position of maximum amplitude.} \label{fig:tonotopic}
		\end{subfigure}
		\hspace{0.8cm}
		\begin{subfigure}{0.5\linewidth}
			\includegraphics[width=\linewidth]{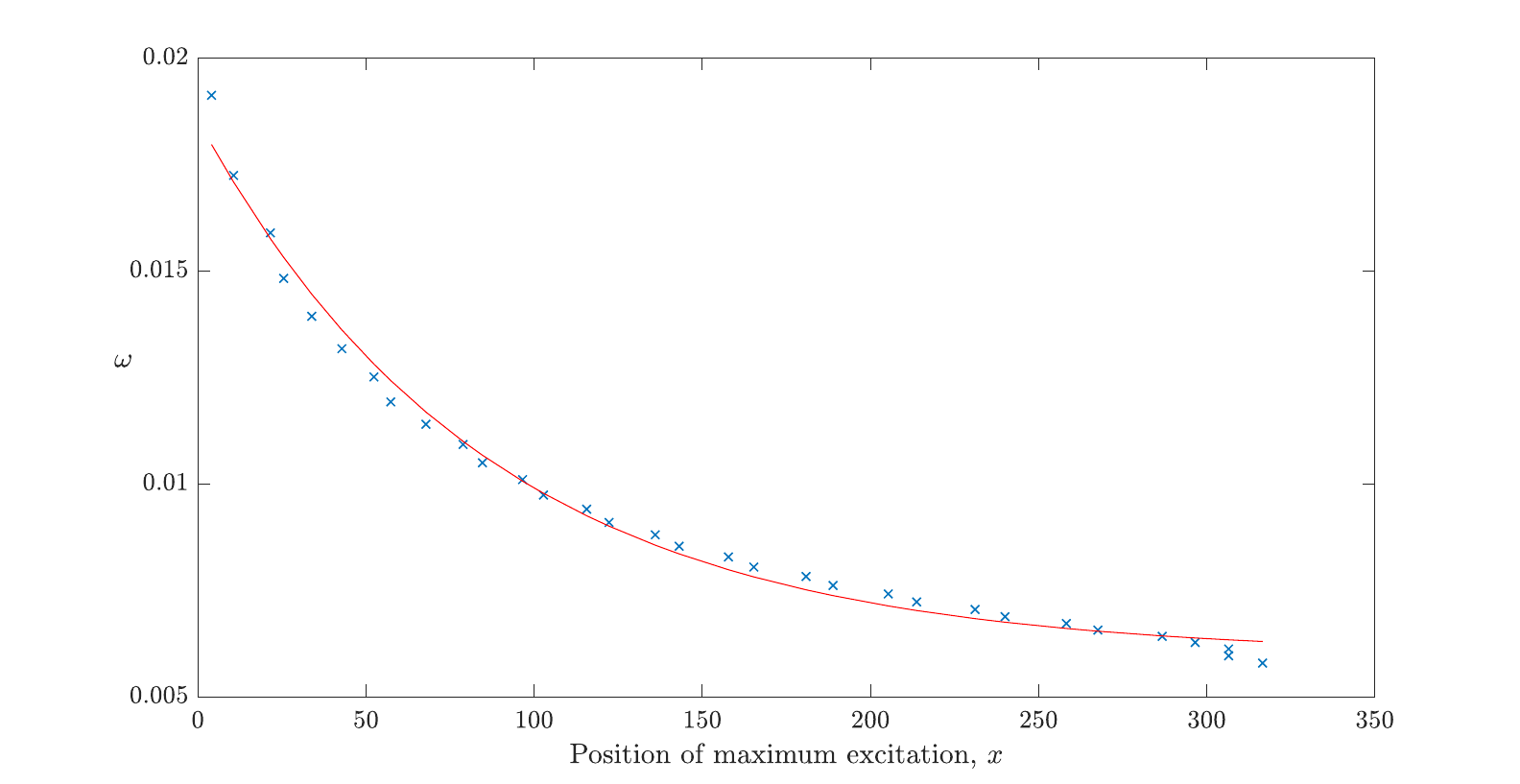}
			\caption{The relationship between the frequency (real part) and the position of maximum amplitude. The crosses are the resonant frequencies of an array of 50 resonators and the line is a relationship of the form observed in the human cochlea.} \label{fig:modes}
		\end{subfigure}
	\end{center}
	\caption{The frequency separation of a graded array of resonators can be chosen to mimic the response of the cochlea. A graded array of 50 resonators is simulated here and its response is chosen to match the relationship that exists in the human cochlea.} \label{fig:cochlea}
\end{figure}

This behaviour is very similar to the action of the cochlea. Devices based on these principles have been used to design biomimemtic hearing devices \cite{ammari2020mimicking, davies2022robustness, rupin2019mimicking, karlos2020cochlea}. For example, a graded resonator array is shown in \Cref{fig:cochlea} to replicate the frequency separation of the cochlea. These structures are useful for building artificial hearing approaches as well as, conversely, learning about the function of human hearing. For example, understanding the details of cochlear amplification is a significant open question and is obscured by the challenges in experimenting on living biological organisms. By designing analogue artificial devices, we are able to test theories and reveal crucial insight into the fundamantal mechanisms that underpin human hearing \cite{ammari2020mimicking, rupin2019mimicking}. See \cite{ammari2023cochleabook} for a review.

\subsection{Sensitivity enhancement using high-order exceptional points} \label{sec:exceptional}
\begin{figure}
	\begin{center}
		\begin{subfigure}{0.5\linewidth}
			\includegraphics[width=\linewidth]{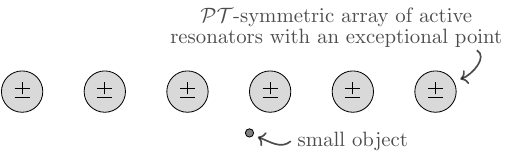}
			\caption{An active system of resonators can be designed such that it supports an exceptional point, meaning that it experiences an enhanced response to small perturbations such as the presence of a small particle.} \label{fig:smallobject}
		\end{subfigure}
	\hspace{0.8cm}
		\begin{subfigure}{0.2\linewidth}
			\includegraphics[width=\linewidth]{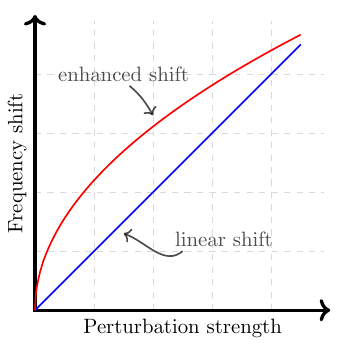}
			\caption{Exceptional sensors can experience enhanced eigenfrequency shifts.} \label{fig:enhancedresponse}
		\end{subfigure}
	\end{center}
	\caption{Resonator arrays with exceptional points can be used to design enhanced sensors.} \label{fig:sensorarray}
\end{figure}

We want to design ``enhanced'' sensors that are strongly influenced by small perturbations in their immediate surroundings. For example, we want to be able to sense the presence of a small object, as depicted in \Cref{fig:smallobject}. This small object might be a virus, for example \cite{vollmer2008single}. Typically, such an occurrence would cause a shift in the resonant frequencies that is proportional to the strength of the perturbation (\emph{i.e.} proportional to the size of the small object). The idea here, however, is to design an array for which this shift is enhanced, particularly for small perturbations \cite{wiersig2016sensors, hodaei2017enhanced}. We achieve this by adding sources of energy gain and loss to the system, represented by material parameters with non-zero imaginary parts. This means that \Cref{lem:eigenbasis} no longer holds and we can fine-tune the material parameters in order to create systems with coincident eigenvectors. The fundamental result guiding this work is the following lemma \cite{ammari2020high} (see, also, \cite{haiexceptional}).

\begin{lemma} \label{lem:enhanced}
	Let $m\in L^\infty(\mathbb{R}^3)$ (possibly complex) be such that $m \equiv 1$ outside a ball $B$. For $k \neq 0$ real, let the Green's function $G_m^k$ be defined as the solution to $\left(\Delta_x + m(x)k^2\right) G_m^k(x,y)  = \delta_y(x)$ in $\mathbb{R}^3$, subject to the Sommerfeld radiation condition.
 Consider  the Hilbert-Schmidt operator $T^k_B: L^2(B) \rightarrow L^2(B)$ defined by 
 $$T^k_B[f](x)= - \int_B (m(y) -1) G^k(x,y) f(y)\, \dx y, \quad x\in B,$$
 where $G^k$ is the free-space Helmholtz Green's function. Then, 
 $$  G_m^k(x,y) = G^k(x,y) + \Bigl(\frac{1}{k^2} I - T^k_B\Bigr)^{-1}\bigl[G^k(\cdot, y)\bigr](x,y), \quad x,y \in B.$$
 Moreover, let  $k^*\in\mathbb{C}$ be a resonant frequency, \ie{} $k^*$ is a characteristic value of $k \mapsto \frac{1}{k^2} I - T^k_B$, and 
 suppose that the system has an $N$\textsuperscript{th}-order singularity in the sense that $G_m^k$
	has the following {pole-pencil decomposition}:
	\begin{equation*} \label{eq:app_Gexcep}
	G_m^k(x,y) =  G^k(x,y) +\sum_{j=1}^N c_j \frac{\varphi_j(x)\varphi_j(y)}{{(k^2 -(k^*)^2)^j}}
	+ R(k,x,y),
	\end{equation*}
	in a neighbourhood of $k^*$, where $\{\varphi_j\}_{j \geq 1}$ form an orthonormal basis of $L^2(B)$ of  generalised eigenfunctions of $T^{k^*}_B$, $\{c_j\}_{1\leq j\leq N}$ are some constants,  and the remainder $R$ is a holomorphic function of $k$ that is smooth in $x,y \in B$. If a small material inclusion $\Omega$ is introduced to the system, then the new system has a resonant frequency $k^{\Omega}$ with the asymptotic behaviour
	\begin{equation*}
	k^{\Omega}=k^*+(\eta_z |\Omega|)^{1/N}+o(|\Omega|), \quad \text{as }|\Omega|\to0,
	\end{equation*}
	with $\eta_z$ is a constant that depends on the position and material parameters of the small particle $\Omega$, the background material parameters $m$ and $\varphi_N$.
\end{lemma}

One way to create $N$\textsuperscript{th}-order singularities, as are required by \Cref{lem:enhanced}, is to design structures with higher-order resonant singularities. That is, structures with \emph{exceptional points}, where eigenvalues and eigenvectors coincide. In this setting, we will search for \emph{asymptotic} exceptional points.
\begin{defn}
Consider the Helmholtz resonance problem \eqref{eq:finite_scattering} for the domain $D=D_1\cup\dots\cup D_N$, where $N\in\N$. A set of material parameter values is said to be an $N$\textsuperscript{th}-order asymptotic exceptional point with respect to $\delta$ if there exist $N$ resonant frequencies $\omega_1,\dots,\omega_N$ and associated eigenmodes $u_1,\dots,u_N$ such that for any $i,j\in\{1,\dots,N\}$
	\begin{equation*}
	\omega_i=\omega_j+O(\delta),\quad \text{as } \delta\to0,
	\end{equation*}
	and for any $i,j\in\{1,\dots,N\}$ there exists some $K\in\mathbb{C}$ such that
	\begin{equation*}
	u_i=K u_j +O(\delta),\quad \text{as } \delta\to0.
	\end{equation*}
\end{defn}

\begin{remark}
	The restriction to considering \emph{asymptotic} exceptional points in this work is not a weakness of the analytic approach but represents the behaviour of the system. The radiation condition means that the symmetry we impose on the resonators is not extended to the far field, meaning that we won’t have exact degeneracy at the exceptional points.
\end{remark}

Exceptional points are a consequence of balanced symmetries in the system, which cause the eigenvectors to align. So that the system already has some underlying symmetry, exceptional points are often sought in structures with \emph{parity--time symmetry}. We will assume that the problem is parity--time symmetric in the sense that each resonator $D_i$ can be uniquely associated to another resonator $D_j$ (possibly with $i=j$) such that
\begin{equation} \label{eq:PT}
D_ i = \P D_j \quad \text{and} \quad v_i^2\delta_i = \T(v_j^2\delta_j),
\end{equation}
where the parity operator $\mathcal{P}: \R^3 \rightarrow \R^3$ and the time-reversal operator $\T: \mathbb{C} \rightarrow \mathbb{C}$ are given, respectively, by
$$\mathcal{P}(x) = -x \quad\text{and}\quad \T(z) = \overline{z}.$$

For a system of two subwavelength resonators, we can find the eigenvalues and eigenvectors of the generalised capacitance matrix explicitly. Using these formulas, we can show that there exists an asymptotic exceptional point for certain parameter values \cite{ammari2020exceptional}.

\begin{thm} \label{thm:EP2}
A $\mathcal{PT}$-symmetric pair of subwavelength resonators 
$D_1$ and $D_2$, \textit{i.e.}, satisfying \eqref{eq:PT},
has an asymptotic exceptional point of order two with respect to $\delta$ in the subwavelength regime. In other words, there is a set of material parameters such that the eigenvalues and eigenvectors of the associated generalised capacitance matrix defined by \eqref{eq:GCM} coincide. In particular, if
\begin{equation*}
\mathrm{Im}(v_1^2\delta_1)=b^*:=\frac{\mathrm{Re}(v_1^2\delta_1)C_{12}}{\sqrt{C_{11}^2-C_{12}^2}},
\end{equation*}
then $\lambda_1=\lambda_2$ and $\textbf{v}_1=K\textbf{v}_2$ for some $K\in\mathbb{C}$, where $(\lambda_i,\textbf{v}_i)$, $i=1,2$, are the eigenpairs of the generalised capacitance matrix $\mathcal{C}$. Further to this,

\vspace{-0.6cm}

\begin{gather*}
\text{if  } \mathrm{Im}(v_1^2\delta_1)<b^* \text{  then  } \sqrt{\lambda_1},\sqrt{\lambda_2} \text{ are real valued and } \sqrt{\lambda_1}\neq\sqrt{\lambda_2}, \\
\text{if  } \mathrm{Im}(v_1^2\delta_1)>b^* \text{  then  } \sqrt{\lambda_1},\sqrt{\lambda_2}\text{ are purely imaginary and } \sqrt{\lambda_1}\neq\sqrt{\lambda_2}.
\end{gather*}
\end{thm}

The approximate nature of the asymptotic exceptional point predicted in \Cref{thm:EP2} is demonstrated by \Cref{fig:EP2}, where the subwavelength resonant frequencies of the full differential system are simulated directly using the multipole expansion method. We can see that there is a critical value of the gain and loss (the imaginary parts of the material parameters) such that the eigenvalues coincide at leading order. Below this critical value the leading-order parts of the resonant frequencies are real and they form a purely imaginary conjugate pair above this value.

\begin{figure}
	\centering
	\begin{subfigure}{0.4\linewidth}
		\centering
		\begin{tikzpicture}
		\draw (-1.5,0) circle (1);
		\draw (1.5,0) circle (1);
		\node at (-2.5,1){$D_1$};
		\node at (2.5,1){$D_2$};
		\node at (-1.5,0){\small $v_1^2\delta_1{=}a{+}\i b$};
		\node at (1.5,0){\small $v_2^2\delta_2{=}a{-}\i b$};
		\node at (0,-1.3){\color{white} .};
		\node at (0,0.7){\small $v$};
		\end{tikzpicture}
		\caption{A $\mathcal{PT}$-symmetric pair of spherical resonators.}
	\end{subfigure}
	\hspace{0.1cm}
	\begin{subfigure}{0.55\linewidth}
		\includegraphics[width=\linewidth]{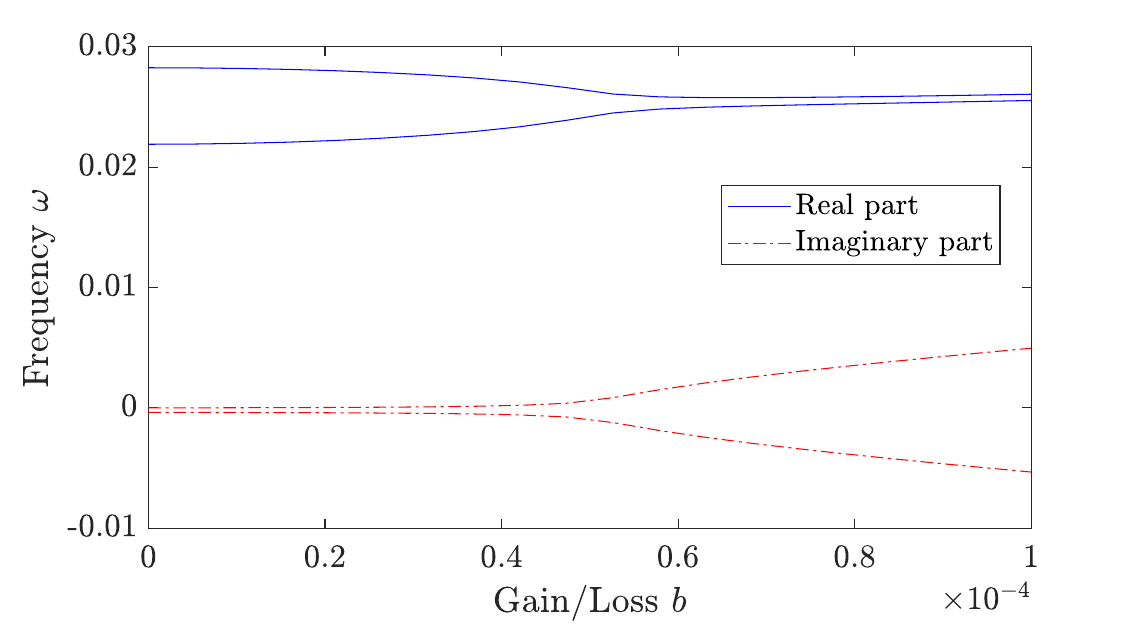}
		\caption{The two subwavelength resonant frequencies.}
	\end{subfigure}
	\caption{The two subwavelength resonant frequencies of a pair of $\mathcal{PT}$-symmetric resonators can be simulated directly from the full differential system using the multipole expansion method. An asymptotic exceptional point occurs at $b^*\approx0.5\times10^{-4}$, where the frequencies coincide at leading order. For $b$ smaller than $b^*$, the leading-order terms of the resonant frequencies are real, while for $b$ larger than $b^*$ they are purely imaginary and are the conjugate of one other, again to leading order.} \label{fig:EP2}
\end{figure}
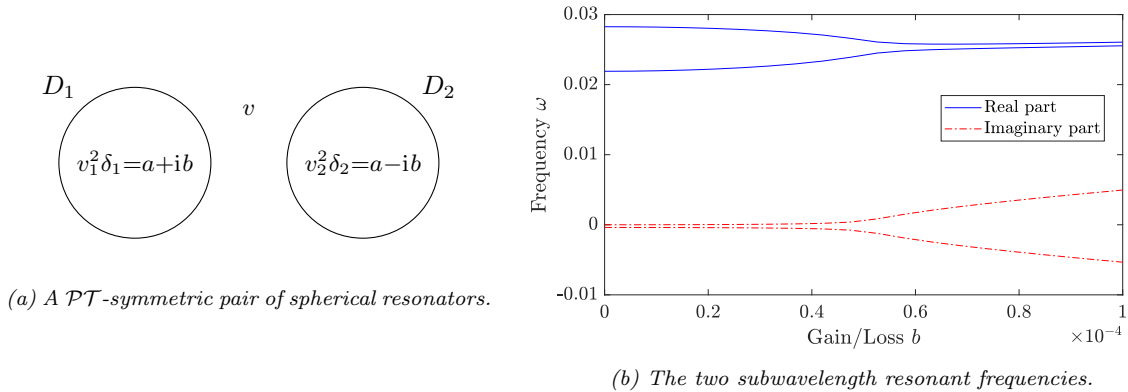

In order to produce higher-order subwavelength exceptional points, we need to study larger systems of resonators. So that the capacitance matrix is easier to work with in this case, we will make an additional assumption that the resonators are relatively far apart, meaning that we can use the dilute approximation given in \Cref{lem:dilutefinite}. This means we can efficiently analyse large systems of resonators, in terms of this leading-order approximation of the capacitance matrix \cite{ammari2020high}. For a system of three subwavelength resonators, we can show that there is one third-order subwavelength asymptotic exceptional point. This is described by \Cref{thm:EP3} and the subwavelength resonant frequencies are depicted in \Cref{fig:EP3}.

\begin{thm} \label{thm:EP3}
	A $\mathcal{PT}$-symmetric system $D$ of three dilute resonators has an asymptotic exceptional point of order $3$ with respect to $\epsilon$ and $\delta$ at the resonant frequency $\omega^*$, which is given as $\epsilon,\delta\to0$ by
	$$\omega^* = \sqrt{\frac{4\pi (3+\epsilon c_1)\mathrm{Re}(v_1^2\delta_1)}{3|D_1|}} + O(\delta + \delta^{1/2}\epsilon),$$
	where $c_1$ is the real root of the polynomial $c_1^3+\frac{27}{4}c_1-\frac{27}{8}=0$ (\textit{i.e.} $c_1\approx0.483\dots$).
\end{thm}

We can continue this process to study higher-order exceptional points in larger systems. We quickly find that as the dimension grows the number of exceptional point similarly grows. In a system of four subwavelength resonators we find that there are four asymptotic exceptional points. The imaginary parts of the material parameters on each resonator at each exceptional point are depicted in \Cref{fig:sols4}. We can see that each exceptional point corresponds to one of the four different combinations of relative magnitude and sign that is possible under the assumption of $\mathcal{PT}$ symmetry.

The symmetry exhibited by the fourth-order asymptotic exceptional points shown in \Cref{fig:sols4} can also be seen in higher-order exceptional points in larger structures. For example, in \Cref{fig:highEPs}, we take the fourth-order exceptional points from Figures~\ref{fig:sol4a}~and~\ref{fig:sol4d} and find asymptotic exceptional points of order 8 and 14 in dilute resonator arrays of the corresponding size. Exceptional points with the same qualitative distribution as in Figures~\ref{fig:EPsymm4}--\ref{fig:EPsymm14} were previously observed in Hamiltonian systems in \cite{zhang2020supersymm}. This analysis demonstrates the value of the generalised capacitance matrix, particularly under an assumption of diluteness, which gives a concise yet rigorous approximation of the behaviour of the differential system.

\begin{figure}
	\centering
	\begin{subfigure}{0.33\linewidth}
		\includegraphics[width=\linewidth]{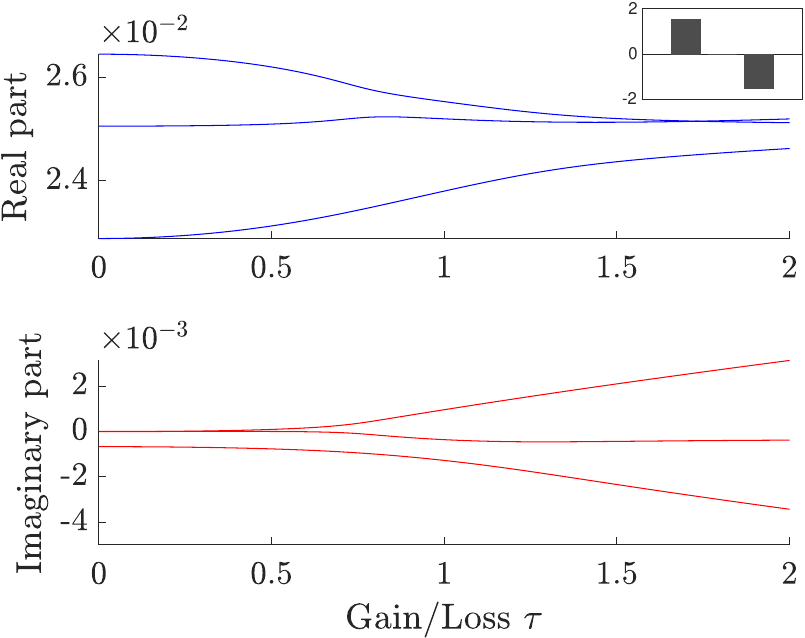}
		\caption{The subwavelength resonant frequencies of the full differential problem. \\}
	\end{subfigure}
	\hspace{1cm}
	\begin{subfigure}{0.33\linewidth}
		\includegraphics[width=\linewidth]{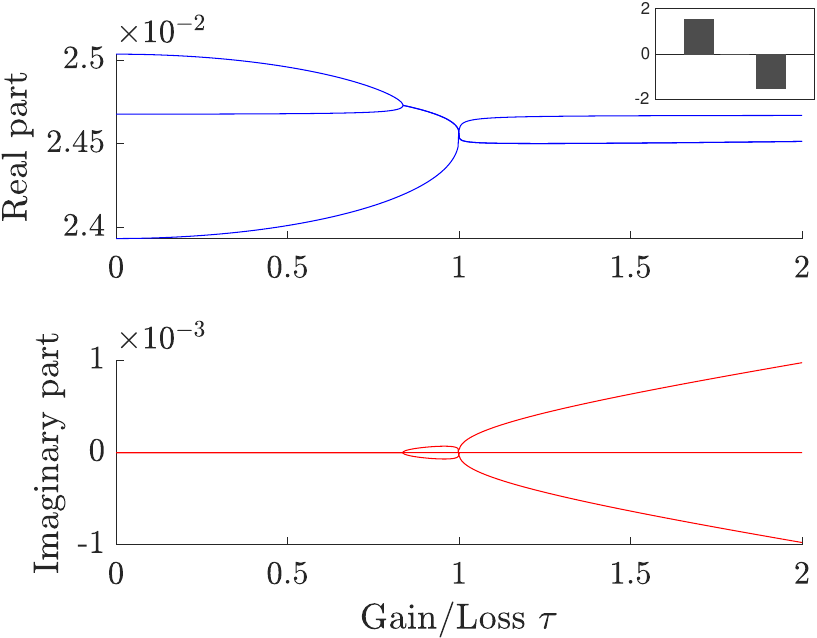}
		\caption{The approximate subwavelength resonant frequencies, under the dilute approximation.}
	\end{subfigure}
	\caption{A $\mathcal{PT}$-symmetric system of three subwavelength resonators supports an asymptotic exceptional point, where the eigenvalues (and corresponding eigenvectors) of the dilute capacitance matrix coincide. We can compare the resonant frequencies of the full differential problem in (a), computed using the multipole expansion method, and the approximate frequencies in (b), computed using the dilute approximation of the generalised capacitance matrix. In both cases, the relative size of the imaginary parts are fixed and they are rescaled by some $\tau$, which is such that the asymptotic exceptional point occurs at $\tau=1$.} \label{fig:EP3}
\end{figure}

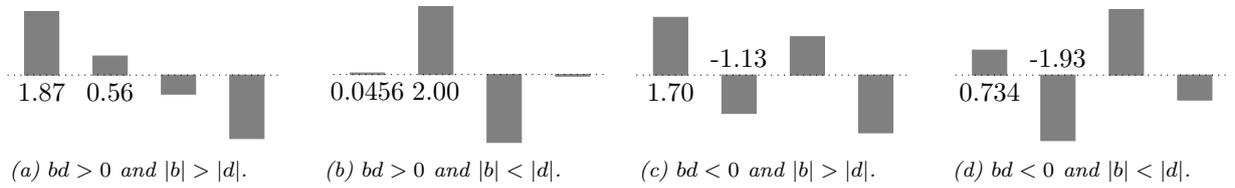
\begin{figure}
	\begin{subfigure}[b]{0.2\linewidth}
		\centering
		\begin{tikzpicture}[scale=0.45]
		\begin{scope}
		\def\b{1.87}
		\def\d{0.56}
		\draw[fill=gray,gray] (0,0) rectangle (1,\b);
		\node at (0.5,-0.5) {\b};
		\draw[fill=gray,gray] (2,0) rectangle (3,\d);
		\node at (2.5,-0.5) {\d};
		\draw[fill=gray,gray] (4,0) rectangle (5,-\d);
		\draw[fill=gray,gray] (6,0) rectangle (7,-\b);
		\draw[dotted] (-0.5,0) -- (7.5,0);
		\path (-0.5,-2) -- (-0.5,2);
		\end{scope}
		\end{tikzpicture}
		%\captionsetup{type=figure}
		\caption{$bd > 0$ and $|b| > |d|$.} \label{fig:sol4a}
	\end{subfigure}
	\hspace{20pt}
	\begin{subfigure}[b]{0.2\linewidth}
		\centering
		\begin{tikzpicture}[scale=0.45]
		\begin{scope}
		\def\b{0.0456}
		\def\d{2.00}
		\draw[fill=gray,gray] (0,0) rectangle (1,\b);
		\node at (0.5,-0.5) {\b};
		\draw[fill=gray,gray] (2,0) rectangle (3,\d);
		\node at (2.5,-0.5) {\d};
		\draw[fill=gray,gray] (4,0) rectangle (5,-\d);
		\draw[fill=gray,gray] (6,0) rectangle (7,-\b);
		\draw[dotted] (-0.5,0) -- (7.5,0);
		\path (-0.5,-2) -- (-0.5,2);
		\end{scope}
		\end{tikzpicture}
		%\captionsetup{type=figure}
		\caption{$bd > 0$ and $|b| < |d|$.}
	\end{subfigure}
	\hspace{20pt}
	\begin{subfigure}[b]{0.2\linewidth}
		\centering
		\begin{tikzpicture}[scale=0.45]
		\begin{scope}
		\def\b{1.70}
		\def\d{-1.13}
		\draw[fill=gray,gray] (0,0) rectangle (1,\b);
		\node at (0.5,-0.5) {\b};
		\draw[fill=gray,gray] (2,0) rectangle (3,\d);
		\node at (2.5,0.5) {\d};
		\draw[fill=gray,gray] (4,0) rectangle (5,-\d);
		\draw[fill=gray,gray] (6,0) rectangle (7,-\b);
		\draw[dotted] (-0.5,0) -- (7.5,0);
		\path (-0.5,-2) -- (-0.5,2);
		\end{scope}
		\end{tikzpicture}
		%\captionsetup{type=figure}
		\caption{$bd < 0$ and $|b| > |d|$.}
	\end{subfigure}
	\hspace{20pt}
	\begin{subfigure}[b]{0.2\linewidth}
		\centering
		\begin{tikzpicture}[scale=0.45]
		\begin{scope}
		\def\b{0.734}
		\def\d{-1.93}
		\draw[fill=gray,gray] (0,0) rectangle (1,\b);
		\node at (0.5,-0.5) {\b};
		\draw[fill=gray,gray] (2,0) rectangle (3,\d);
		\node at (2.5,0.5) {\d};
		\draw[fill=gray,gray] (4,0) rectangle (5,-\d);
		\draw[fill=gray,gray] (6,0) rectangle (7,-\b);
		\draw[dotted] (-0.5,0) -- (7.5,0);
		\path (-0.5,-2) -- (-0.5,2);
		\end{scope}
		\end{tikzpicture}
		%\captionsetup{type=figure}
		\caption{$bd < 0$ and $|b| < |d|$.} \label{fig:sol4d}
	\end{subfigure}
	\caption{A system of four $\mathcal{PT}$-symmetric resonators supports four asymptotic exceptional points, each with its own symmetry. Here, we plot the leading order coefficients of the imaginary parts of the material parameters (the gain or loss) at each of the four asymptotic exceptional points, using the notation $b=\mathrm{Im}(v_1^2\delta_1)$ and $d=\mathrm{Im}(v_2^2\delta_2)$.}
	\label{fig:sols4}
\end{figure}

\begin{figure}
	\centering
	\begin{subfigure}{0.3\linewidth}
		\includegraphics[width=\linewidth]{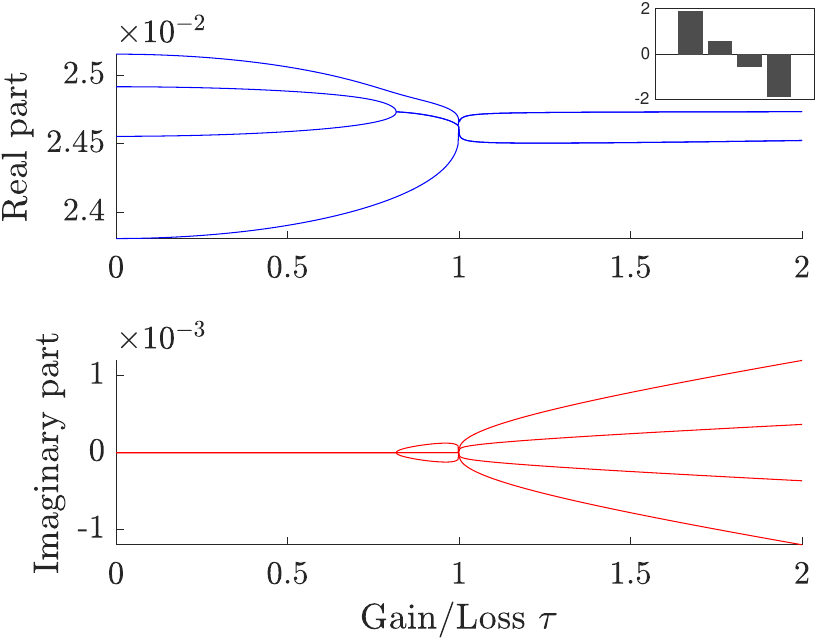}
		\caption{$N=4$}
	\end{subfigure}
	\hspace{0.2cm}
	\begin{subfigure}{0.3\linewidth}
		\includegraphics[width=\linewidth]{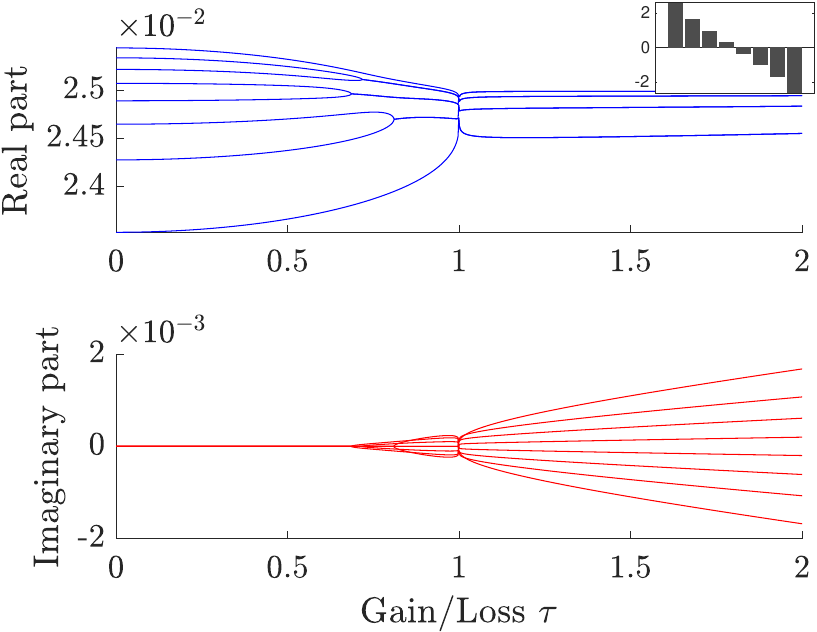}
		\caption{$N=8$}
	\end{subfigure}
	\hspace{0.2cm}
	\begin{subfigure}{0.3\linewidth}
		\includegraphics[width=\linewidth]{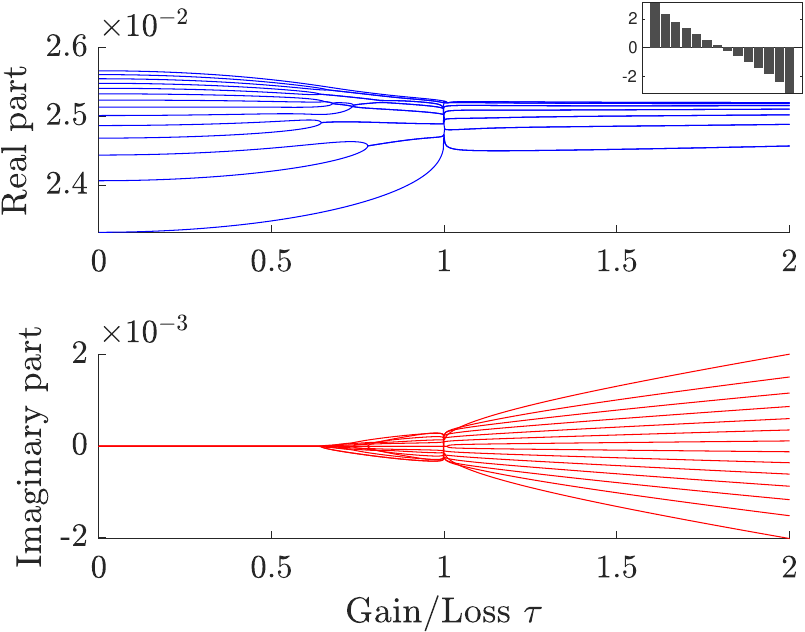}
		\caption{$N=14$}
	\end{subfigure}

	\vspace{0.3cm}

	\begin{subfigure}{0.28\linewidth}
		\includegraphics[width=\linewidth]{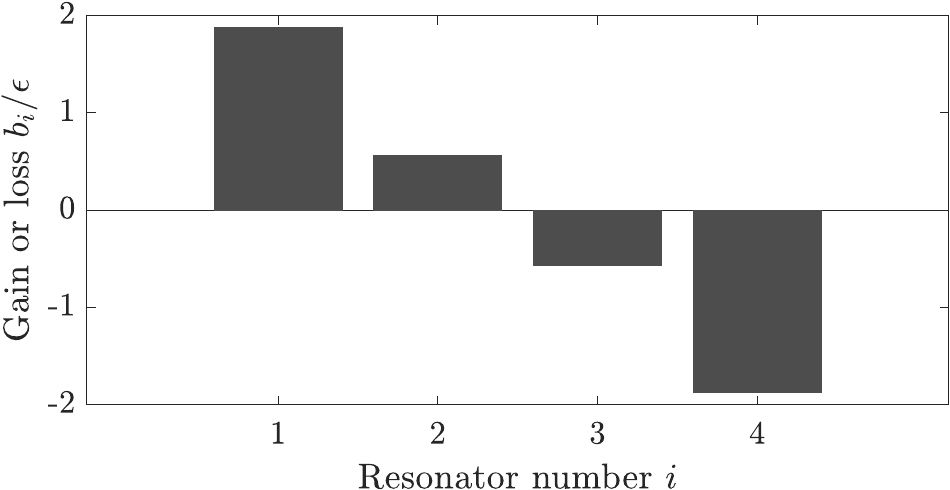}
		\caption{$N=4$} \label{fig:EPsymm4}
	\end{subfigure}
	\hspace{0.4cm}
	\begin{subfigure}{0.28\linewidth}
		\includegraphics[width=\linewidth]{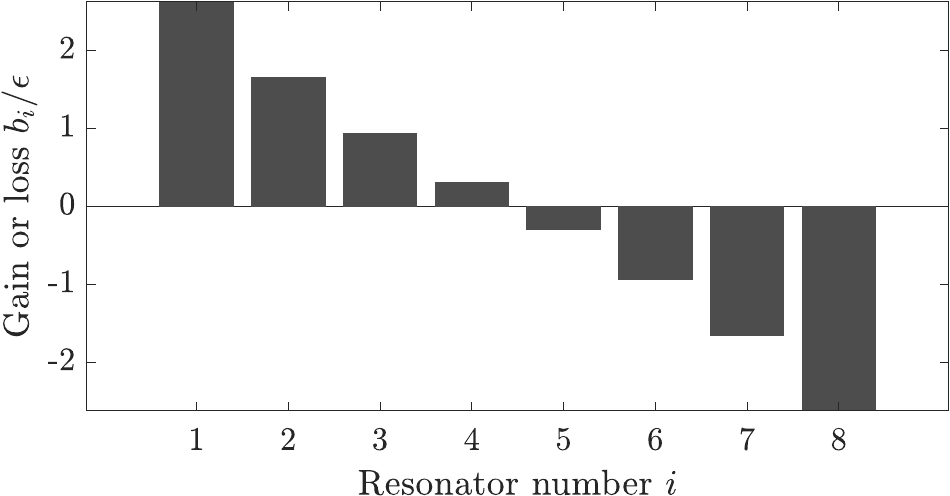}
		\caption{$N=8$} \label{fig:EPsymm8}
	\end{subfigure}
	\hspace{0.4cm}
	\begin{subfigure}{0.28\linewidth}
		\includegraphics[width=\linewidth]{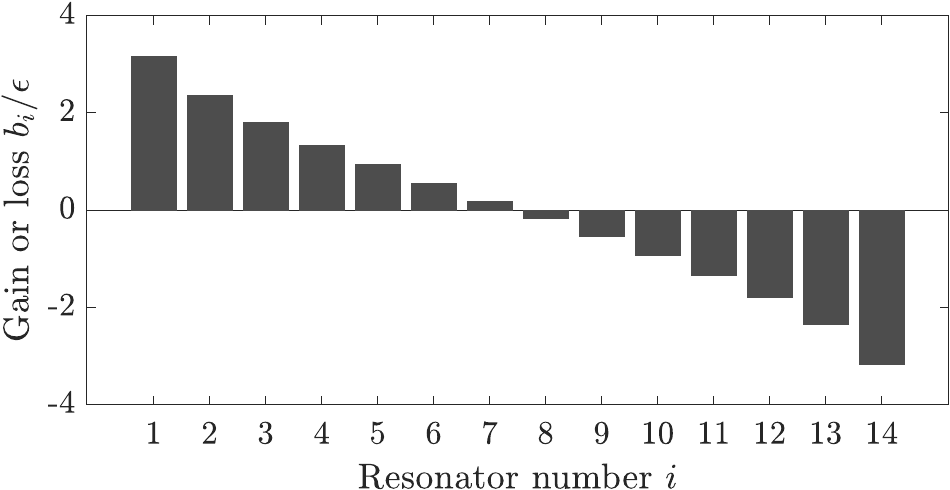}
		\caption{$N=14$} \label{fig:EPsymm14}
	\end{subfigure}

	\vspace{0.3cm}

	\begin{subfigure}{0.28\linewidth}
		\includegraphics[width=\linewidth]{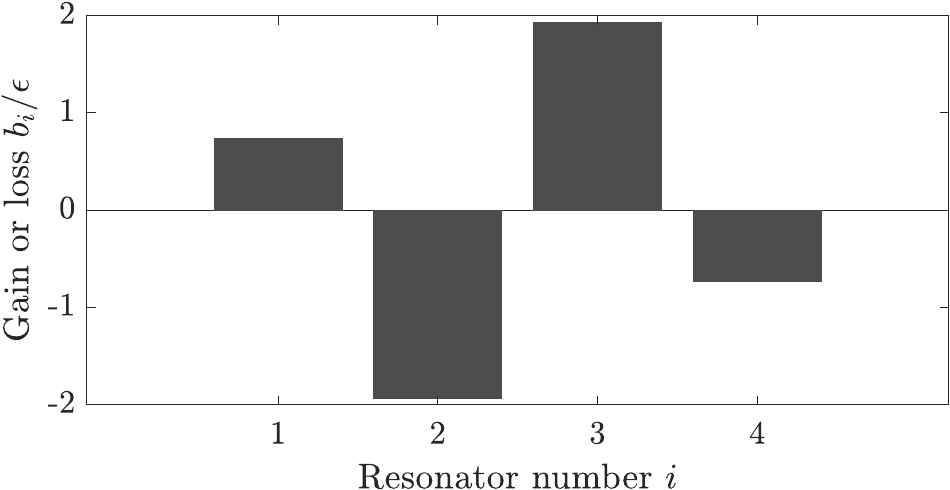}
		\caption{$N=4$}
	\end{subfigure}
	\hspace{0.4cm}
	\begin{subfigure}{0.28\linewidth}
		\includegraphics[width=\linewidth]{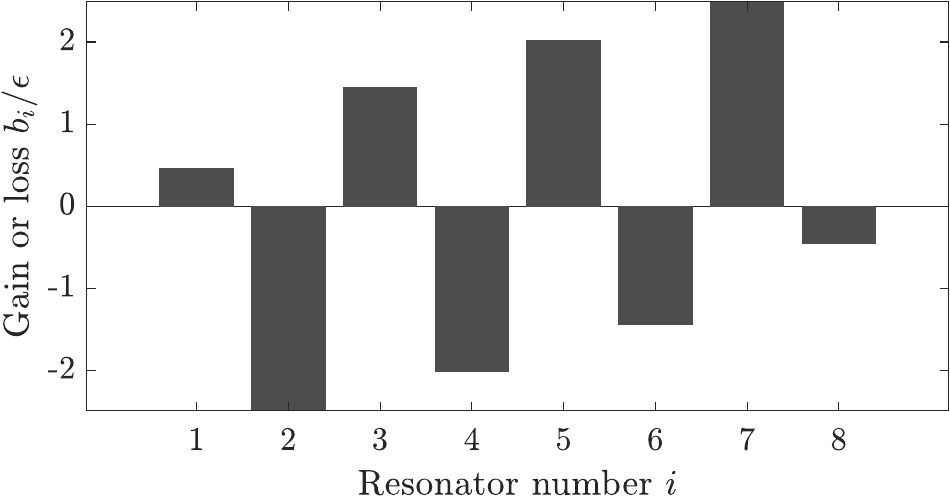}
		\caption{$N=8$}
	\end{subfigure}
	\hspace{0.4cm}
	\begin{subfigure}{0.28\linewidth}
		\includegraphics[width=\linewidth]{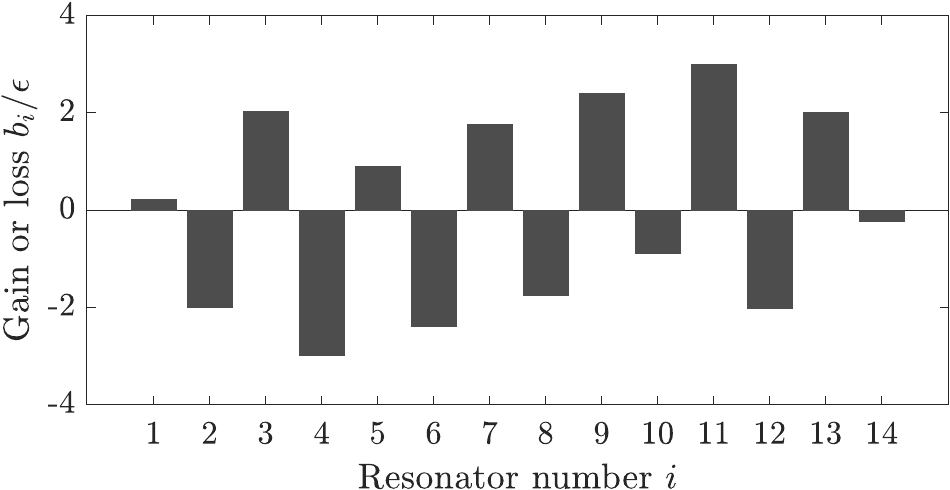}
		\caption{$N=14$}
	\end{subfigure}
	\caption{Higher-order asymptotic exceptional points can be found in larger arrays of resonators.} \label{fig:highEPs}
\end{figure}

\subsection{Subwavelength guiding of waves} \label{sec:defects}
The study of point defects in band gap materials has had immense impact on technological applications such as wireless communications, biomedical super-resolution imaging and quantum computing \cite{pointdefect1,pointdefect2}. The most notable examples are doped semiconductors, where conducting modes are induced by impurities in the semiconductor material. As we have seen, band gaps can be found in any type of wave-propagation systems, and defects in a subwavelength band gap material can be used to enable trapped or guided waves on very small length-scales.

To fix the setting, we consider the Helmholtz resonance problem in the fully periodic case $d=d_l=2$. We study the equation \eqref{eq:scattering_periodic}, where $D$ consists of a single, circular resonator inside $Y$. However, instead of considering the periodic crystal $\Dc$, we detune the size of some resonators, thereby creating a defect. We will consider the two defects illustrated in \Cref{fig:defect}, where either a single resonator or a line of resonators are detuned. We choose a square lattice with unit cell $Y$
$$Y =\left[-\frac12, \frac12\right)\times\left[-\frac12, \frac12\right),$$
and let $D$ be a circle of radius $R$ and $D_d$ a circle of radius $R+\epsilon$ for some $-R<\epsilon<1-R$. Then we define the defect crystals
$$\Dc_\mathrm{pt} = \Bigg(\bigcup_{m\in \Z^2\setminus\{(0,0)\}} D+m\Bigg)\cup D_d,\qquad \Dc_\mathrm{ln} = \Bigg(\bigcup_{m_1\in \Z, m_2\in \Z\setminus\{0\}} D+(m_1,m_2)\Bigg)\cup \Bigg(\bigcup_{m\in \Z\times\{0\}} D_d+m\Bigg).$$
We will always consider defect structures in relation to a corresponding unperturbed, periodic structure, in this case $\Dc$.

\begin{defn}[Subwavelength band-gap frequency] \label{defn:localized}
	A subwavelength resonant frequency of a resonator structure with a defect is called a subwavelength band-gap frequency if it lies inside a band gap of the unperturbed structure.
\end{defn}
We remark that since we consider the fully periodic case, both the subwavelength bands and the  subwavelength band-gap frequencies are real \cite{ammari2017subwavelength,ammari2018subwavelength}.

\begin{figure}[htb]
	\begin{subfigure}[b]{0.45\linewidth}
		\centering
	\begin{tikzpicture}[scale=1.2]
		\draw[dashed] (0,0) circle (10pt) node{$D_d$};
		\draw (0,0) circle (6pt) node[yshift=-13pt, xshift=11pt ]{$D$};
		\draw (1,0) circle (10pt);
		\draw (0,1) circle (10pt);
		\draw (1,1) circle (10pt);
		\draw (-1,0) circle (10pt);
		\draw (0,-1) circle (10pt);
		\draw (1,-1) circle (10pt);
		\draw (-1,1) circle (10pt);
		\draw (-1,-1) circle (10pt);
		\draw (1.6,0) node{$\cdots$};
		\draw (-1.55,0) node{$\cdots$};
		\draw (0,1.6) node{$\vdots$};
		\draw (0,-1.5) node{$\vdots$};
		\draw (1.7,1) node[xshift=2pt]{$v$};
		\draw (1.35,1.3) node{$\delta$};
		\draw (1,1) node{$v_\mathrm{r}$};
		\draw[opacity = 0.3] (0.5,0.5) -- (0.5,-0.5) node[yshift=-2pt, xshift=5pt ]{$Y$} -- (-0.5,-0.5) -- (-0.5,0.5) -- cycle;
	\end{tikzpicture}
	\caption{Point-defect crystal.} \label{fig:defect_point}
	\end{subfigure}\hspace{0.1cm}
	\begin{subfigure}[b]{0.45\linewidth}
	\centering
	\begin{tikzpicture}[scale=1.2]
		\draw[dashed] (0,0) circle (10pt) node[yshift=11pt, xshift=11pt ]{$D$};
		\draw (0,0) circle (6pt) node{$D_d$};
		\draw[dashed] (-2,0) circle (10pt);
		\draw[dashed] (-1,0) circle (10pt);
		\draw[dashed] (2,0) circle (10pt);
		\draw[dashed] (1,0) circle (10pt);
		\draw (-2,0) circle (6pt);
		\draw (-1,0) circle (6pt);
		\draw (2,0) circle (6pt);
		\draw (1,0) circle (6pt);

		\draw (0,1) circle (10pt);
		\draw (1,1) circle (10pt);
		\draw (0,-1) circle (10pt);
		\draw (1,-1) circle (10pt);
		\draw (-1,1) circle (10pt);
		\draw (-1,-1) circle (10pt);
		\draw (2,1) circle (10pt);
		\draw (2,-1) circle (10pt);
		\draw (-2,1) circle (10pt);
		\draw (-2,-1) circle (10pt);
		\draw (2.65,0) node{$\cdots$};
		\draw (-2.6,0) node{$\cdots$};
		\draw (0,1.65) node{$\vdots$};
		\draw (0,-1.5) node{$\vdots$};
		\draw (2.8,1) node{$v$};
		\draw (2,1) node{$v_\mathrm{r}$};
		\draw (2.35,1.3) node{$\delta$};

		\draw[opacity = 0.5] (0.5,0.5) -- (0.5,-0.5) node[yshift=-5pt, xshift=-7pt ]{$Y$} -- (-0.5,-0.5) -- (-0.5,0.5) -- cycle;
		\draw[opacity = 0.5] (-0.5,-1.8) -- (-0.5,1.8)
		(0.5,-1.8) node[right]{$Y_\mathrm{strip}$} -- (0.5,1.8);
	\end{tikzpicture}
	\caption{Line-defect crystal.} \label{fig:defect_line}
\end{subfigure}
\caption{Illustration of the two different defect crystals, and the material parameters, in the case of smaller defect resonator.} \label{fig:defect}
\end{figure}
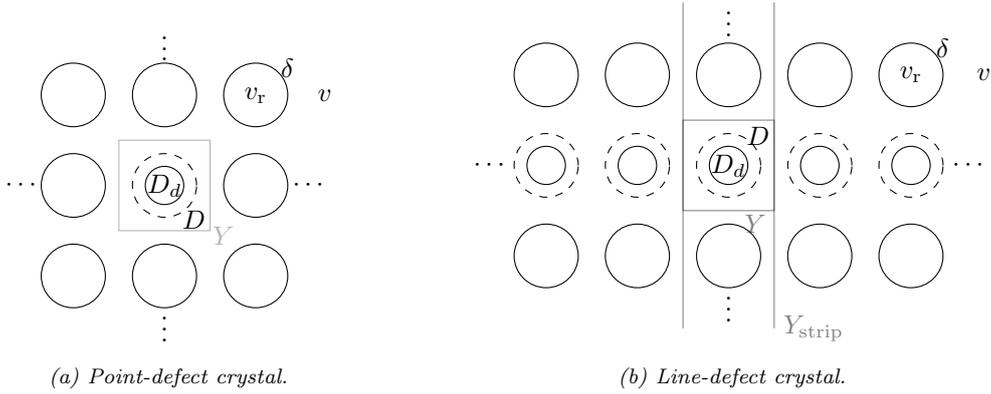

\subsubsection{Fictitious source method for point defects}
We wish to solve \eqref{eq:scattering_periodic} where the periodic crystal $\Dc$ is replaced by $\Dc_\mathrm{pt}$ or $\Dc_\mathrm{ln}$. Firstly, we observe that $\Dc_\mathrm{ln}$ is periodic in the $x_1$-direction, so corresponding problem can be reduced to the quasiperiodic problem in the strip $Y_\mathrm{strip} = \left[-\frac12,\frac12\right)\times \R$ by Floquet-Bloch theory. This makes the analysis of the two problems ($\Dc_\mathrm{pt}$ on $\R^2$ or $\Dc_\mathrm{ln}$ on $Y_\mathrm{strip}$) conceptually very similar. We will outline the method of \emph{fictitious sources} in the case of a point defect. For $d=d_l=2$, we are studying the problem
\begin{equation} \label{eq:scattering_defect}
	\left\{
	\begin{array} {ll}
		\ds \Delta u+ k^2 u  = 0  &\text{in } \R^2 \backslash \Dc_\mathrm{pt}, \\
		\nm
		\ds \Delta u+ k_\mathrm{r}^2 u  = 0  &\text{in } \Dc_\mathrm{pt}, \\
		\nm
		\ds  u|_{+} -u|_{-}  =0   &\text{on } \partial \Dc_\mathrm{pt}, \\[0.3em]
		\nm
		\ds  \delta \frac{\partial u}{\partial \nu} \bigg|_{+} - \frac{\partial u}{\partial \nu} \bigg|_{-} =0 &\text{on }\partial \Dc_\mathrm{pt} ,
	\end{array}
	\right.
\end{equation}
where $k_\mathrm{r}=\omega/v_\mathrm{r}$ and $v_\mathrm{r}$ is the wave speed inside the resonators. The idea is to replace the defected resonator with an unperturbed resonator, along with fictitious sources which are designed to make the new problem equivalent to the original. More precisely, we consider the problem
\begin{equation} \label{eq:scattering_fictitious}
	\left\{
	\begin{array} {ll}
		\ds \Delta \widetilde u+ k^2\widetilde u  = 0  &\text{in } \R^2 \backslash \Dc, \\
		\nm
		\ds \Delta  \widetilde u+ k_\mathrm{r}^2\widetilde u  = 0  &\text{in } \Dc, \\
		\nm
		\ds  \widetilde u|_{+} - \widetilde u|_{-}  = f\chi_{\p D} &\text{on } \partial \Dc, \\[0.3em]
		\nm
		\ds  \delta \frac{\partial\widetilde u}{\partial \nu} \bigg|_{+} - \frac{\partial\widetilde u}{\partial \nu} \bigg|_{-} =g\chi_{\p D} &\text{on }\partial \Dc,
	\end{array}
	\right.
\end{equation}
where $f,g$ are the source terms and $\chi_{\p D}$ is, as usual, the characteristic function of $\p D$. Note that \eqref{eq:scattering_fictitious} is posed on the periodic geometry $\Dc$, and that the non-zero sources are present only on the boundary of the central resonator $D$. Note also that since we are  in the fully-periodic case ($d=d_l)$, we do not assume any radiation condition in \eqref{eq:scattering_defect} and \eqref{eq:scattering_fictitious}. We seek solutions $u$ of \eqref{eq:scattering_defect} and $\widetilde u$ of  \eqref{eq:scattering_fictitious} for $\omega$ inside the band gap of the periodic problem. Inside the unit cell $Y$, we can represent the solution $\widetilde u$ as
\begin{align}\label{eq:u_rep_effective}
	\widetilde{u}
	=
	\begin{cases}
		H + \S_D^{k}[\phi] &\quad \mbox{in } Y\setminus \overline{D},\\
		\S_D^{k_\mathrm{r}}[\psi] &\quad \mbox{in } D,
	\end{cases}
\end{align}
where $H$ satisfies the homogeneous equation $\Delta H + k^2H = 0$ in $Y$. Using this approach, we can explicitly compute $f,g$ in terms of $\psi,\phi$ so that $\widetilde u$ coincides with $u$ in $Y\setminus \left(D\cap D_d\right)$.
\begin{lemma}
	The density pair $(\psi,\phi)$ and the effective sources $(f,g)$ satisfy the following relation
	\begin{equation}\label{eq:effective}
		(\mathcal{A}^\epsilon - \mathcal{A})
		\begin{pmatrix}
			\psi\\ \phi
		\end{pmatrix} =
		\begin{pmatrix}
			f\\
			g
		\end{pmatrix},
	\end{equation}
	where $\mathcal{A}^\epsilon$ is defined as
	\begin{equation}\label{eq:ADd}
		\mathcal{A}^\epsilon  := (\mathcal{P}_2)^{-1}\mathcal{A}_{D_d} \mathcal{P}_1.
	\end{equation}
	Here, $\A$ and $\A_{D_d}$ are the operators defined in \Cref{lem:BIE_finite} for the domains $D$ and $D_d$, respectively, and the operators $\mathcal{P}_1: L^2(\p D)^2\rightarrow L^2(\p D_d)^2$ and $\mathcal{P}_2: L^2(\p D)^2\rightarrow L^2(\p D_d)^2$ are defined by
	\begin{align*}
		\mathcal{P}_1\begin{pmatrix} e^{\i n \theta} \\ e^{\i m \theta}\end{pmatrix} &= \delta_{mn}
		\frac{R}{R_d}\begin{pmatrix}
			\ds\frac{H_n^{(1)}(k_\mathrm{r} R) }{H_n^{(1)}(k_\mathrm{r} R_d)}e^{\i n \theta}
			\\[0.8em]
			\ds\frac{J_n(k R) }{J_n(k R_d)}e^{\i n \theta}
		\end{pmatrix}, &
		\mathcal{P}_2\begin{pmatrix} e^{\i n \theta} \\ e^{\i m \theta}\end{pmatrix} &= \delta_{mn} \begin{pmatrix}
			\ds \frac{J_n(k R_d) }{J_n(k R)}e^{\i n \theta}
			\\[0.8em]
			\ds\frac{J_n'(k R_d) }{J_n'(k R)}e^{\i n \theta}
		\end{pmatrix}.
	\end{align*}
\end{lemma}
As in previous settings, we can use an integral equation formulation of the problem. In the defect crystal setting, however, the formulation is slightly different than previously  \cite{ammari2018subwavelength}.
\begin{lemma} \label{prop:fg}
	The subwavelength band-gap frequencies of \eqref{eq:scattering_defect} are precisely the characteristic values $\omega = \omega^\epsilon(\delta)$ of the operator
	\begin{align} \label{eq:Mdensity}
		\mathcal{M}^{\epsilon}(\omega,\delta)	=I+\left(\frac{1}{(2\pi)^2}\int_{Y^*}\A^\alpha(\omega,\delta)^{-1} \dx\alpha\right)\big(\A_{D}^\epsilon(\omega,\delta) - \A(\omega,\delta)\big)
	\end{align}
	inside the band gap of $\Dc$, such that $\omega^\epsilon \to 0$ as $\delta\to0$. Here, $\A^\alpha$ is the operator defined in  \Cref{lem:BIE_periodic}, which is invertible for $\omega$ inside the band gap.
\end{lemma}
We let $\omega^{*} = \omega^{*}(\delta)$ be the maximum of the first band:
$$\omega^{*}(\delta) = \max_{\alpha\in Y^*} \omega_1^\alpha(\delta).$$
The following theorem demonstrates the existence of band-gap frequencies in the point defect crystal \cite{ammari2018subwavelength}.
\begin{thm} \label{thm:defect_pt}
	Assume that $\delta$ and $\epsilon$ are small enough and the pair $(R,\epsilon)$ satisfies one of the two assumptions:
	\begin{itemize}
		\item[(i)] $R$ small enough and $\epsilon<0$ (dilute regime);
		\item[(ii)] $R$ close enough to $1/2$ and $\epsilon>0$ (non-dilute regime).
	\end{itemize}
	Then there exists a subwavelength band-gap frequency $\omega^\epsilon$ of \eqref{eq:scattering_defect}. In both cases we have the asymptotic expansion
	\begin{align*}
		\omega^\epsilon-\omega^* =
		\exp\left(-\frac{\mu}{\delta\epsilon}+O\left(\frac{1}{\epsilon\ln\delta}\right)\right),
	\end{align*}
	when $\epsilon$ and $\delta$ go to zero, for some constant $\mu>0$.
\end{thm}
\begin{remark}
	An explicit expression for $\mu$ is given in \cite{ammari2018subwavelength}.
\end{remark}
\begin{remark}
	In the reverse cases, \textit{i.e.} if $R$ is close to $1/2$ and $\epsilon<0$ or if $R$ is small and $\epsilon>0$, there will be no band-gap frequencies. In these cases, the resonance frequency of the defected resonator is shifted downwards, and will therefore lie outside of the band gap.
\end{remark}

\begin{figure}[tbh]
	\begin{subfigure}[b]{0.45\linewidth}
\includegraphics[trim={0 0.29cm 0 0.9cm}, clip,width=0.97\linewidth]{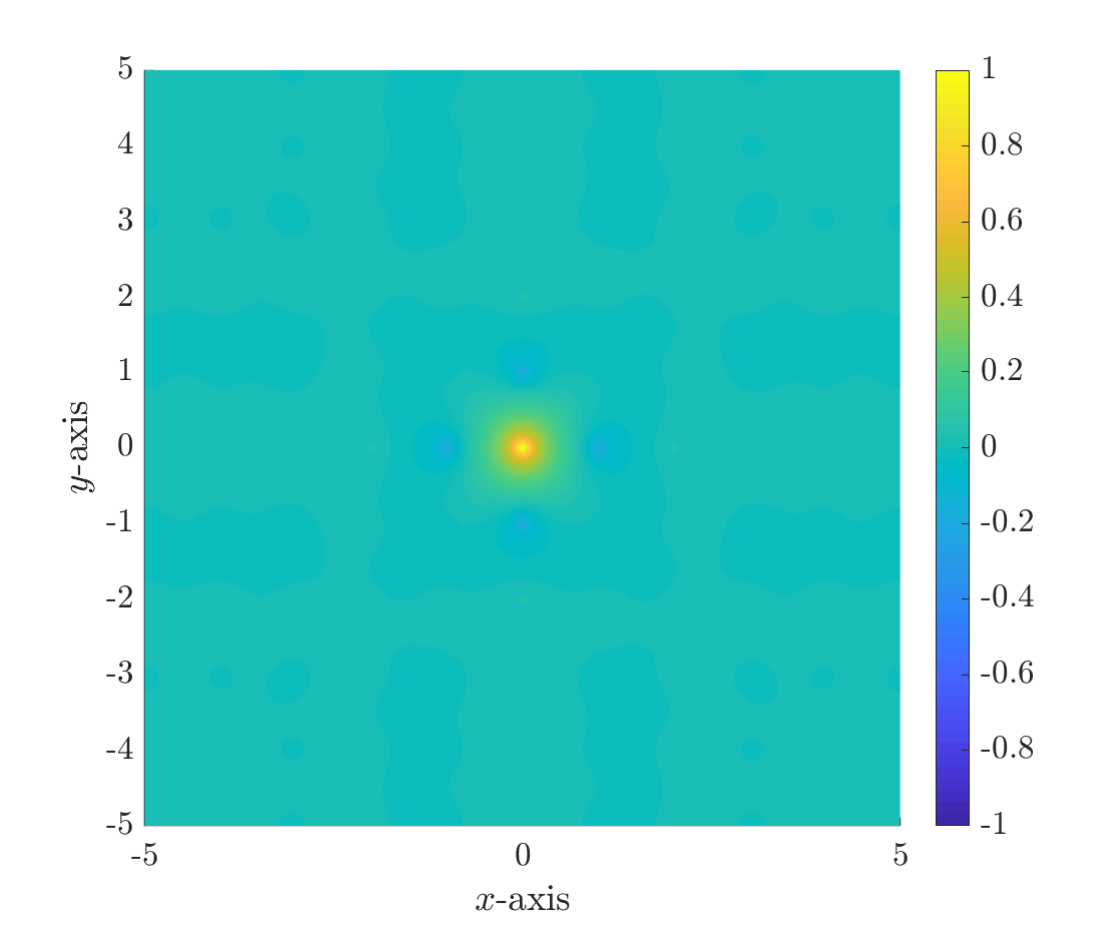}
\caption{Localised mode.}\label{fig:loc}
\end{subfigure}\hfill
\begin{subfigure}[b]{0.45\linewidth}
\includegraphics[width=\linewidth]{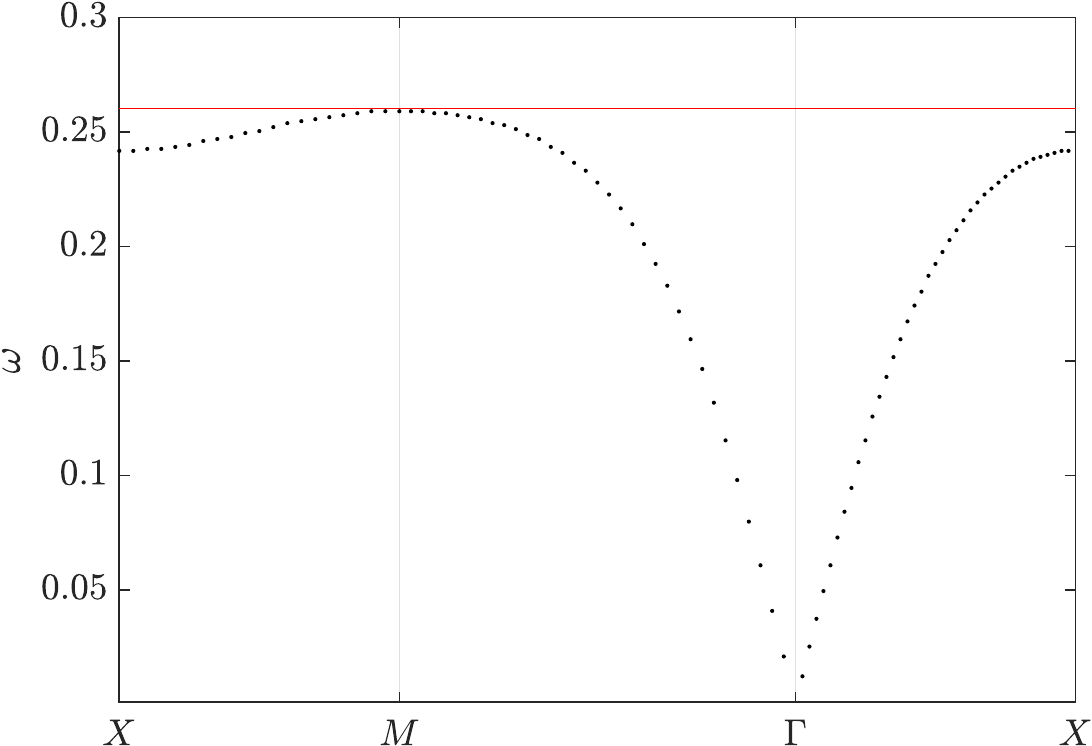}
\caption{Subwavelength bulk band (dotted black) and band-gap frequency (solid red).}\label{fig:bandgapf}
\end{subfigure}
\caption{Localised mode on subwavelength scales (a) and corresponding band-gap frequency (b). These point defect structures can provide a good degree of localisation. However, as is apparent from (b), the band-gap frequency is exponentially close to the edge of the bulk band and is therefore sensitive to imperfections of the structure.}\label{fig:ptdefect}
\end{figure}

\subsubsection{Guided waves in line defects}
We let $\Dc_\mathrm{ln}$ be the line defect crystal (see \Cref{fig:defect_line}) and study the problem
\begin{equation} \label{eq:scattering_ln}
	\left\{
	\begin{array} {ll}
		\ds \Delta u+ k^2 u  = 0  &\text{in } \R^2 \backslash \Dc_\mathrm{ln}, \\
		\nm
		\ds \Delta u+ k_\mathrm{r}^2 u  = 0  &\text{in } \Dc_\mathrm{ln}, \\
		\nm
		\ds  u|_{+} -u|_{-}  =0   &\text{on } \partial \Dc_\mathrm{ln}, \\[0.3em]
		\nm
		\ds  \delta \frac{\partial u}{\partial \nu} \bigg|_{+} - \frac{\partial u}{\partial \nu} \bigg|_{-} =0 &\text{on }\partial \Dc_\mathrm{ln}, \\
		\nm
		\ds u(x_1,x_2)e^{-\i \alpha_1 x_1} \text{ is periodic in $x_1$,}
	\end{array}
	\right.
\end{equation}
where $\alpha=(\alpha_1,\alpha_2).$
In this case, since the defect is periodic in the $x_1$-direction, we can reduce \eqref{eq:scattering_ln} to a Helmholtz resonance  problem posed on the infinite strip $Y_{\mathrm{strip}}$ shown in \Cref{fig:defect_line} and have a very similar characterisation of the band-gap frequencies to the one given in \Cref{prop:fg}.
\begin{lemma} \label{prop:M_ln}
	The subwavelength band-gap frequencies of \eqref{eq:scattering_defect} are precisely the characteristic values $\omega = \omega^\epsilon(\delta,\alpha_1)$ of the operator
	\begin{equation}\label{eq:Mline}
		\mathcal{M}^{\epsilon,\alpha_1}(\omega,\delta)	=I+\left(\frac{1}{2\pi}\int_{-\pi}^\pi\A^{(\alpha_1,\alpha_2)}(\omega,\delta)^{-1} \dx\alpha_2\right)\big(\A_{D}^\epsilon(\omega,\delta) - \A(\omega,\delta)\big)
	\end{equation}
	inside the band gap of $\Dc$, such that $\omega^\epsilon \to 0$ as $\delta\to0$.
\end{lemma}
We now let $\omega^{\alpha_1,*}$ be the maximum of the first band at $\alpha=(\alpha_1,\alpha_2)$:
$$\omega^{\alpha_1,*}(\delta) = \max_{\alpha_2\in[-\pi,\pi]} \omega_1^{(\alpha_1,\alpha_2)}(\delta).$$
Then we have the following result \cite{ammari2020subwavelength}.
\begin{thm} \label{thm:defect_ln} 	Assume that $\delta$ and $\epsilon$ are small enough and the pair $(R,\epsilon)$ satisfies one of the two assumptions:
	\begin{itemize}
		\item[(i)] $R$ small enough and $\epsilon<0$ (dilute regime);
		\item[(ii)] $R$ close enough to $1/2$ and $\epsilon>0$ (non-dilute regime).
	\end{itemize}
	Then there exists a subwavelength resonant frequency $\omega^\epsilon$ of \eqref{eq:scattering_ln} satisfying $\omega^\epsilon > \omega^{\alpha_1,*}$. Moreover, as $\delta,\epsilon \rightarrow 0$ we have
	\begin{equation} \label{eq:defect-freq}
		\omega^\epsilon(\delta,\alpha_1) = \omega^{\alpha_1,*}(\delta) + \mu(\alpha_1)\sqrt{\delta}\epsilon^2 +
		O\left(\epsilon^2\sqrt{\delta}\left(\frac{1}{\ln\delta} + \epsilon\right)\right)
%		 O\left(\frac{\epsilon^2\sqrt{\delta}}{\ln\delta} + \epsilon^3\sqrt{\delta}\right),
	\end{equation}
	for some $\mu=\mu(\alpha_1)>0$ which is independent of $\epsilon$ and $\delta$.
\end{thm}
We call $\omega^\epsilon$, viewed as a function of $\alpha_1$, a \emph{defect band}. We emphasize that parts of the defect band might not correspond to band-gap frequencies in the sense of \Cref{defn:localized}: from \Cref{thm:defect_ln} we know that $\omega^\epsilon(\alpha_1) > \omega^{\alpha_1,*}$ while in order for $\omega^\epsilon$ to lie in the band gap of $\Dc$ we need $\omega^\epsilon(\alpha_1) > \omega^*$ where, as before, $\omega^* = \max_{\alpha_1\in[-\pi,\pi]} \omega_1^{\alpha_1,*}$.

In order for the whole defect band to lie in the subwavelength band gap, we need a sufficiently large $\epsilon$. Since \Cref{thm:defect_ln} is based on asymptotically expanding $\M$ for small $\epsilon$, a different analysis is needed to handle this case. The following result is based on asymptotics in the dilute regime, and is valid even for $\epsilon$ with large magnitude \cite{ammari2020subwavelength}.
\begin{thm} \label{thm:dilute}
	For $\delta$ and $R$ small enough, and for fixed $\epsilon\in(-R,0)$, there exists a unique subwavelength resonant frequency $\omega^\epsilon$ of \eqref{eq:scattering_ln} satisfying $\omega^\epsilon > \omega^{\alpha_1,*}$. For $\alpha_1 \neq 0$, $$\omega^\epsilon(\alpha_1) = \hat\omega+O\left(R^2 + \delta\right),$$ where $\hat\omega$ is the root of the following equation:
	\begin{equation} \label{eq:dilute}
		1 + \left(\frac{\hat\omega^2R^2}{2\delta}\ln\frac{R}{R_d} + \left(1-\frac{R^2}{R_d^2} \right)\right)\frac{1}{2\pi}\int_{-\pi}^\pi \frac{(\omega^\alpha)^2}{\hat\omega^2-(\omega^\alpha)^2}\dx \alpha_2 = 0.
	\end{equation}
	For $\delta$ and $R$ small enough, and for fixed $\epsilon \in[0,1-R)$, there are no resonant frequencies satisfying $\omega^\epsilon > \omega^{\alpha_1,*}$.
\end{thm}

\begin{prop}
For $R$ and $\delta$ small enough, there exists an $\epsilon_0>0$ such that for any $\epsilon\in(-R,-\epsilon_0)$ we have
$$\omega^\epsilon(\alpha_1) > \omega^*$$
for all $\alpha_1\in [-\pi,\pi]$.
\end{prop}
In order for the line defect crystal to be useful as a waveguide, we need the localised modes to propagate along the defect line. In other words, we must exclude the case of \emph{bound modes}, which are modes that are localised along the direction of the line. If there is such a mode $u$, corresponding to a frequency $\omega$, we can apply the Floquet transform so that $u^\alpha(x):=\F[u](x,\alpha)$ solves \eqref{eq:scattering_defect} for any $\alpha$. Corresponding band function attains the same value $\omega$ for any $\alpha$, so we conclude that  \emph{bound modes are associated to flat band functions}. The next result from \cite{ammari2020subwavelength} shows that the defect modes in our case are not bound along the defect line.
\begin{prop} \label{prop:palpha1}
	For $\delta$ and $R$ small enough, and for $\alpha_1 \notin \{0,\pi\} $, the subwavelength resonant frequency $\omega^\epsilon = \omega^\epsilon(\alpha_1)$ satisfies
	$$\frac{\p \omega^\epsilon}{\p \alpha_1} \neq 0.$$
\end{prop}

\begin{figure}[tbh]
	\begin{subfigure}[b]{0.45\linewidth}
		\includegraphics[trim={0 0.29cm 0 0.9cm}, clip,width=0.9\linewidth]{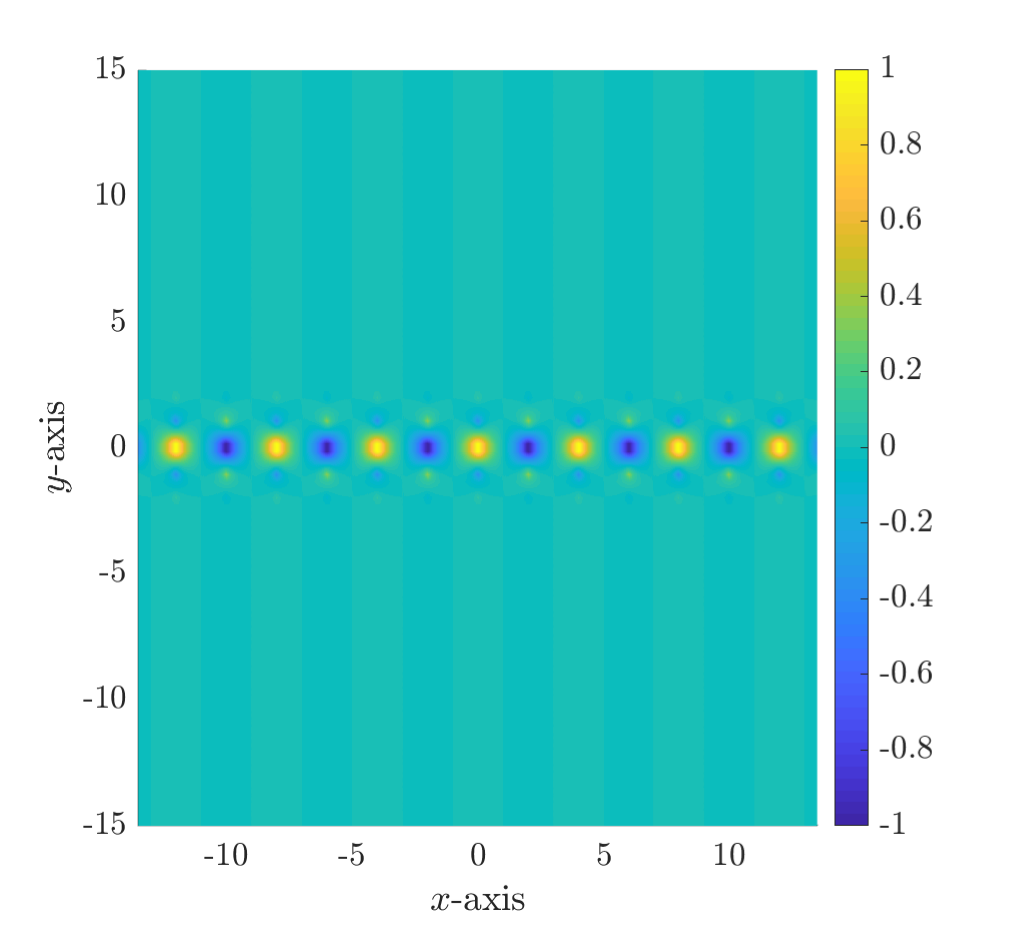}
		\caption{Guided mode.}\label{fig:loc_line}
	\end{subfigure}\hfill
	\begin{subfigure}[b]{0.45\linewidth}
		\includegraphics[width=\linewidth]{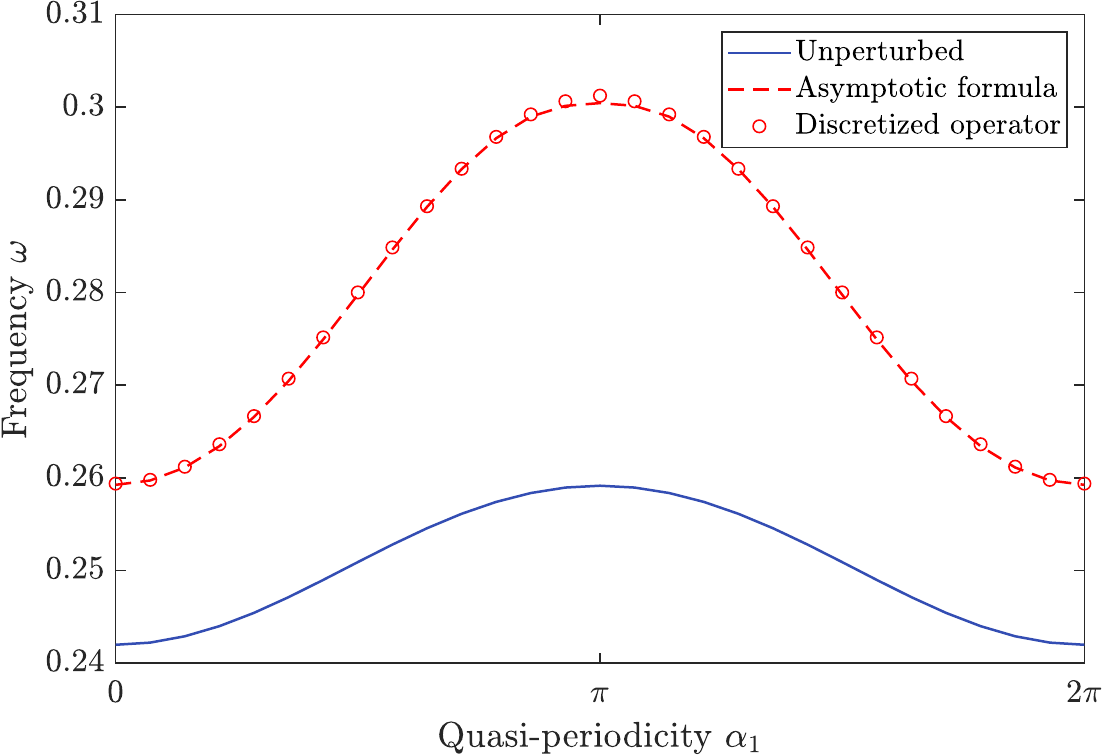}
		\caption{Subwavelength bulk band (solid blue) and defect band (dashed red).}\label{fig:bandgapline}
	\end{subfigure}
	\caption{Guided mode at subwavelength scales (a) and corresponding band-gap frequency (b). Here $\epsilon<-\epsilon_0$, which means that the whole defect band will be above the bulk band. Moreover, we see that the defect band is not flat, and corresponding modes are propagating along the line defect. The defect band was computed using the asymptotic formula \eqref{eq:dilute} (red dashed) and by discretizing the operator $\mathcal{M}^{\epsilon,\alpha_1}(\omega,\delta)$ (defined in \eqref{eq:Mline}) using the multipole method (red circles).}\label{fig:linedefect}
\end{figure}
The integral representation for the localised and guided waves (\Cref{prop:fg} and \Cref{prop:M_ln}) can be discretised using the multipole method, which provides an efficient method to compute the band-gap frequency and corresponding localised mode (see \Cref{fig:ptdefect} and \Cref{fig:linedefect}).

\subsection{Robust guiding at subwavelength scales}
There is a fundamental restriction of the practical applicability of the localised modes studied in \Cref{sec:defects}. Taking the point defect as example, we have from \Cref{thm:defect_pt} that the band-gap frequency is exponentially close to the edge of the bulk bands (see \Cref{fig:ptdefect}). When fabricating such structures, there is a large risk that small imperfections will cause the band-gap frequency to be lost inside the bulk bands, so that the desirable localisation property disappears.

For structures to have physically achievable localisation properties, such properties must be robust against imperfections of their design. To accomplish this, we take inspiration from the field of topological insulators. These are materials which can conduct current along the edges of the material, while are insulating in the bulk. The striking property of topological insulators is that these conducting edge modes originate from the structure of the bulk, rather than the edge itself. Moreover, they are localized to the edges and exhibit a remarkable robustness against perturbations of the system. These concepts have been recently studied mathematically in a variety of settings; see, for instance, \cite{feffermandislocated, fefferman1d, haibulk,haiinterface,gontier1,gontier2}.

We study a high-contrast resonator analogue of the Su-Schrieffer-Heeger (SSH) model \cite{SSH}. We consider the case of a chain of resonators, corresponding to $d=3$ and $d_l = 1$, where the unit cell is given by
$$Y=\left[-\frac{L}2,\frac{L}2\right) \times\R^2.$$
Moreover, we consider the case of a dimer of resonators $D=D_1\cup D_2$. We need two assumptions of symmetry for the analysis that follows. The first assumption is that each individual resonator is symmetric in the sense that there exists some $x_1\in\mathbb{R}$ such that
\begin{equation} \label{resonator_symmetry}
	R_1D_1 = D_1, \quad R_2D_2 = D_2,
\end{equation}
where $R_1$ and $R_2$ are the reflections in the planes $p_1=\{-x_1\}\times \R^2$ and $p_2=\{x_1\}\times \R^2$, respectively. We also assume that the dimer is symmetric in the sense that
\begin{equation} \label{dimer_symmetry}
	\P D_1= D_2,
\end{equation}
where, as before $\P(x) = -x$. We denote the resonator separation within the cell as $l$, and between the cells as $l'$, \textit{i.e.}
$$l = 2x_1, \qquad l' = L-l,$$
see \Cref{fig:SSH}. For simplicity, we choose $v=v_1=v_2 =1$, and study effects originating from the geometry of the structure. As we shall see, a topological phase transition occurs when $l$ changes across the symmetry point $l=L/2$ (corresponding to $l=l'$).

We begin by considering the periodic equation
\begin{equation} \label{eq:scattering_ssh}
	\left\{
	\begin{array} {ll}
		\ds \Delta u^\alpha+ \omega^2 {u^\alpha}  = 0 &\text{in } \R^3 \setminus \p \Dc, \\
		\nm
		\ds  {u^\alpha}|_{+} -{u^\alpha}|_{-}  =0  & \text{on } \partial \Dc, \\
		\nm
		\ds  \delta \frac{\partial {u^\alpha}}{\partial \nu} \bigg|_{+} - \frac{\partial {u^\alpha}}{\partial \nu} \bigg|_{-} =0 & \text{on } \partial \Dc, \\
		\nm
		\ds e^{-\i  \alpha_1 x_1}  u^\alpha(x_1,x_2,x_3)  \,\,\,&  \mbox{is periodic in } x_1, \\
		\nm
		\ds u^\alpha(x_1,x_2,x_3)& \text{satisfies the $\alpha$-quasiperiodic outgoing radiation condition} \\ &\hspace{0.5cm} \text{as } \sqrt{x_2^2+x_3^2} \rightarrow \infty,
	\end{array}
	\right.
\end{equation}
where $\alpha=(\alpha_1,0,0)$ and the  $\alpha$-quasiperiodic outgoing radiation condition says that the wave is the sum of a finite number of outgoing propagating plane waves and an infinite number of evanescent waves (\ie, exponentially decaying  in the directions $x_2$ and $x_3$); see, for instance, \cite[Chapter 4]{ammari2018mathematical}.
\begin{figure}[tbh]
	\centering
	\begin{tikzpicture}[scale=2]
		\pgfmathsetmacro{\rb}{0.25pt}
		\pgfmathsetmacro{\rs}{0.2pt}
		\coordinate (a) at (0.25,0);
		\coordinate (b) at (1.05,0);

		\draw[dashed, opacity=0.5] (-0.5,0.85) -- (-0.5,-1);
		\draw[dashed, opacity=0.5]  (1.8,0.85) -- (1.8,-1)node[yshift=4pt,xshift=-7pt]{};
		\draw[{<[scale=1.5]}-{>[scale=1.5]}, opacity=0.5] (-0.5,-0.6) -- (1.8,-0.6)  node[pos=0.5, yshift=-7pt,]{$L$};
		\draw plot [smooth cycle] coordinates {($(a)+(210:\rb)$) ($(a)+(270:\rs)$) ($(a)+(330:\rb)$) ($(a)+(30:\rs)$) ($(a)+(90:\rb)$) ($(a)+(150:\rs)$) }; \draw (a) node{$D_1$};
		\draw plot [smooth cycle] coordinates {($(b)+(30:\rb)$) ($(b)+(90:\rs)$) ($(b)+(150:\rb)$) ($(b)+(210:\rs)$) ($(b)+(270:\rb)$) ($(b)+(330:\rs)$) }; \draw (b) node{$D_2$};
		\draw[{<[scale=1.5]}-{>[scale=1.5]}, opacity=0.5] (0.25,0.6) -- (1.05,0.6) node[pos=0.5, yshift=-5pt,]{$l$};
		\draw[dotted,opacity=0.5] (0.25,0.7) -- (0.25,-0.8) node[at end, yshift=-0.2cm]{$p_1$};
		\draw[dotted,opacity=0.5] (1.05,0.7) -- (1.05,-0.8) node[at end, yshift=-0.2cm]{$p_2$};
		\draw[{<[scale=1.5]}-{>[scale=1.5]}, opacity=0.5] (1.05,0.6) -- (2.55,0.6) node[pos=0.6, yshift=-5pt,]{$l'$};
		%\draw[{<[scale=1.5]}-{>[scale=1.5]}, opacity=0.5] (0.25,0.6) -- (1.05,0.6) node[pos=0.5, yshift=-5pt,]{$l$};
		%\draw[dotted,opacity=0.5] (0.25,0.7) -- (0.25,-0.8) node[at end, yshift=-0.2cm]{$p_1$};
		%\draw[dotted,opacity=0.5] (1.05,0.7) -- (1.05,-0.8) node[at end, yshift=-0.2cm]{$p_2$};

		\begin{scope}[xshift=-2.3cm]
			\coordinate (a) at (0.25,0);
			\coordinate (b) at (1.05,0);
			\draw plot [smooth cycle] coordinates {($(a)+(210:\rb)$) ($(a)+(270:\rs)$) ($(a)+(330:\rb)$) ($(a)+(30:\rs)$) ($(a)+(90:\rb)$) ($(a)+(150:\rs)$) };
			\draw plot [smooth cycle] coordinates {($(b)+(30:\rb)$) ($(b)+(90:\rs)$) ($(b)+(150:\rb)$) ($(b)+(210:\rs)$) ($(b)+(270:\rb)$) ($(b)+(330:\rs)$) };
			\begin{scope}[xshift = 1.2cm]
				\draw (-1.6,0) node{$\cdots$};
			\end{scope};
		\end{scope}
		\begin{scope}[xshift=2.3cm]
			\coordinate (a) at (0.25,0);
			\coordinate (b) at (1.05,0);
			\draw plot [smooth cycle] coordinates {($(a)+(210:\rb)$) ($(a)+(270:\rs)$) ($(a)+(330:\rb)$) ($(a)+(30:\rs)$) ($(a)+(90:\rb)$) ($(a)+(150:\rs)$) };
			\draw plot [smooth cycle] coordinates {($(b)+(30:\rb)$) ($(b)+(90:\rs)$) ($(b)+(150:\rb)$) ($(b)+(210:\rs)$) ($(b)+(270:\rb)$) ($(b)+(330:\rs)$) };
			\draw[dotted,opacity=0.5] (0.25,0.7) -- (0.25,-0.8);
			\begin{scope}[xshift = 1.1cm]
			\end{scope}
			\draw (1.7,0) node{$\cdots$};
		\end{scope}

		\begin{scope}[yshift=0.9cm]
			\draw [decorate,opacity=0.5,decoration={brace,amplitude=10pt}]
			(-0.5,0) -- (1.8,0) node [black,midway]{};
			\node[opacity=0.5] at (0.67,0.35) {$Y$};
		\end{scope}
	\end{tikzpicture}
	\caption{Example of the array, drawn to illustrate the symmetry assumptions.} \label{fig:SSH}
\end{figure}
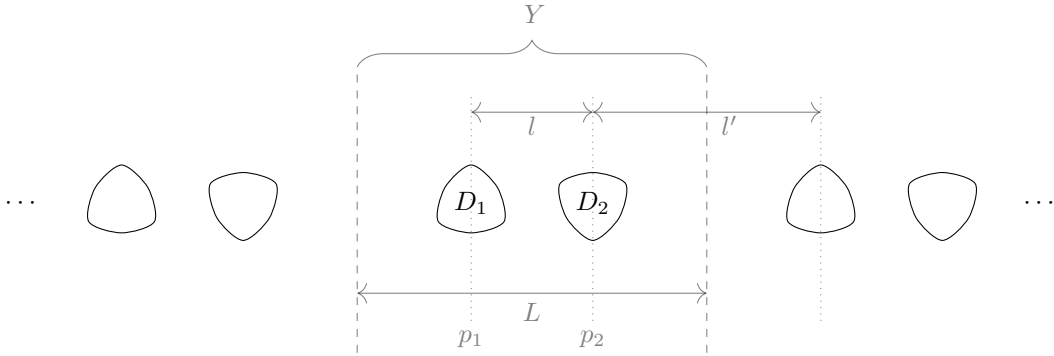

To enable explicit computations of asymptotic expansions and topological properties, we will assume that the resonators are dilute, in the sense that they occupy a small volume compared to the surrounding medium. As in \Cref{sec:exceptional}, we assume that the resonators can be obtained by rescaling fixed domains $B_1, B_2$ as follows:
\begin{equation}\label{eq:dilute_SSH}
	D_1=\epsilon B_1 - \left(\frac{l}{2},0,0\right), \quad  D_2=\epsilon B_2 + \left(\frac{l}{2},0,0\right),
\end{equation}
for some small parameter $\epsilon > 0$. The first result shows that, in addition to the band gap above the second band, there is also a band gap between the first two bands in the case $l\neq L/2$ \cite{ammari2020robust}.
\begin{thm} \label{cor:bandgap}
	In the dilute regime and with $\delta$ sufficiently small, there exists a subwavelength band gap between the first two Bloch band functions $\alpha \mapsto \omega^\alpha_j, j=1,2,$ if $l\neq L/2$, \ie{}
	\begin{equation*}
		\max_{\alpha \in Y^*} \omega_1^\alpha < \min_{\alpha \in Y^*} \omega_2^\alpha,
	\end{equation*}
	for $\epsilon$ and $\delta$ small enough.
\end{thm}
In the case $l=l'=L/2$, the first and the second bands will have a degeneracy at $\alpha=\pi/L$. Moreover, the two bands will intersect linearly in a so-called \emph{Dirac cone}. By altering $l$, this Dirac cone will open a band gap. As we shall see, the nature of this band gap is fundamentally different in the two cases $l<l'$ and $l>l'$.

\subsubsection{Topological indices and  band inversion} \label{sec:bandinv}
Assuming that the dimension of periodicity $d_l=1$, we will now define an index which quantifies the topological properties of the Bloch eigenbundle. Such index is the one-dimensional analogue of Chern numbers in higher periodicity dimensions; see, for instance, \cite{topologychern,topologychern2}.

\begin{defn}[Zak phase] \label{defn:zak}
For a non-degenerate band $\omega_j^\alpha$, we let $u_j^\alpha$ be a family of normalised eigenmodes which depends continuously on $\alpha$. We then define the Zak phase $\varphi_j^\mathrm{zak}$ as
$$\varphi_j^\mathrm{zak} := \i \int_{Y^*} \big\langle u_j^\alpha, \frac{\p }{\p \alpha} u_j^\alpha \big\rangle \dx \alpha,$$
where $\langle \cdot, \cdot \rangle$ denotes the $L^2(D)$-inner product.
\end{defn}
Qualitatively, a non-zero Zak phase means that the crystal has undergone \emph{band inversion}, meaning that at some point in the Brillouin zone the monopole/dipole nature of the first/second Bloch eigenmodes has swapped. In the current setting, a monopole mode is a mode with even parity, while a dipole mode is a mode with odd parity. For the two-resonator chain studied here, the first mode will always be of monopole nature at the origin $\alpha=0$  (which is a consequence of the fact that the eigenvectors of the periodic capacitance matrix $C^0$ are, respectively, $\left(\begin{smallmatrix}1\\1\end{smallmatrix}\right)$ and $\left(\begin{smallmatrix}1\\-1\end{smallmatrix}\right)$). At $\alpha = \pi/L$, the modes will also be of monopole/dipole nature, but which of these corresponds to the first or second mode depends on the geometry of the structure \cite{ammari2020topological}.
\begin{thm} \label{thm:phase}
	We assume that $D$ is in the dilute regime specified by \eqref{eq:dilute_SSH}. Then the Zak phase satisfies
	$$ \varphi_j^z = \begin{cases}
		0, \quad &\text{if} \ \ l < l', \\
		\pi, \quad &\text{if} \ \ l > l',
	\end{cases}$$
	for $\epsilon$ and $\delta$ small enough.
\end{thm}
The Zak phase predicts band inversion in the case $l>l'$, which is the statement of the next result.
\begin{prop} \label{prop:bandinv}
	For $\epsilon$ and $\delta$ small enough, the crystal has undergone band inversion in the case $l>l'$, but not if $l<l'$. In other words, the eigenfunctions associated with the first and second bands at $\alpha = \pi/L$ satisfy
	$$u_1^{\pi/L}(-x) = u_1^{\pi/L}(x), \qquad u_2^{\pi/L}(-x) = - u_2^{\pi/L}(x),\qquad \text{ when $l<l'$,}$$
	%while
	$$u_1^{\pi/L}(-x) = -u_1^{\pi/L}(x), \qquad u_2^{\pi/L}(-x) = u_2^{\pi/L}(x),\qquad \text{ when $l>l'$}.$$
\end{prop}
	As $\delta \to 0$, these eigenfunctions are given, respectively, by
	\begin{align*}
	u_1^{\pi/L}(x) &=
		\left(\begin{smallmatrix}1\\1\end{smallmatrix}\right)\cdot\textbf{S}_D^{\pi/L,\omega}(x)+O(\delta), & u_2^{\pi/L}(x)&=
		\left(\begin{smallmatrix}1\\-1\end{smallmatrix}\right)\cdot\textbf{S}_D^{\pi/L,\omega}(x)+O(\delta),& \text{ when $l<l'$,} \\[0.3em]
u_1^{\pi/L}(x)&=
\left(\begin{smallmatrix}1\\-1\end{smallmatrix}\right)\cdot\textbf{S}_D^{\pi/L,\omega}(x)+O(\delta), & u_2^{\pi/L}(x)&=
\left(\begin{smallmatrix}1\\1\end{smallmatrix}\right)\cdot\textbf{S}_D^{\pi/L,\omega}(x)+O(\delta),& \text{ when $l<l'$.}
\end{align*}
Here, $\textbf{S}_D^{\alpha,\omega}(x)$ is the function defined in \Cref{prop:eigenvector_quasi}.

\subsubsection{Robustness of edge modes}
We now study a finite chain of resonators which supports topologically protected edge modes. Specifically, we assume that $D$ has the form
\begin{equation} \label{finite_form}
	D = \left(\bigcup_{n=-M}^{M} D_0 + n(l+l',0,0)  \right) \bigcup \left( \bigcup_{n=-M+1}^M D_0 + n(l+l',0,0) - (l',0,0)\right),
\end{equation}
where $D_0$ is a single repeating resonator. In other words, $D$ consists of an odd number  $N$ of identical resonators ($N = 4M+1$) with alternating distances $l$ and $l'$ that are swapped at the middle resonator. An example of such a configuration is depicted in \Cref{fig:finite_SSH}. In this figure, it is shown how to associate different Zak phases on either side of the central resonator (which constitutes the ``edge''). Based on the principle of \emph{bulk-boundary correspondence}, we thereby expect robust localised modes around this edge.

\begin{figure}[p]
	\centering
	\begin{tikzpicture}[scale=1.1]
		\begin{scope}
			\draw (0.65,0) coordinate (start1) circle (8pt);
			%	\shade[ball color = gray!50, opacity = 0.4](0.65,0) circle (8pt);
			%\draw[fill] (0.2,0) circle (0.5pt) node[yshift=-8pt]{$-x_1$};
			%\draw[fill] (1.1,0) circle (0.5pt) node[yshift=-8pt]{$x_1$};
			\draw[<->, opacity=0.5] (0.65,0) -- (2.05,0) node[pos=0.5, yshift=-7pt,]{$l'$};
			\draw[<->, opacity=0.5] (-0.75,0) -- (0.65,0) node[pos=0.5, yshift=-7pt,]{$l'$};
			%\draw[dashed, opacity=0.5] (0.2,-0.2) -- (0.6,-0.2) coordinate (end);
			%\pic [draw, ->, "$\theta$", angle eccentricity=1.5] {angle = end--start1--start2};
		\end{scope}

		\begin{scope}[xshift=-1.85cm]
			\draw (0.2,0) circle (8pt);
			\draw[<->, opacity=0.5] (0.2,0) -- (1.1,0) node[pos=0.5, yshift=-7pt,]{$l$};
			\begin{scope}[xshift = 1.3cm]
				\draw (-0.2,0) circle (8pt);
			\end{scope};
		\end{scope};

		\begin{scope}[xshift=1.85cm]
			\draw[<->, opacity=0.5] (0.2,0) -- (1.1,0) node[pos=0.5, yshift=-7pt,]{$l$};
			\draw (0.2,0) circle (8pt);
			%\draw[fill] (0.2,0) circle (0.5pt);
			\begin{scope}[xshift = 1.3cm]
				\draw (-0.2,0) circle (8pt);
			\end{scope};
		\end{scope};

		\begin{scope}[xshift=4.05cm]
			\draw (0.2,0) circle (8pt);
			%\draw[<->, opacity=0.5] (0.2,0) -- (1.1,0) node[pos=0.5, yshift=-5pt,]{$d$};
			%	\draw[opacity=0.5] (0.65,1) -- (0.65,-1);
			\begin{scope}[xshift = 1.3cm]
				\draw (-0.2,0) circle (8pt);
			\end{scope};
		\end{scope};

		\begin{scope}[xshift=6.25cm]
			\draw (0.2,0) circle (8pt);
			%	\draw[opacity=0.5] (0.65,1) -- (0.65,-1) node[above left]{$\phi^z = \pi$};
			%\draw[<->, opacity=0.5] (0.2,0) -- (1.1,0) node[pos=0.5, yshift=-5pt,]{$d$};
			\begin{scope}[xshift = 1.3cm]
				\draw (-0.2,0) circle (8pt);
			\end{scope};
		\end{scope};

		\begin{scope}[xshift=-4.05cm]
			\draw (0.2,0) circle (8pt);
			%\draw[<->, opacity=0.5] (0.2,0) -- (1.1,0) node[pos=0.5, yshift=-5pt,]{$d$};
			\begin{scope}[xshift = 1.3cm]
				\draw (-0.2,0) circle (8pt);
			\end{scope};
		\end{scope};

		\begin{scope}[xshift=-6.25cm]
			\draw (0.2,0) circle (8pt);
			%	\draw[opacity=0.5] (-0.5,1) -- (-0.5,-1);
			%	\draw[opacity=0.5]  (1.8,1) -- (1.8,-1) node[above left]{$\phi^z = 0$};
			%\draw[<->, opacity=0.5] (0.2,0) -- (1.1,0) node[pos=0.5, yshift=-5pt,]{$d$};
			\begin{scope}[xshift = 1.3cm]
				\draw (-0.2,0) circle (8pt);
			\end{scope};
		\end{scope};

		%	\begin{scope}[xshift=-6.25cm]
		%	\draw [decorate,opacity=0.5,decoration={brace,amplitude=10pt},xshift=-4pt,yshift=0pt]
		%	(-0.25,0.5) -- (1.8,0.5) node [black,midway,xshift=-0.6cm]{};
		%	\end{scope}

		\begin{scope}[xshift=-4.05cm]
			\draw [decorate,opacity=0.5,decoration={brace,amplitude=10pt},xshift=-4pt,yshift=0pt]
			(-0.25,0.5) -- (1.8,0.5) node [black,midway,xshift=-0.6cm]{};
		\end{scope}

		\begin{scope}[xshift=-1.85cm]
			\draw [decorate,opacity=0.5,decoration={brace,amplitude=10pt},xshift=-4pt,yshift=0pt]
			(-0.25,0.5) -- (1.8,0.5) node [black,midway,xshift=-0.6cm]{};
		\end{scope}

		\begin{scope}[xshift=0.72cm]
			\draw [decorate,opacity=0.5,decoration={brace,amplitude=10pt},xshift=-4pt,yshift=0pt]
			(-0.25,0.5) -- (1.8,0.5) node [black,midway,xshift=-0.6cm]{};
		\end{scope}

		\begin{scope}[xshift=2.92cm]
			\draw [decorate,opacity=0.5,decoration={brace,amplitude=10pt},xshift=-4pt,yshift=0pt]
			(-0.25,0.5) -- (1.8,0.5) node [black,midway,xshift=-0.6cm]{};
		\end{scope}

		%	\begin{scope}[xshift=5.12cm]
		%	\draw [decorate,opacity=0.5,decoration={brace,amplitude=10pt},xshift=-4pt,yshift=0pt]
		%	(-0.25,0.5) -- (1.8,0.5) node [black,midway,xshift=-0.6cm]{};
		%	\end{scope}

		\begin{scope}[xshift=-4.05cm-4pt]
			\draw[opacity=0.5]    (0.775,0.9) to[out=90,in=-100] (1.54,1.5);
		\end{scope}

		\begin{scope}[xshift=-1.85cm-4pt]
			\draw[opacity=0.5]    (0.775,0.9) to[out=100,in=-80] (-0.64,1.5);
		\end{scope}

		%	\begin{scope}[xshift=-1.85cm-4pt]
		%	\draw[opacity=0.5]    (0.775,0.9) to[out=100,in=-45] (-2.25,1.5);
		%	\end{scope}

		\begin{scope}[xshift=0.72cm-4pt]
			\draw[opacity=0.5]    (0.775,0.9) to[out=90,in=-100] (2.04,1.5);
		\end{scope}

		\begin{scope}[xshift=2.92cm-4pt]
			\draw[opacity=0.5]    (0.775,0.9) to[out=100,in=-80] (-0.14,1.5);
		\end{scope}

		\begin{scope}[xshift=-2.4cm]
			\draw [decorate,opacity=0.5,decoration={brace,amplitude=10pt,mirror},xshift=-4pt,yshift=0pt]
			(-1.125,1.9) -- (0.925,1.9) node [black,midway,xshift=-0.6cm]{};
			\begin{scope}[xshift=-4pt]
				\draw[opacity=0.5, dotted]    (-1.25,2) -- (-1.25,2.7);
				\draw[opacity=0.5, dotted]    (1.05,2) -- (1.05,2.7);
			\end{scope}
			\begin{scope}[xshift = -4pt-0.75cm]
				\draw[opacity=0.5] (0.2,2.3) circle (8pt);
				\node[opacity=0.5] at (0.65,2.9) {$\varphi_j^z = 0$};
				%	\draw[opacity=0.5] (0.65,1) -- (0.65,-1) node[above left]{$\phi^z = \pi$};
				%\draw[<->, opacity=0.5] (0.2,0) -- (1.1,0) node[pos=0.5, yshift=-5pt,]{$d$};
				\begin{scope}[xshift = 1.3cm]
					\draw[opacity=0.5] (-0.2,2.3) circle (8pt);
				\end{scope}
			\end{scope}
		\end{scope}

		\begin{scope}[xshift=2.87cm]
			\draw [decorate,opacity=0.5,decoration={brace,amplitude=10pt,mirror},xshift=-4pt,yshift=0pt]
			(-1.125,1.9) -- (0.925,1.9) node [black,midway,xshift=-0.6cm]{};
			\begin{scope}[xshift=-4pt]
				\draw[opacity=0.5, dotted]    (-1.25,2) -- (-1.25,2.7);
				\draw[opacity=0.5, dotted]    (1.05,2) -- (1.05,2.7);
			\end{scope}
			\begin{scope}[xshift = -4pt-0.75cm]
				\draw[opacity=0.5] (-0.05,2.3) circle (8pt);
				\node[opacity=0.5] at (0.65,2.9) {$\varphi_j^z = \pi$};
				%	\draw[opacity=0.5] (0.65,1) -- (0.65,-1) node[above left]{$\phi^z = \pi$};
				%\draw[<->, opacity=0.5] (0.2,0) -- (1.1,0) node[pos=0.5, yshift=-5pt,]{$d$};
				\begin{scope}[xshift = 1.3cm]
					\draw[opacity=0.5] (0.05,2.3) circle (8pt);
				\end{scope}
			\end{scope}
		\end{scope}
	\end{tikzpicture}
	\caption{Two-dimensional cross-section of a finite dimer chain with 13 resonators, heuristically showing how to identify unit cells with different Zak phases on either side of the edge.} \label{fig:finite_SSH}
\end{figure}
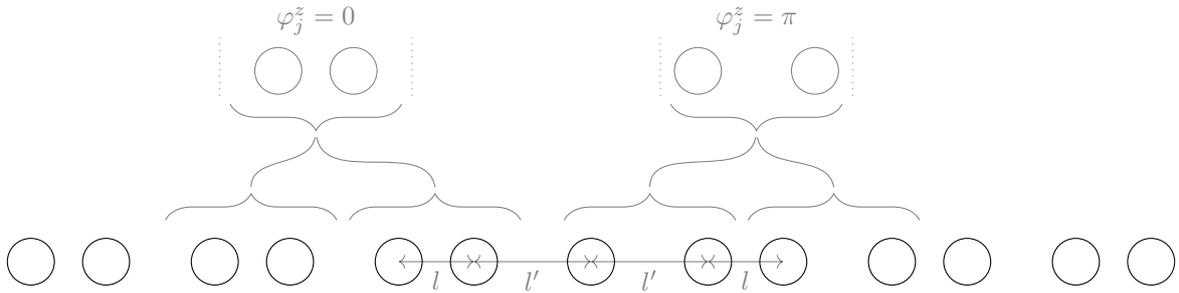

We model wave propagation in the crystal $D$ by the Helmholtz problem
\begin{equation} \label{eq:scattering_finite}
	\left\{
	\begin{array} {ll}
		\ds \Delta {u}+ \omega^2 {u}  = 0 \quad &\text{in } \R^3 \setminus \p D, \\
		\nm
		\ds  {u}|_{+} -{u}|_{-}  =0  \quad &\text{on } \partial D, \\
		\nm
		\ds  \delta \frac{\partial {u}}{\partial \nu} \bigg|_{+} - \frac{\partial {u}}{\partial \nu} \bigg|_{-} =0 \quad& \text{on } \partial D, \\
		\nm
		\ds |x| \left(\tfrac{\p}{\p|x|}-\i \omega\right)u \to 0
		&\text{as } {|x|} \rightarrow \infty.
	\end{array}
	\right.
\end{equation}
\Cref{fig:perturbation1} shows a comparison between the band-gap frequencies of a chain with a topological defect (as in \Cref{fig:finite_SSH}), compared against a chain with a point defect (analogously as in \Cref{sec:defects}).  Not only is the topological band-gap frequency further away from the edges of the band gap, but it also exhibits a lower variance when random Gaussian errors is imposed on the resonator locations. For large error, we see that the band-gap frequency of the topological defect chain (\Cref{fig:nondilute}) is much more robust than the band-gap frequency of the point-defect chain.

In \Cref{rmk:tightbind} we observed that, under the dilute assumption, the capacitance formulation is analogous to the tight-binding model commonly utilised in studies of quantum-mechanical systems. The tight-binding model is often coupled with a nearest-neighbour approximation, whereby long-range interactions are neglected. \Cref{fig:tight_bind} shows the bulk and band-gap frequencies computed using such nearest-neighbour approximation. Compared to \Cref{fig:dilute}, where all interactions are taken into account, we observe that the nearest-neighbour approximation is not accurate. This discrepancy shows a fundamental difference between topological edge modes in the setting of classical waves and quantum-mechanical waves. Moreover, we see that the band-gap frequency in \Cref{fig:tight_bind} is unaffected by the error. This can be seen as a consequence of chiral symmetry, which is not present without the nearest-neighbour approximation.

\begin{figure}[p]
	\begin{center}
		\begin{subfigure}[b]{0.45\linewidth}
			\includegraphics[height=5.0cm]{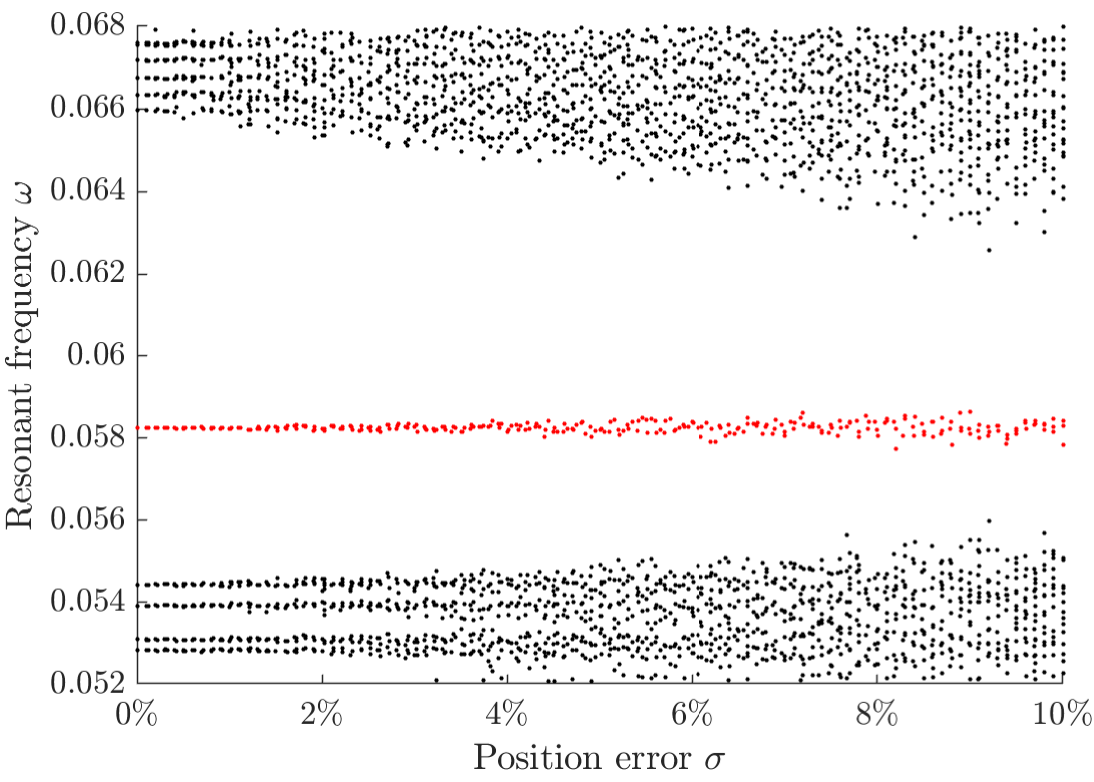}
			\caption{Dimer chain with 41 resonators, separation distances $d=3, d' = 6$.} \label{fig:nondilute}
		\end{subfigure}
		\hspace{10pt}
		\begin{subfigure}[b]{0.45\linewidth}
			\includegraphics[height=4.1cm]{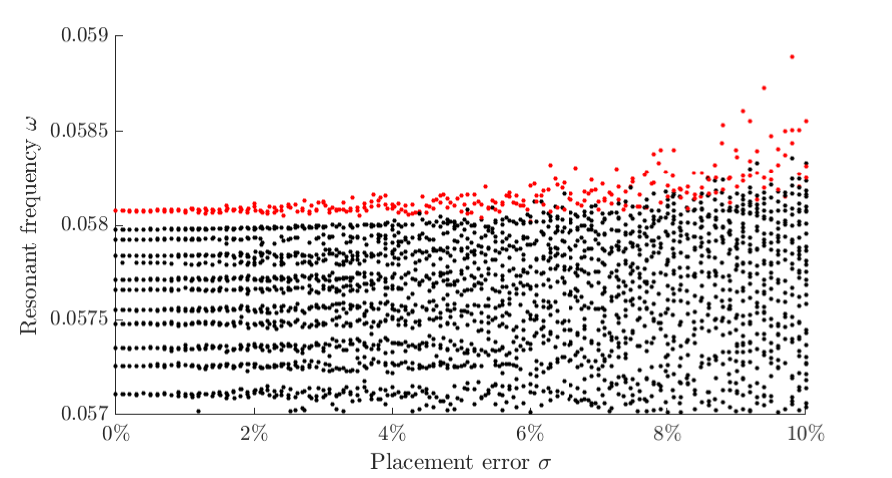}
			\caption{Point-defect chain with 41 resonators, separation distance $d=12$ and defect radius $R_d = 0.99$.} \label{fig:point_defect_errors}
		\end{subfigure}
		\caption{Simulation of band-gap frequency (red) and bulk frequencies (black) of a topological defect (a) and point defect (b) chain, with Gaussian $\mathcal{N}(0,\sigma^2)$ errors added to the resonator positions. The standard deviation $\sigma$ is expressed as a percentage of the average resonator separation. %In all cases, the resonator radius was $R=1$.
		}
		\label{fig:perturbation1}
	\end{center}
\end{figure}

\begin{figure}[p]
	\begin{center}
		\begin{subfigure}[b]{0.45\linewidth}
			\includegraphics[height=5.0cm]{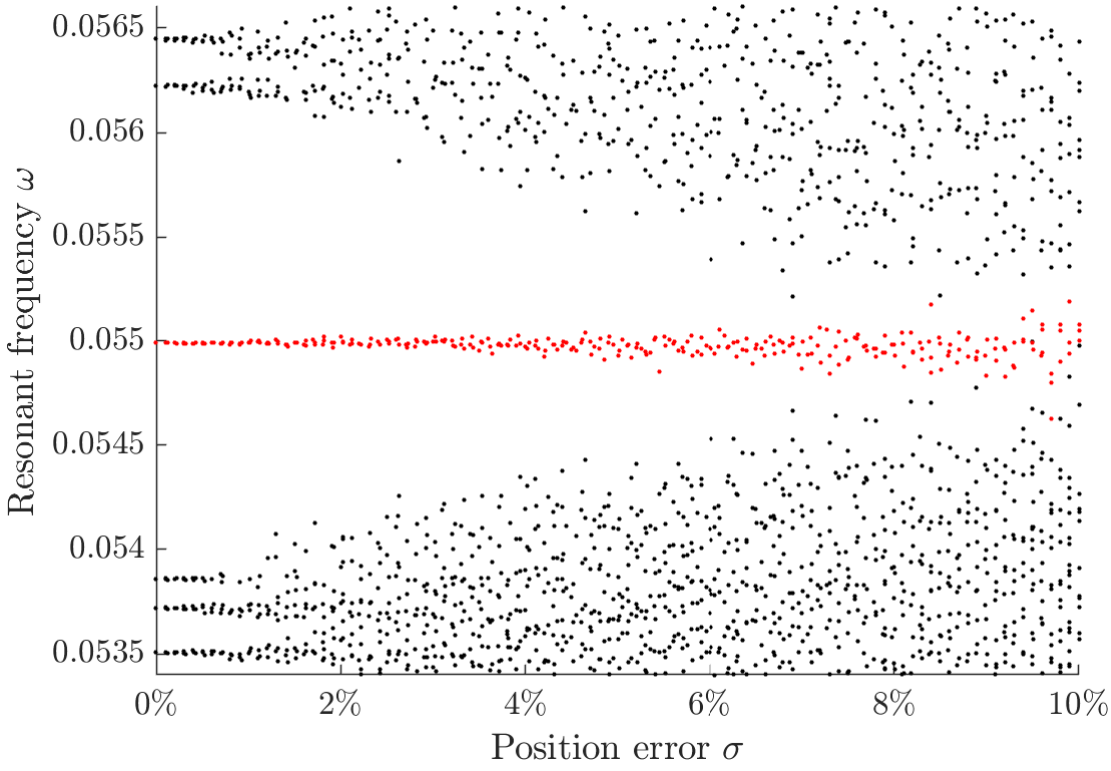}
			\caption{Dilute dimer chain with 41 resonators, separation distances $d=12, d' = 42$.} \label{fig:dilute}
		\end{subfigure}
		\hspace{10pt}
		\begin{subfigure}[b]{0.45\linewidth}
			\includegraphics[height=5.0cm]{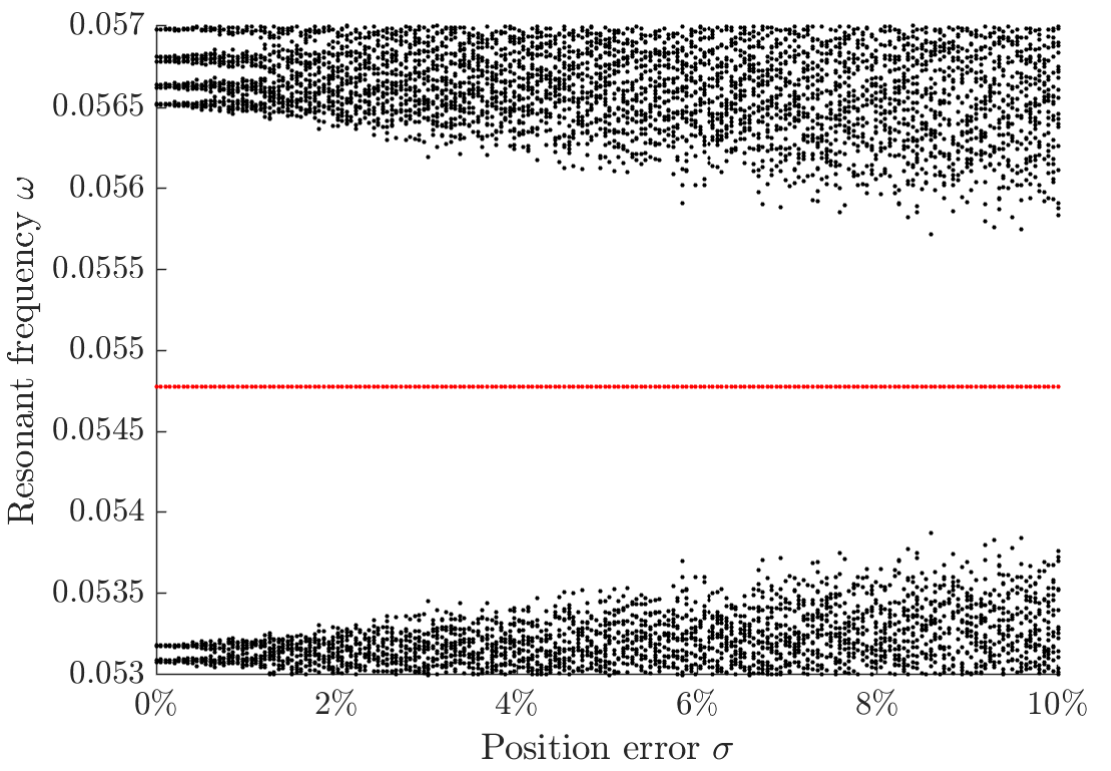}
			\caption{Nearest-neighbour approximation for the dilute dimer chain from (a).% with 41 resonators, separation distances $d=12, d' = 42$.
			} \label{fig:tight_bind}
		\end{subfigure}
		\caption{Similar simulations as in \Cref{fig:perturbation1}, but for a dilute chain. (a): fully-coupled simulations, where all interactions are taken into account. (b): nearest-neighbour approximation, where long-range interactions are neglected. The nearest-neighbour approximation provides a chirally symmetric problem, in which the centre frequency is preserved. This approximation is not accurate, and the centre frequency will have a small but non-zero variance corresponding to small but non-zero long-range interactions.
		}
		\label{fig:perturbation2}
	\end{center}
\end{figure}

\subsubsection{Edge modes in a dislocated chain}
In \Cref{sec:bandinv} we saw that the structure under consideration has different topological properties in the different cases $l>l'$ and $l<l'$. These two regimes only differ by a choice of unit cell: if we shift the unit cell $Y$ by $L/2$ the regimes will be swapped. Correspondingly, edges between the different topological regimes are created by introducing a shift between the unit cells on either side of the edge. This is exemplified in \Cref{fig:finite_SSH}. Another example of this principle is to translate, or \emph{dislocate}, half the chain by some distance $d$ as illustrated in \Cref{fig:infinite}. In this case, we consider a chain of resonators given by
$$\Dc_\mathrm{disloc} = \Bigg(\bigcup_{m=-\infty}^{-1} D+(mL,0,0)\Bigg)\cup\Bigg(\bigcup_{m=0}^\infty D+(mL+d,0,0)\Bigg),$$
where $D$ is the resonator dimer defined as in the beginning of this section.

Qualitatively, it is straightforward to understand how the dislocation will affect the existence and behaviour of band-gap frequencies. As $d$ increases from $0$, we will detune the dimer coupling, which is responsible for the band gap, and band-gap frequencies will therefore appear from each edge of the band gap. On the other hand, when $d$ is very large the two half-space chains will decouple, and the bulk-boundary correspondence suggests that each of these half-chains will support a single band-gap frequency. We therefore expect only a single band-gap frequency in this limit. As $d$ varies between $0$ and $\infty$ these two frequencies hybridize, and will together cover the whole band gap. This is sketched in \Cref{fig:gapcrossing}.

We now set out to prove that the picture outlined above is correct. Similarly to \Cref{sec:defects}, we will apply a fictitious source method to model the dislocated structure $\Dc_\mathrm{disloc}$ in terms of the original structure $\Dc$ along with fictitious sources on the dislocated resonators. The analysis is separated into three parts, depending on the value of $d$:
\begin{itemize}
	\item $d\ll1$. In this case, we can use asymptotic expansions in terms of $d$ to prove that there is a band-gap frequency emerging from each edge of the band gap.
	\item $d=mL$ for some $m>0$. This particular dislocation is equivalent to removing $m$ dimers from $\Dc$. This observation simplifies the problem and allows explicit computations of the band-gap frequencies in terms of the eigenvalue problem of certain Toeplitz matrices.
	\item $d>d_0$, where $d_0$ is the width of one resonator. Due to technical reasons, we restrict the remaining values of $d$, so that the dislocated and the original resonators are not overlapping.
\end{itemize}
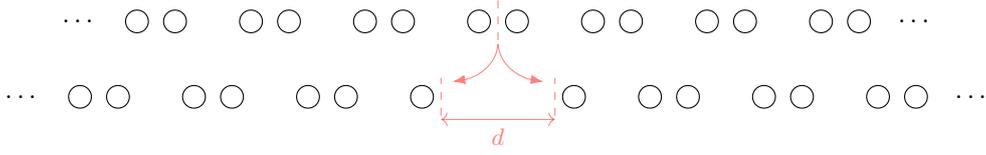
\begin{figure}
	\centering
	\begin{tikzpicture}[scale=0.5]
		\begin{scope}[shift={(1.5,2)}]
			\draw (-6,0) circle (0.3);
			\draw (-5,0) circle (0.3);
			\draw (-3,0) circle (0.3);
			\draw (-2,0) circle (0.3);
			\draw (0,0) circle (0.3);
			\draw (1,0) circle (0.3);
			\draw (3,0) circle (0.3);
			\draw (4,0) circle (0.3);
			\draw (6,0) circle (0.3);
			\draw (7,0) circle (0.3);
			\draw (9,0) circle (0.3);
			\draw (10,0) circle (0.3);
			\draw (12,0) circle (0.3);
			\draw (13,0) circle (0.3);
			\node at (14.5,0) {$\dots$};
			\node at (-7.5,0) {$\dots$};
			\draw[dashed,opacity=0.5,red] (3.5,-0.5) -- (3.5,0.6);
			\draw[-Latex,opacity=0.5,red] (3.5,-0.6) to[out=-95,in=10] (2.3,-1.6);
			\draw[-Latex,opacity=0.5,red] (3.5,-0.6) to[out=-85,in=170] (4.7,-1.6);
		\end{scope}
		\draw (-6,0) circle (0.3);
		\draw (-5,0) circle (0.3);
		\draw (-3,0) circle (0.3);
		\draw (-2,0) circle (0.3);
		\draw (0,0) circle (0.3);
		\draw (1,0) circle (0.3);
		\draw (3,0) circle (0.3);
		\draw[dashed,opacity=0.5,red] (3.5,-0.5) -- (3.5,0.5);
		\draw[dashed,opacity=0.5,red] (6.5,-0.5) -- (6.5,0.5);
		\draw (7,0) circle (0.3);
		\draw (9,0) circle (0.3);
		\draw (10,0) circle (0.3);
		\draw (12,0) circle (0.3);
		\draw (13,0) circle (0.3);
		\draw (15,0) circle (0.3);
		\draw (16,0) circle (0.3);
		\node at (17.5,0) {$\dots$};
		\node at (-7.5,0) {$\dots$};
		\draw[<->,opacity=0.5,red] (3.5,-0.6) -- (6.5,-0.6) node[pos=0.5,below]{\small $d$};
		%	\draw[dotted] (-6.5,-1.5) rectangle (3.5,1.5);
		%	\draw[dotted] (6.5,-1.5) rectangle (16.5,1.5);
		% 	\draw[<->] (3,0.5) -- (7,0.5) node[pos=0.5, yshift=7pt]{$l+d$};
		% 	\draw[<->] (0,0.5) -- (1,0.5) node[pos=0.5, yshift=7pt]{$l$};
		% 	\draw[<->] (0,-0.5) -- (3,-0.5) node[pos=0.5, yshift=-7pt]{$L$};
		% 	\draw[<->] (9,0.5) -- (10,0.5) node[pos=0.5, yshift=7pt]{$l$};
		% 	\draw[<->] (7,-0.5) -- (10,-0.5) node[pos=0.5, yshift=-7pt]{$L$};
		%	\draw[<->, opacity=0.5] (0,-0.8) -- (5,-0.8) node[pos=0.5, yshift=-7pt,]{$L$};
	\end{tikzpicture}
	\caption{We start with an array of pairs of subwavelength resonators, known to have a subwavelength band gap. A dislocation (with size $d>0$) is introduced to create band-gap frequencies.} \label{fig:infinite}
\end{figure}

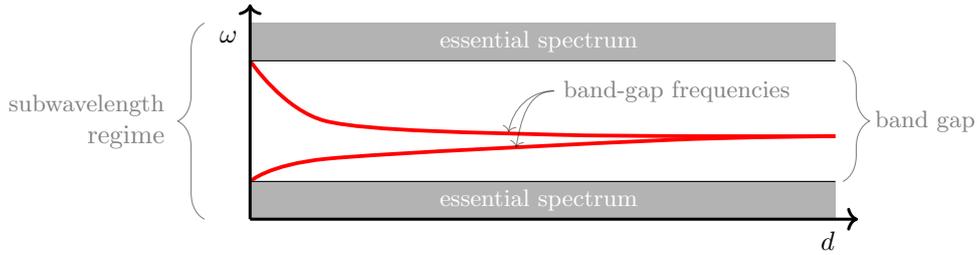
\begin{figure}[htb]
	\centering
	\begin{tikzpicture}
		\draw[line width=0.5mm,red] plot [smooth] coordinates {(0,0.5) (1,0.8) (5.5,1.05) (7.7,1.1)};
		\draw[line width=0.5mm,red] plot [smooth] coordinates {(0,2.1) (1,1.3) (4,1.12) (7.7,1.1)};
		\draw[gray!60!white,fill=gray!60!white] (0,0) -- (7.7,0) -- (7.7,0.5) --  (0,0.5);
		% 	\draw plot [smooth] coordinates {(7.7,0.5) (4,0.55) (2,0.5) (0,0.5)};
		\draw (0,0.5) -- (7.7,0.5);
		\draw[gray!60!white,fill=gray!60!white] (0,2.6) -- (7.7,2.6) -- (7.7,2.1) -- (0,2.1);
		% 	\draw plot [smooth] coordinates {(7.7,2.1) (4,2.1) (2,2.05) (0,2)} -- (0,2.5);
		\draw (7.7,2.1) -- (0,2.1);
		\draw[->,line width=0.4mm] (0,0) -- (8,0) node[pos=0.95, yshift=-8pt]{$d$};
		\draw[->,line width=0.4mm] (0,0) -- (0,2.85) node[pos=0.85, xshift=-8pt]{$\omega$};
		\draw[<-,opacity=0.5] (3.5,0.95) to[out=80,in=180] (4,1.7);
		\node[right,opacity=0.5] at (4,1.7) {\small band-gap frequencies};
		\draw[<-,opacity=0.5] (3.4,1.15) to[out=80,in=180] (4,1.7);
		\draw [decorate,opacity=0.5,decoration={brace,amplitude=10pt}]
		(7.8,2.1) -- (7.8,0.5) node [midway,right,xshift=3mm]{\small band gap};
		\node at (3.8,2.35){\color{white}\small essential spectrum};
		\node at (3.8,0.25){\color{white}\small essential spectrum};
		\draw [decorate,opacity=0.5,decoration={brace,amplitude=10pt}]
		(-0.6,0) -- (-0.6,2.6) node [midway,xshift=-4mm,align=right,left]{\small subwavelength\\ regime};
	\end{tikzpicture}
	\caption{As the dislocation size $d$ increases from zero, a band-gap frequency appears from each edge of the subwavelength band gap. These two frequencies converge to a single value within the subwavelength band gap as $d\to\infty$.} \label{fig:gapcrossing}
\end{figure}

With these ideas at hand, the following two theorems were proved in \cite{ammari2020robust}. The first result, valid for small $d$, shows the emergence of a band-gap frequency from each edge of the band gap.
\begin{thm}\label{prop:smalld}
	Assume that $D_1$ and $D_2$ are strictly convex. For small enough $d$ and $\delta$, and in the case $l_0 > 1/2$, there are two band-gap frequencies $\omega_1(d), \omega_2(d)$ such that $\omega_j(d) \rightarrow \omega_j^\diamond, j=1,2$ as $d\rightarrow 0$. In the case $l_0 < 1/2$, there are no band-gap frequencies as $d, \delta \rightarrow 0$.
\end{thm}
\begin{thm} \label{thm:main}
	Assume that the resonators are in the dilute regime specified by \eqref{eq:dilute_SSH} and that $l_0>1/2$. Then, for small enough $\delta$ and $\epsilon$, there exists some $d_0=O(\epsilon)$ such that there are two band-gap frequencies $\omega_1(d)$ and $\omega_2(d)$
	for all $d\in [d_0,\infty)$, both of which converge to the same value $\omega_\infty$ as $d\rightarrow \infty$.
\end{thm}
In particular, \Cref{thm:main} states that the band-gap frequencies will cover an interval $\mathcal{I}:=[\omega_1(d_0), \omega_2(d_0)]$ inside the band gap, and therefore allows us to fine-tune the system to achieve optimal robustness. In \Cref{fig:eigenmodes} we see numerically computed eigenmodes, clearly showing how the edge modes of the ``half system'' hybridize to create the two band-gap frequencies.

\begin{figure}[htb]
		\begin{minipage}[b]{0.48\linewidth}
			\vspace{0pt}
			\centering
			\includegraphics[width=\linewidth,trim={1cm 1cm 1cm 0},clip]{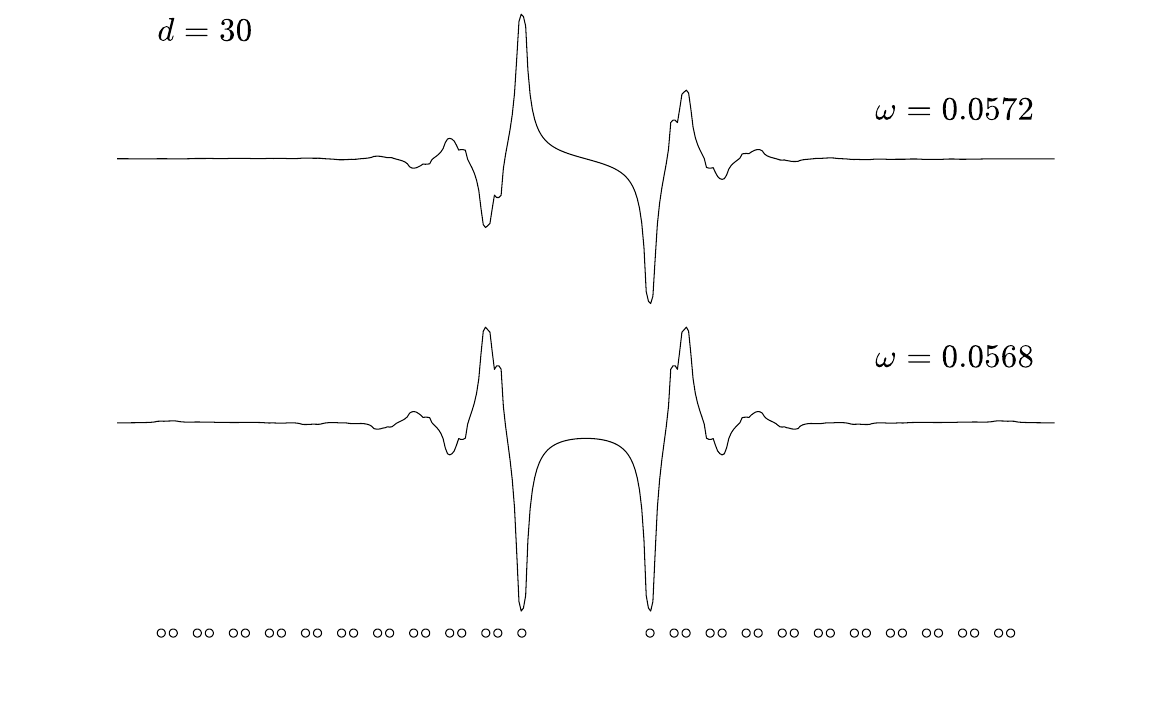}
		\end{minipage}
	\begin{minipage}[b]{0.48\linewidth}
		\centering
		\includegraphics[width=\linewidth,trim={1cm 1cm 1cm 5cm},clip]{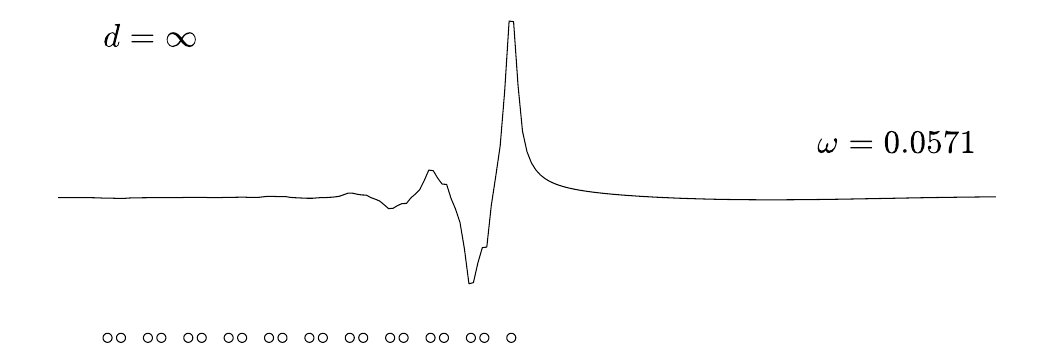}
		\vspace{10pt}
	\end{minipage}
	\caption{\textit{Left:} The two edge modes for an array of 42 spherical resonators of radius $1$. Here, we simulate an array with parameters $L=9$, $l=6$, $d=30$ and $\delta=1/7000$. \textit{Right:} For comparison, the edge mode of the corresponding `half system' is shown, which can be thought of as the $d=\infty$ case. In both cases, the eigenmodes are shown directly above the corresponding system of resonators.} \label{fig:eigenmodes}
\end{figure}

\subsubsection{Non-Hermitian band inversion and edge modes}\label{sec:pt_top}
In the previous sections, we assumed that the wave speeds satisfy $v=v_1=v_2=1$, and only considered topological phenomena originating from the geometry of the resonator chain. In this section we do the opposite, namely we study topologically protected modes in structures where the geometry is periodic, and where topological edges are introduced in terms of the resonator wave speeds. Most importantly, we allow the wave speeds to be complex, giving a non-Hermitian capacitance formulation similarly as in \Cref{sec:exceptional}. Details of this analysis are found in \cite{ammari2020edge}.

We let $D$ be the dimer $D=D_1\cup D_2$ as defined previously in \Cref{sec:bandinv}. For $i=1,2$ and $m\in \Z$, we introduce the notation
$$D_i^m := D_i + (mL,0,0), \qquad k_i^m := \frac{\omega}{v_i^m}.$$
Here, $v_i^m$ denotes the wave speed inside $D_i^m$. We will assume that the contrast parameter is given by $\delta\in \R$ for all resonators. We then consider
\begin{equation} \label{eq:scattering_pttop}
	\left\{
	\begin{array} {ll}
		\ds \Delta {u}+ k^2 {u}  = 0 & \text{in } \R^3 \setminus \Dc, \\
		\nm
		\ds \Delta {u}+ (k_i^m)^2 {u}  = 0 & \text{in } D_i^m, \\
		\nm
		\ds  {u}|_{+} -{u}|_{-}  =0  & \text{on } \p \Dc, \\
		\nm
		\ds  \delta \frac{\partial {u}}{\partial \nu} \bigg|_{+} - \frac{\partial {u}}{\partial \nu} \bigg|_{-} =0 & \text{on } \p \Dc, \\
		\nm
		\ds u(x_1,x_2,x_3) & \text{satisfies the outgoing radiation condition as } \sqrt{x_2^2+x_3^2} \rightarrow \infty.
	\end{array}
	\right.
\end{equation}
In the general setting, this equation is not periodic and cannot be approached using Floquet-Bloch theory. For prescribed values $v_1,v_2\in \CC$ we will consider two different cases, depending on $v_i^m$:
\begin{align}\label{eq:vperiodic}
	\text{Periodic structure:} &  \qquad v_1^m = v_1, \ m\in \Z, & &v_2^m = v_2,  \ m\in \Z,
\\ \label{eq:vdefect}
	\text{Defect structure:} & \qquad v_1^m = \begin{cases}v_1, & m \leq 0, \\ v_2, &m > 0,\end{cases} & &v_2^m = \begin{cases}v_2, & m \leq 0, \\ v_1, & m > 0.\end{cases}
\end{align}
We emphasize that complex values of $v_i$ correspond to non-Hermitian structures. In the case $v_1 = \overline{v_2}$, the periodic structure is $\P\T$-symmetric.

We begin by studying the periodic structure. In \Cref{sec:bandinv} we described the concept of band inversion in the Hermitian case. In the present case, since the quasiperiodic generalised capacitance matrix $\C^\alpha$ is non-Hermitian, the band inversion phenomenon is formulated in a slightly different manner. For $j = 1,2,$ we let $\mathbf{u}_j$ and $\mathbf{w}_j$ denote a bi-orthogonal system of eigenvectors of $\C^\alpha$. In other words, $\mathbf{u}_j$ and $\mathbf{w}_j$ are eigenvectors of $\C^\alpha$ and $\left(\C^\alpha\right)^*$, respectively, and satisfy $\mathbf{w}_i\cdot \mathbf{u}_j = \delta_{ij}$.

We let the right Bloch eigenmodes $u_j^\alpha$  be the modes of \eqref{eq:scattering_pttop}, and the left Bloch eigenmodes  $w_j^\alpha$ be the modes of \eqref{eq:scattering_pttop} with wave speeds given by $\overline{v_1}$ and $\overline{v_2}$, respectively. From \Cref{prop:eigenvector_quasi}, we then have
	\begin{align*}
		u_j^\alpha &= \mathbf{u}_{j}\cdot  \textbf{S}_D^{\alpha,k}(x) + O(\delta^{1/2}), \\
		w_j^\alpha &= \mathbf{w}_{j}\cdot\textbf{S}_D^{\alpha,k}(x) + O(\delta^{1/2}).
	\end{align*}

\begin{defn}[non-Hermitian Zak phase] \label{defn:zakNH}
	For a non-degenerate band $\omega_j^\alpha$ of the periodic structure \eqref{eq:vperiodic}, we let $u_j^\alpha$ and $w_j^\alpha$ be a family of normalised right, respectively left, eigenmodes which depend continuously on $\alpha$.
	We define the (non-Hermitian) Zak phase $\varphi_j^{\mathrm{zak}}$ by
	\begin{align*}
		\varphi_j^{\mathrm{zak}} &:= \frac{\i}{2} \int_{Y^*}\left( \Big\langle w_j^\alpha, \frac{\p u_j^\alpha  }{\p \alpha}\Big\rangle + \Big\langle u_j^\alpha, \frac{\p  w_j^\alpha}{\p \alpha} \Big\rangle \right) \dx \alpha.
	\end{align*}
\end{defn}
In the case $v_1 = v_2\in \R$, this definition coincides with \Cref{defn:zak}, and we therefore choose the same notation for these two definitions. In the sequel, we will occasionally write $\varphi_j^{\mathrm{zak}}(v_1, v_2)$ to denote the Zak phase defined with wave speed $v_1$ inside $D_1$ and $v_2$ inside $D_2$.

In the Hermitian case, a non-zero Zak phase is equivalent to an inverted band structure. The fact that the Hermitian Zak phase is quantized originates from the fact that the eigenmodes are purely monopole and dipole modes at $\alpha = 0$ and $\alpha = \pi/L$. Unlike the Hermitian case, the non-Hermitian Zak phase is not quantized. A non-integer value of the Zak phase can be attributed to a ``partial'' band inversion, where the eigenmodes are expressed as (complex) linear combinations of monopole and dipole modes, which swap as $\alpha$ traverses $Y^*$.

For the next result, we will assume that the Hermitian counterpart of the structure is topologically trivial. More precisely, we assume
\begin{equation}\label{eq:trivial}
	\varphi_j^{\mathrm{zak}}\big(\Re(v_1), \Re(v_2)\big) = 0.
\end{equation}
Then the following result from \cite{ammari2020edge} holds.
\begin{prop} \label{prop:zak}
	Assume that the chain is periodic and topologically trivial, in other words that it satisfies \eqref{eq:vperiodic} and \eqref{eq:trivial}. Then we have
	$$\varphi_j^{\mathrm{zak}}(v_1, v_2) = - \varphi_j^{\mathrm{zak}}(v_2, v_1) + O(\delta) \quad \text{and} \quad \varphi_j^{\mathrm{zak}}(\overline{v_1}, \overline{v_2}) = \varphi_j^{\mathrm{zak}}(v_1, v_2) + O(\delta).$$
	In particular, if $v_2 = \overline{v_1}$, we have $\varphi_j^{\mathrm{zak}}(v_1, \overline{v_1}) = O(\delta).$
\end{prop}
\begin{remark}
	\Cref{prop:zak} provides intuition on how to create structures supporting edge modes. \Cref{prop:zak} shows that distinct Zak phases can, in general, be achieved by swapping $v_1$ and $v_2$ while keeping the geometry fixed. Therefore, the defect specified in \eqref{eq:vdefect} results in Zak phases with opposite signs on the different sides of the edge. The reason we assume \eqref{eq:trivial} is to emphasize that distinct Zak phases can originate as a pure non-Hermitian effect, which disappears in the Hermitian limit as $\mathrm{Im}(v_1), \mathrm{Im}(v_2) \rightarrow 0$.
\end{remark}

We now turn to the analysis of the defect problem, specified in \eqref{eq:vdefect} and illustrated in \Cref{fig:edge}.
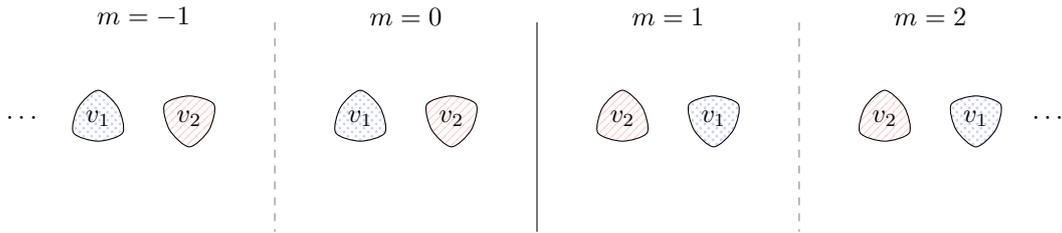
\begin{figure}[tbh]
	\centering
	\begin{tikzpicture}[scale=1.5]
		\pgfmathsetmacro{\rb}{0.25pt}
		\pgfmathsetmacro{\rs}{0.2pt}
		\coordinate (a) at (0.25,0);
		\coordinate (b) at (1.05,0);

		\draw (-0.5,0.85) -- (-0.5,-1);
		\def\Done{ plot [smooth cycle] coordinates {($(a)+(210:\rb)$) ($(a)+(270:\rs)$) ($(a)+(330:\rb)$) ($(a)+(30:\rs)$) ($(a)+(90:\rb)$) ($(a)+(150:\rs)$)
		}};
		\def\Dtwo{ plot [smooth cycle] coordinates {($(b)+(30:\rb)$) ($(b)+(90:\rs)$) ($(b)+(150:\rb)$) ($(b)+(210:\rs)$) ($(b)+(270:\rb)$) ($(b)+(330:\rs)$) }};
		\draw \Done;
		\draw \Dtwo;
		\pattern[pattern=north east lines, opacity=0.45, pattern color = red] \Done;
		\pattern[pattern=crosshatch dots, opacity=0.45, pattern color = blue] \Dtwo;
		\draw (a) node{$v_2$};
		\draw (b) node{$v_1$};
		\draw (0.65,0.9) node{$m=1$};

		\begin{scope}[xshift=-2.3cm]
			\coordinate (a) at (0.25,0);
			\coordinate (b) at (1.05,0);
			\draw[dashed, opacity=0.5] (-0.5,0.85) -- (-0.5,-1);
			\draw \Done;
			\draw \Dtwo;
			\pattern[pattern=north east lines, opacity=0.45, pattern color = red] \Dtwo;
			\pattern[pattern=crosshatch dots, opacity=0.45, pattern color = blue] \Done;
			\draw (a) node{$v_1$};
			\draw (b) node{$v_2$};
			\draw (0.65,0.9) node{$m=0$};
		\end{scope}
		\begin{scope}[xshift=-4.6cm]
			\coordinate (a) at (0.25,0);
			\coordinate (b) at (1.05,0);
			\draw \Done;
			\draw \Dtwo;
			\pattern[pattern=north east lines, opacity=0.45, pattern color = red] \Dtwo;
			\pattern[pattern=crosshatch dots, opacity=0.45, pattern color = blue] \Done;
			\draw (a) node{$v_1$};
			\draw (b) node{$v_2$};
			\draw (0.65,0.9) node{$m=-1$};
			\begin{scope}[xshift = 1.2cm]
				\draw (-1.6,0) node{$\cdots$};
			\end{scope};
		\end{scope}
		\begin{scope}[xshift=2.3cm]
			\coordinate (a) at (0.25,0);
			\coordinate (b) at (1.05,0);
			\draw[dashed, opacity=0.5] (-0.5,0.85) -- (-0.5,-1);
			\draw \Done;
			\draw \Dtwo;
			\pattern[pattern=north east lines, opacity=0.45, pattern color = red] \Done;
			\pattern[pattern=crosshatch dots, opacity=0.45, pattern color = blue] \Dtwo;
			\draw (a) node{$v_2$};
			\draw (b) node{$v_1$};
			\draw (0.65,0.9) node{$m=2$};
			\draw (1.7,0) node{$\cdots$};
		\end{scope}
	\end{tikzpicture}
	\caption[Single]{Illustration of the edge. The special case $v_1 = \overline{v_2}$ corresponds to local $\P \T$ symmetry. Legend:
		\raisebox{-2pt}{\tikz{
				\draw[pattern=crosshatch dots, opacity=0.45, pattern color = blue] (0,0) circle (5pt);
				\draw (0,0) circle (5pt);}} material parameter $v_1$,
		\raisebox{-2pt}{\tikz{
				\draw[pattern=north east lines, opacity=0.45, pattern color = red] (0,0) circle (5pt);
				\draw (0,0) circle (5pt);}} material parameter $v_2$. \label{fig:edge} }
\end{figure}
Since the geometry of this defect chain is periodic, we can utilise a Floquet-Bloch approach to derive a capacitance matrix characterisation of any localised mode. We define the matrices
$$A = \begin{pmatrix} 1 & b \\[0.3em] b & 1\end{pmatrix}, \quad B = |D_1|\begin{pmatrix} \ds \left(\delta v_2^2\right)^{-1} & \ds b\left(\delta v_1^2\right)^{-1}  \\[0.3em] \ds b\left(\delta v_1^2\right)^{-1}  & \ds \left(\delta v_2^2\right)^{-1} \end{pmatrix},$$
for some parameter $b\in \CC$, which can be interpreted as the decay of the localised mode between two consecutive resonators. The corresponding capacitance matrix, whose eigenvalue problem provides a discrete approximation to the localised mode and corresponding frequency, is given by
$$\C_\mathrm{edge} = B^{-1}C^\alpha A.$$
The eigenvalues of this matrix depend on $\alpha \in Y^*$. In order for localised modes to exist, there must be an eigenvalue $\mu = \mu_j^\alpha(b_0)$ of $\C_\mathrm{edge}$ which is constant in $\alpha$ for some $b=b_0$.
We can then compute $b$ as $b=b_\pm$, where
\begin{equation} \label{eq:b}
	b_\pm = \frac{1}{2}\left(l\left(1 - \frac{v_1^2}{v_2^2}\right) \pm \sqrt{l^2\left(1 - \frac{v_1^2}{v_2^2}\right)^2 + \frac{4v_1^2}{v_2^2}}\right), \qquad l = \frac{\lambda_2+\lambda_1}{\lambda_2-\lambda_1},
\end{equation}
with $\lambda_1 = C_{11}^{\pi/L} + C_{12}^{\pi/L}$ and $\lambda_2 = 2C_{11}^0$. Depending on the values of $v_1$ and $v_2$, we either have $|b_+| = |b_-| = 1$ (in which case there are no localised modes) or that $|b_\pm|<1$ while $|b_\mp| >1$. Based on these ideas, we can prove the following result \cite{ammari2020edge}.
\begin{thm} \label{thm:nonH}
	Suppose that the array of resonators has a defect in the material parameters specified by \eqref{eq:vdefect}. Then, for small $\delta$,
	\begin{itemize}
		\item if $v_1 = \overline{v_2}$ with $|\mathrm{Im}(v_1^2)| \leq \frac{\mathrm{Re}(v_1^2)}{\sqrt{l^2-1}}$ (unbroken $\P\T$-symmetry), the structure does not support simple localised modes in the subwavelength regime.
		\item if $v_1 = \overline{v_2}$ with $|\mathrm{Im}(v_1^2)| > \frac{\mathrm{Re}(v_1^2)}{\sqrt{l^2-1}}$ (broken $\P\T$-symmetry) or if $v_1 \neq \overline{v_2}$ (no $\P\T$-symmetry), the frequency $\omega$ of a simple localised mode in the subwavelength regime must satisfy
		$$\omega = \sqrt{\mu_j^\alpha(b_0)} + O(\delta).$$
		Here, $b_0$ is the value of $b$ specified by \eqref{eq:b} satisfying $|b_0| < 1$.
	\end{itemize}
\end{thm}
\begin{remark}
	\Cref{thm:nonH} provides a characterisation of the possible frequency $\omega$ and decay length $b$ of simple localised modes in the subwavelength regime. In order to prove that such mode indeed exists, we would have to prove that one eigenvalue $\mu_j^\alpha(b_0)$ is indeed constant in $\alpha$. Analytically, this is obscured by the fact that the capacitance coefficients have a complicated dependency on $\alpha$. Numerically, however, the eigenvalues $\mu$ are straightforward to compute, providing convincing evidence for such modes to exist. Moreover, the localised modes can easily be computed in large finite chains (see \Cref{fig:edgemodes}), demonstrating excellent agreement between the numerical and theoretically predicted values.
\end{remark}

	\begin{figure}[htb]
	\begin{subfigure}[b]{0.44\linewidth}
		\centering
		\includegraphics[width=0.95\linewidth]{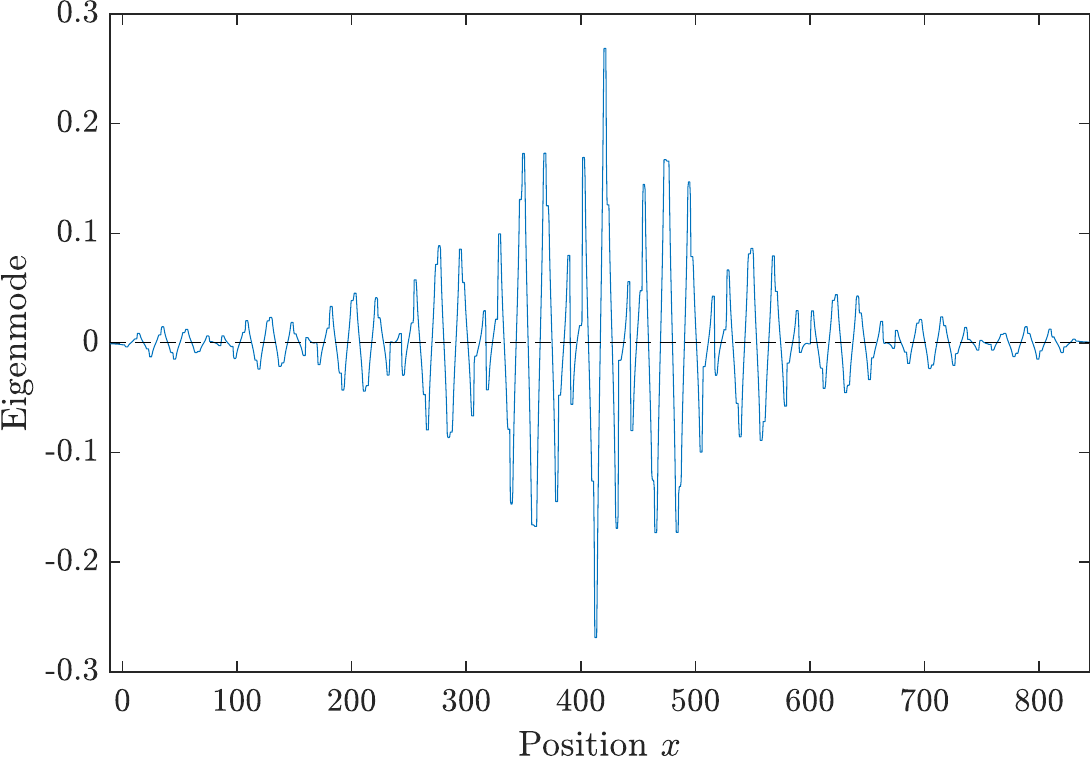}
		\caption{Localised mode in the case of small imaginary part.} \label{fig:model}
	\end{subfigure}\hfill
	\begin{subfigure}[b]{0.44\linewidth}
		\centering
		\includegraphics[width=0.95\linewidth]{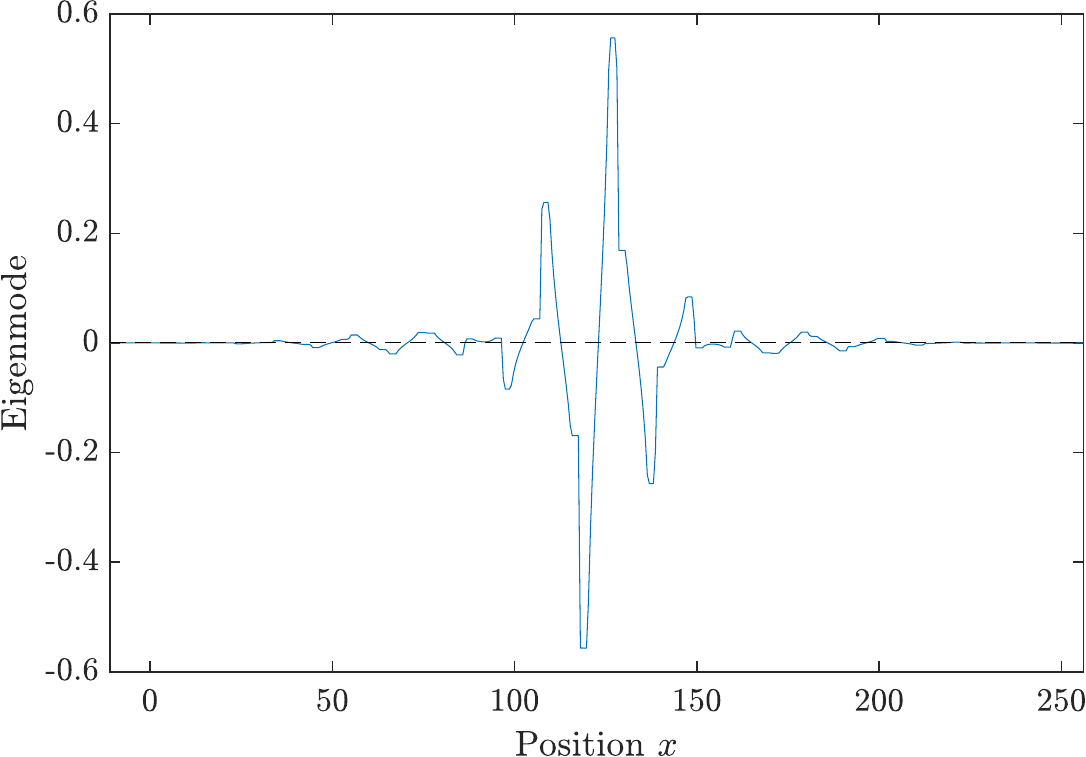}
		\caption{Localised mode in the case of large imaginary part.} \label{fig:mode}
	\end{subfigure}
	\caption{Plots of localised modes in finite, large, arrays of resonators. Observe the different $x$-axis scales in the two subfigures, demonstrating significantly different degrees of localisation. The localised modes have purely non-Hermitian origin and disappear in the Hermitian limit when $\Im(v_i)\to0$. }  \label{fig:edgemodes}
\end{figure}

\begin{figure}[h]
	\begin{minipage}{0.15\linewidth}
	\centering\small  	Hermitian:
	\end{minipage}\hfill	\begin{minipage}{0.8\linewidth}
	\begin{center}
	\begin{tikzpicture}[scale=1.9]
		\draw[->] (-0.8,0) -- (0.8,0) node[below]{\footnotesize$\alpha$};
		\draw [->] (0,-0.1) -- (0,1) node[left]{\footnotesize$\omega$};

		\draw plot[smooth] coordinates {(-0.8,0.7) (-0.55,0.6) (-0.2,0.75) (0.0,0.8) (0.2,0.75)  (0.55,0.6) (0.8,0.7)};
		\draw plot[smooth] coordinates {(0.0,0.0) (0.3,0.35)  (0.55,0.48) (0.8,0.4)};
		\draw plot[smooth] coordinates {(0.0,0.0) (-0.3,0.35)  (-0.55,0.48) (-0.8,0.4)};
		\draw (0,-0.23) node{\footnotesize $l_0 < 1/2$};
		\draw (0,-0.4) node{\footnotesize No band inversion};
		\begin{scope}[xshift=2.25cm]
			\coordinate (d1) at (0.55,0.55);
			\coordinate (d2) at (-0.55,0.55);

			\draw[->] (-0.8,0) -- (0.8,0) node[below]{\footnotesize$\alpha$};
			\draw [->] (0,-0.1) -- (0,1) node[left]{\footnotesize$\omega$};

			\draw plot[smooth] coordinates {(-0.8,0.4) (d2) (-0.2,0.75) (0.0,0.8) (0.2,0.75)  (d1) (0.8,0.4)};
			\draw plot[smooth] coordinates {(0.0,0.0) (0.3,0.35)  (d1) (0.8,0.7)};
			\draw plot[smooth] coordinates {(0.0,0.0) (-0.3,0.35)  (d2) (-0.8,0.7)};
			\draw (0,-0.4) node{\footnotesize Dirac cone degeneracy};
			\draw (0,-0.23) node{\footnotesize $l_0 = 1/2$};
		\end{scope}

		\begin{scope}[xshift=4.5cm]
			\draw[->] (-0.8,0) -- (0.8,0) node[below]{\footnotesize$\alpha$};
			\draw [->] (0,-0.1) -- (0,1) node[left]{\footnotesize$\omega$};

			\draw plot[smooth] coordinates {(-0.8,0.7) (-0.55,0.6) (-0.2,0.75) (0.0,0.8) (0.2,0.75)  (0.55,0.6) (0.8,0.7)};
			\draw plot[smooth] coordinates {(0.0,0.0) (0.3,0.35)  (0.55,0.48) (0.8,0.4)};
			\draw plot[smooth] coordinates {(0.0,0.0) (-0.3,0.35)  (-0.55,0.48) (-0.8,0.4)};

			\draw (0,-0.23) node{\footnotesize$l_0 > 1/2$};
			\draw (0,-0.4) node{\footnotesize Full band inversion};
		\end{scope}
	\end{tikzpicture}
\end{center}
\end{minipage}
\begin{minipage}{0.15\linewidth}
\centering	 \small Non-Hermitian:
\end{minipage}	\hfill
\begin{minipage}{0.8\linewidth}
	\begin{subfigure}{0.33\linewidth}
		\begin{center}
			\begin{tikzpicture}
				\draw (0,0) node{\includegraphics[clip,width=0.75\linewidth]{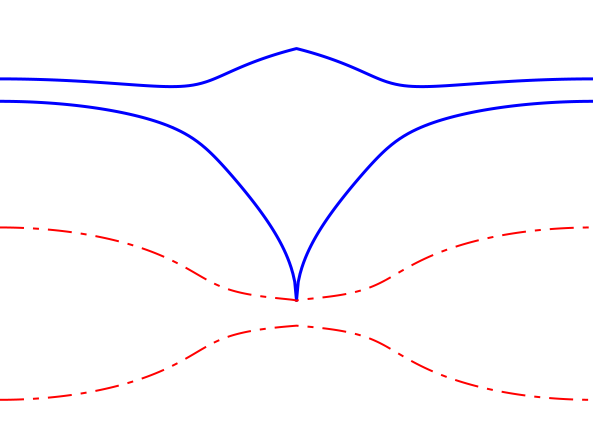}};
				\draw[->] (-1.65,-0.56) -- (1.65,-0.56) node[below]{\footnotesize $\alpha$};
				\draw[->] (0,-0.71) -- (0,1.2) node[right]{\footnotesize $\omega$};
				\draw (0,-1.4) node{\footnotesize Partial inversion};
				\draw (0,-1.1) node{\footnotesize $|v_1| < |v_2|$};
			\end{tikzpicture}
		\end{center}
	\end{subfigure}\begin{subfigure}{0.33\linewidth}
		\vspace{4pt}
		\begin{center}
			\begin{tikzpicture}
				\draw (0,0) node{\includegraphics[width=0.75\linewidth]{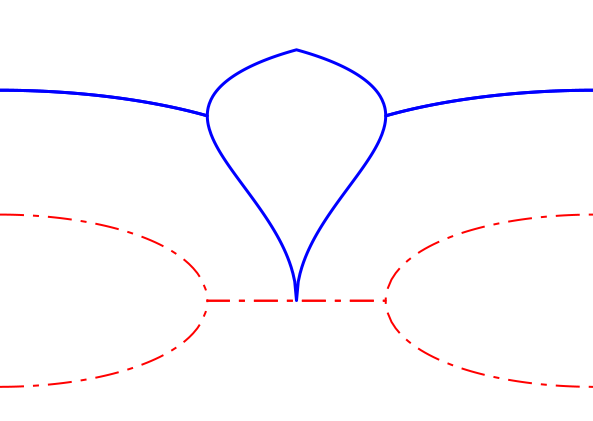}};
				\draw[->] (-1.65,-0.5) -- (1.65,-0.5) node[below]{\footnotesize $\alpha$};
				\draw[->] (0,-0.65) -- (0,1.2) node[right]{\footnotesize $\omega$};
				\draw (0,-1.4) node{\footnotesize Exceptional point};
				\draw (0,-1.1) node{\footnotesize $v_1 = \overline{v_2}$};
			\end{tikzpicture}
		\end{center}
	\end{subfigure}\begin{subfigure}{0.33\linewidth}
		\vspace{-4pt}
		\begin{center}
			\begin{tikzpicture}
				\draw (0,0) node{\includegraphics[width=0.75\linewidth]{bandnpt.png}};
				\draw[->] (-1.65,-0.56) -- (1.65,-0.56) node[below]{\footnotesize $\alpha$};
				\draw[->] (0,-0.71) -- (0,1.2) node[left]{\footnotesize $\omega$};
				\draw (0,-1.4) node{\footnotesize Partial inversion};
				\draw (0,-1.1) node{\footnotesize $|v_1| > |v_2|$};
				\draw (0.85,1.2) node{\includegraphics[width=0.33\linewidth]{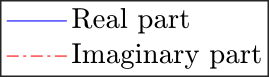}};
			\end{tikzpicture}
		\end{center}
	\end{subfigure}
\end{minipage}
\caption{Comparison between the topological phase transitions studied in \Cref{sec:bandinv} and \Cref{sec:pt_top}. In the Hermitian case studied in \Cref{sec:bandinv}, the Dirac cone found in the symmetric case $l_0 = 1/2$ can open into topologically distinct band gaps. In the non-Hermitian case, the degeneracy corresponds to an exceptional point. This exceptional point can open into separable bands which have distinct, albeit non-quantized, Zak phases.}\label{fig:phasetransition}
\end{figure}
In both the Hermitian and non-Hermitian cases, there are degeneracies associated to symmetric structures. Breaking the symmetry can open the degeneracy into topologically distinct band gaps as illustrated in \Cref{fig:phasetransition}. In the Hermitian case, the degeneracy is a linear intersection known as a \emph{Dirac cone}. In the non-Hermitian, $\P\T$-symmetric case, there can be exceptional point degeneracies which open when the $\P\T$-symmetry is detuned. In both the Hermitian and the non-Hermitian case, edge modes are created when the two topologically different phases are joined along an interface.

\subsection{Bound states in the continuum and Fano resonances}
Sometimes, localised modes can exist in periodic structures without a defect. This typically happens when the structure has certain symmetries, resulting in resonant modes in the radiation continuum whose far-field radiation vanishes. Such states are known as \emph{bound states in the continuum}.
\begin{defn}
	A resonant mode $u_n^\alpha$ of \eqref{eq:scattering_quasi} is said to be a bound state in the continuum if the corresponding resonant frequency $\omega^\alpha_n$ is real, satisfies $|\alpha| < \omega^\alpha_n/v$ and the mode satisfies,
	$$u_n^\alpha(x_l,x_0) = O(e^{-K|x_0|}), \quad |x_0| \to \infty, \quad K>0.$$
\end{defn}
As we shall see, there are two symmetry conditions required to achieve bound states in the continuum in the subwavelength resonator arrays. The first condition is a symmetry condition of the structure $D$, while the second condition is that $\alpha = 0$, which corresponds to modes that radiate perpendicularly to the structure. We then have the following result from \cite{ammari2021bound}.
\begin{thm}
	Assume that $d-d_l = 1$ and that $D = D_1\cup D_2 \subset Y$. We assume that $D$ satisfies the symmetry conditions $$\P D_1 = D_2, \quad \P_0 D_1 = D_1 \quad \P_0 D_2 = D_2,$$
	and $\delta_1v_2^2 = \delta_2 v_2^2 \in \R$, where $\P,\P_0:\R^d\to\R^d$ are the parity operators
	$$\P(x) = -x\quad\text{and}\quad \P_0(x_l,x_0) = (x_l,-x_0).$$
	Moreover, we assume that $\alpha = 0$. Then, for small enough $\delta$, the second resonant mode $u_2^0$ is a bound state in the continuum.
\end{thm}
Interestingly, we can never have the first resonant mode $u_1^0$ as a bound state in the continuum. This is due to the fact that the first mode corresponds to the broad, low-frequency response of the screen which is not due to the local resonance of the structure.

An interesting problem is now to describe the behaviour when the symmetry conditions are no longer satisfied. In particular, when $\P_0 D_i \neq D_i$, we can use \Cref{thm:higher_periodic} to conclude that $\omega_2^\alpha$ has a small but non-zero imaginary part. As we will see next, this causes an interesting \emph{Fano-type} transmission scattering. This phenomena emerges from the interference between the two coupled resonant frequencies of the pair $D$ of resonators, which have significantly different imaginary parts. The first resonant frequency corresponds to the ``continuum'' of states while the second resonant frequency originates from the resonant behaviour of the periodic structure and has a comparatively sharp response, corresponding to a ``discrete state''. The interaction between these two states leads to the creation of a Fano-type asymmetric transmission anomaly \cite{ammari2021bound,haifano}.

Given a unit vector $\w=(w_l,w_0)\in \R^d$ with $w_0>0$ we define the wave vectors
$$\k_+ = \frac{\omega}{v}\begin{pmatrix}w_l\\w_0\end{pmatrix}, \qquad \k_- = \frac{\omega}{v}\begin{pmatrix}w_l\\-w_0\end{pmatrix}.$$
We also let $\alpha = \frac{\omega}{v}w_l$. We now assume that
$$\uin(x) = c_1 e^{\i \k_-\cdot x} + c_2 e^{\i\k_+\cdot x},$$
and seek the behaviour of the solution $u$ of \eqref{eq:scattering_p}. In the first radiation continuum, the scattered field $u-\uin$ consists of a single propagating mode as $|x_0| \rightarrow \infty$.
We will write $ f\sim g$ to denote that two functions $f,g$ are equal up to exponentially decaying factors, in the sense that there is some constant $K>0$ such that
$$|f(x_l,x_0) - g(x_l,x_0)| = O(e^{-Kx_0}) \text{ as } x_0\rightarrow \infty.$$ We therefore have
\begin{equation}\label{eq:utot}
	u \sim \begin{cases}c_1 e^{\i\k_-\cdot x} + d_1 e^{\i\k_+\cdot x}, & x_0\rightarrow \infty, \\ c_2 e^{\i\k_+\cdot x} + d_2 e^{\i\k_-\cdot x}, & x_0\rightarrow -\infty,\end{cases}
\end{equation}
where
\begin{equation} \label{eq:S_def}
	\begin{pmatrix} d_1 \\ d_2 \end{pmatrix} =  S\begin{pmatrix} c_1 \\ c_2 \end{pmatrix}, \qquad S = \begin{pmatrix} r_+ & t_- \\ t_+ & r_- \end{pmatrix}.
\end{equation}
$S$ is known as the scattering matrix.
The reflection and transmission coefficients $r_+, t_+$ are the coefficients of the outgoing part of the field in the case $\uin(x) =  e^{\i\k_-\cdot x}$, \ie{} when the incident field is a plane wave from the positive $x_0$ direction (and reversely for $r_-, t_-$). The following theorem was proved in \cite{ammari2021bound}.
\begin{thm} \label{thm:phano}
	Assume that  $\P D = D$ and that $0\leq \omega \leq K\sqrt{\delta}$ for some constant $K$. Let $c=\int_{\D}y_0\psi_1^0(t)\de \sigma(y)$ and assume that $c \neq 0$. Then we have the following asymptotic expansion of the scattering matrix as $\delta \to 0$
	\begin{equation}
		S = \frac{\omega_1^\alpha}{\omega_1^\alpha-\omega}\begin{pmatrix} 1&1\\1&1\end{pmatrix} +\frac{2\i\omega\Im(\omega_2^\alpha)}{(\omega_2^\alpha)^2-\omega^2}\begin{pmatrix} 1&-1\\-1&1\end{pmatrix} - \begin{pmatrix} 1&0\\0&1\end{pmatrix} + O(\delta^{1/2}),
	\end{equation}
	where the error term is uniform with respect to $\omega$.
\end{thm}
The scattering matrix $S$ contains two transmission peaks originating from the resonances $\omega_1^\alpha$ and $\omega_2^\alpha$. Since the imaginary part of $\omega_2^\alpha$ is very small, this corresponds to a sharp peak which will interfere with the broader peak associated to $\omega_1^\alpha$ to create an asymmetric, Fano-type, transmission peak as illustrated in \Cref{fig:fano}.

\begin{figure}[tbh]
	\begin{subfigure}[b]{0.48\linewidth}
	\centering
	\begin{tikzpicture}[scale=0.4]
		\begin{scope}[rotate=40]
			\draw (0.8,0) circle (0.5);
			\draw (-0.8,0) circle (0.5);
		\end{scope}
		\coordinate (oo) at (220:1.8);
		\coordinate (a) at ($(0:3)+(oo)$);
		\coordinate (b) at ($(40:3)+(oo)$);
		\coordinate (o) at ($(0,0)+(oo)$);
		\draw[opacity=0.7] (a) -- (o) -- (b);
		\pic [draw, opacity=0.7, ->, "$\theta$", angle eccentricity=1.3,angle radius=0.7cm] {angle = a--o--b};
		\node at (7,0){$\cdots$};
		\node at (-7,0){$\cdots$};
		\draw[->,thick] (-2,-4) -- (-1,-2) node[pos=0.7,below right]{$\uin$};
		\draw[dashed] (-2,3) -- (-2,-2.5) (2,3) -- (2,-2.5);
		\draw[<->] (-2,2.5) -- (2,2.5) node[pos=0.5,above]{$L$};
		\begin{scope}[xshift=4cm,rotate=40]
			\draw (0.8,0) circle (0.5);
			\draw (-0.8,0) circle (0.5);
		\end{scope}
		\begin{scope}[xshift=-4cm,rotate=40]
			\draw (0.8,0) circle (0.5);
			\draw (-0.8,0) circle (0.5);
		\end{scope}
	\end{tikzpicture}
	\vspace{0.5cm}
	\caption{Sketch of a screen of resonators with an incident plane wave $\uin$. In this case, we have resonators arranged in a $\mathcal{P}$-symmetric dimer that is inclined at an angle of $\theta$ to the plane of the metascreen.} \label{fig:2D}
	\end{subfigure}\hfill
	\begin{subfigure}[b]{0.48\linewidth}
		\begin{center}
			\includegraphics[width=1\linewidth]{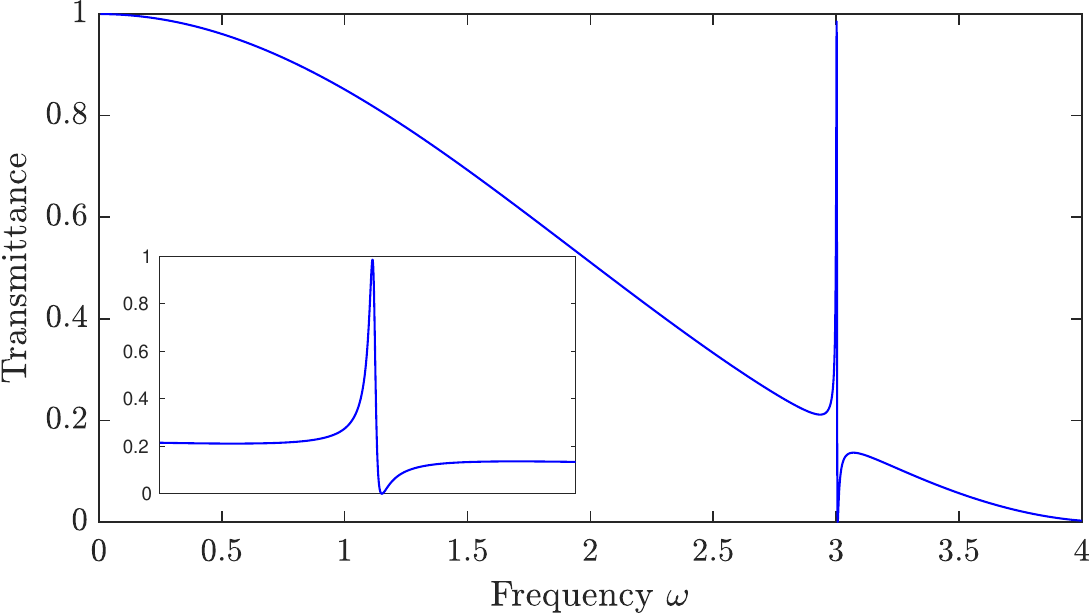}
		\end{center}
		\caption{Transmission spectrum for the system in (a) with $\theta = 0.025\pi$. } \label{fig:fanoT}
	\end{subfigure}
\caption{Sketch of a resonator screen (a) and corresponding transmission spectrum (b). If $\theta=0$, the resonant modes at $\alpha=0$ are bound states in the continuum and the transmission spectrum does not show a resonant peak at $\omega \approx 3$. For small but non-zero $\theta$, a characteristic Fano-type transmission peak appears at the resonant frequency.} \label{fig:fano}
\end{figure}
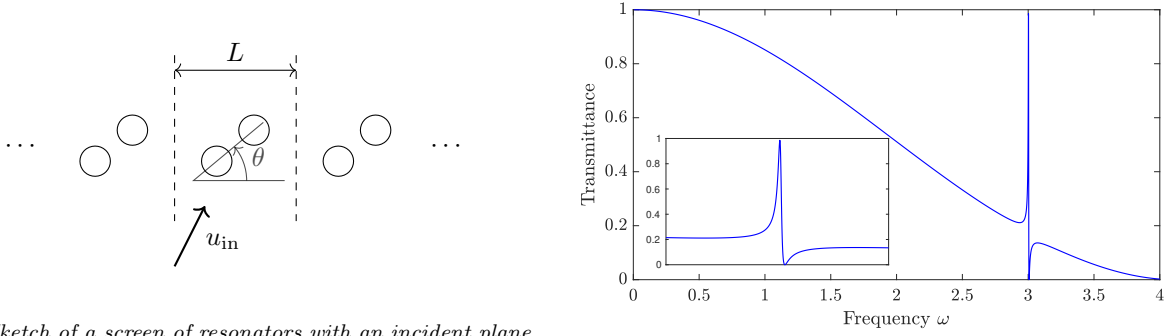

\subsection{Extraordinary transmission and unidirectional reflection}
We now investigate further scattering phenomena of resonator screens. In particular, we will extend the analysis in the previous section (valid in the Hermitian case when $v_i\in \R$) to the non-Hermitian case with a $\P\T$-symmetric screen. We assume that
$$\delta_1v_1^2 = a+\i b, \quad \delta_2v_2^2 = a-\i b,$$
in other words that $\delta_1v_1^2 = \T(\delta_2v_2^2 )$. Then we can apply \Cref{thm:higher_periodic} to conclude that
$$\omega_2^\alpha = \sqrt{\frac{2aC_{11}^0}{|D_1|}} + \frac{\i k_0}{4C_{11}^0}\left(\frac{b^2|Y_l|}{a^2} - \frac{c^2}{|Y_l|}\right) + O(\delta^{3/2}),$$
where $k_0$ and $c$ are defined as in the previous section. The following theorem, proved in \cite{ammari2020exceptional}, describes the scattering behaviour of this $\P\T$-symmetrical screen of resonators.
\begin{thm} \label{thm:unidir}
	Assume that $D$ is $\P\T$-symmetric, so that $\P D_1 = D_2$ while $\delta_1v_1^2 = a+\i b$ and $\delta_2v_2^2 = a-\i b$ for $a>0,b\geq0$. Moreover, assume that $\P_3D_1=D_2$. Let $c=\int_{\D}y_0\psi_1^0(t)\de \sigma(y)$ and $\omega_* = \sqrt{\frac{2aC_{11}^0}{|D_1|}}$. Assume that $b|Y_l| \neq a|c| $ and that $\omega\in \R$ is in the subwavelength regime such that $\omega-\omega_* = O(\delta)$.	We then have the following asymptotic expansion of the scattering matrix:
	$$S =\frac{2\i \omega\Im(\omega_2^\alpha )}{(\omega_2^\alpha)^2 - \omega^2}\begin{pmatrix} 1&-1\\-1&1\end{pmatrix} + \frac{2k_0 b}{a|D_1|\big((\omega_2^\alpha)^2 - \omega^2\big)}\begin{pmatrix} - ac & \i b|Y_l| \\ \i b|Y_l| & ac \end{pmatrix} - \begin{pmatrix} 1&0\\0&1\end{pmatrix} + O(\delta^{1/2}),$$
	where the error term is uniform with respect to $\omega$ in a neighbourhood of $\omega_*$.	In particular, we have
	$$r_\pm =- \frac{\omega_*^2-\omega^2 \pm \frac{2 k_0 bc}{|D_1|}}{(\omega_2^\alpha)^2 - \omega^2}+ O(\delta^{1/2}),$$
	and, at leading order, $r_+$ and $r_-$ vanish at $\omega = \omega_+$ and $\omega = \omega_-$, respectively, which are given by
	$$\omega_+^2 = \omega_*^2 +  \frac{2 k_0 bc}{|D_1|}, \quad \omega_-^2 = \omega_*^2-\frac{2 k_0 bc}{|D_1|}.$$
\end{thm}
Comparing \Cref{thm:phano} (which is valid only for $b=0$) and \Cref{thm:unidir}, we see that there is an extra term in \Cref{thm:unidir} corresponding to the non-zero gain and loss $b\neq 0$. This term is responsible for the approximate zeros of $r_+$ and $r_-$ leading to unidirectional reflection. This is numerically verified in \Cref{fig:unidir}, where we observe that the two reflectances vanish on different sides of the critical frequency $\omega_*$. Moreover, at least formally, we see that when $\Im(\omega_2^\alpha )=0$, the singularity of the scattering matrix $S$ will not vanish in the case $b\neq0$. In the Hermitian case, real resonances correspond to bound states in the continuum, which decouple from the far-field and therefore cannot be excited by incoming waves. In the non-Hermitian case, however, we can have real resonances with modes which are excited by incoming waves. Such resonances correspond to extraordinary transmission, where the transmitted field is greatly amplified. This amplification, which is impossible in the Hermitian case due to energy conservation, is possible due to the energy input in the non-Hermitian case \cite{ammari2020exceptional}.

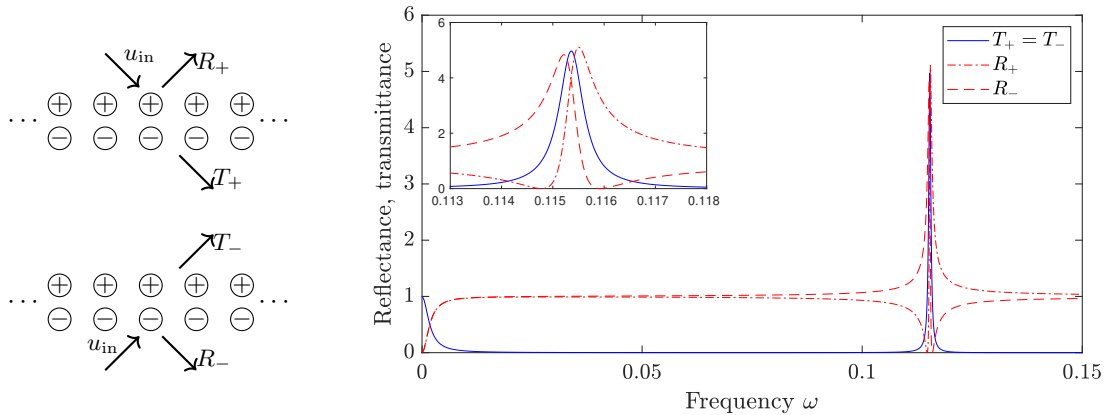
\begin{figure}[h]
	\begin{subfigure}{0.35\linewidth}
			\begin{center}
			\begin{tikzpicture}[scale=0.15]
				\draw (0,-1.5) circle (1);
				\draw (0,1.5) circle (1);
				\node at (0,-1.5){\small $-$};
				\node at (0,1.5){\small $+$};
				\node at (11,0){$\cdots$};
				\node at (-11,0){$\cdots$};
				\draw[->,thick] (-4,6) -- (-1,3) node[pos=0.7,xshift=4pt,yshift=8pt]{\small $\uin$};
				\draw[->,thick] (1,3) -- (4,6) node[pos=0.7,right]{\small $R_+$};
				\draw[->,thick] (2.5,-3) -- (5.5,-6) node[pos=0.7,right]{\small $T_+$};
				\begin{scope}[xshift=4cm]
					\draw (0,-1.5) circle (1);
					\draw (0,1.5) circle (1);
					\node at (0,-1.5){\small $-$};
					\node at (0,1.5){\small $+$};
				\end{scope}
				\begin{scope}[xshift=-4cm]
					\draw (0,-1.5) circle (1);
					\draw (0,1.5) circle (1);
					\node at (0,-1.5){\small $-$};
					\node at (0,1.5){\small $+$};
				\end{scope}
				\begin{scope}[xshift=8cm]
					\draw (0,-1.5) circle (1);
					\draw (0,1.5) circle (1);
					\node at (0,-1.5){\small $-$};
					\node at (0,1.5){\small $+$};
				\end{scope}
				\begin{scope}[xshift=-8cm]
					\draw (0,-1.5) circle (1);
					\draw (0,1.5) circle (1);
					\node at (0,-1.5){\small $-$};
					\node at (0,1.5){\small $+$};
				\end{scope}
				\begin{scope}[yshift=-16cm]
					\draw (0,-1.5) circle (1);
					\draw (0,1.5) circle (1);
					\node at (0,-1.5){\small $-$};
					\node at (0,1.5){\small $+$};
					\node at (11,0){$\cdots$};
					\node at (-11,0){$\cdots$};
					\draw[->,thick] (-4,-6) -- (-1,-3) node[pos=0.7,left]{\small $\uin$};
					\draw[->,thick] (2.5,3) -- (5.5,6) node[pos=0.7,right]{\small $T_-$};
					\draw[->,thick] (1,-3) -- (4,-6) node[pos=0.7,right]{\small $R_-$};
					\begin{scope}[xshift=4cm]
						\draw (0,-1.5) circle (1);
						\draw (0,1.5) circle (1);
						\node at (0,-1.5){\small $-$};
						\node at (0,1.5){\small $+$};
					\end{scope}
					\begin{scope}[xshift=-4cm]
						\draw (0,-1.5) circle (1);
						\draw (0,1.5) circle (1);
						\node at (0,-1.5){\small $-$};
						\node at (0,1.5){\small $+$};
					\end{scope}
					\begin{scope}[xshift=8cm]
						\draw (0,-1.5) circle (1);
						\draw (0,1.5) circle (1);
						\node at (0,-1.5){\small $-$};
						\node at (0,1.5){\small $+$};
					\end{scope}
					\begin{scope}[xshift=-8cm]
						\draw (0,-1.5) circle (1);
						\draw (0,1.5) circle (1);
						\node at (0,-1.5){\small $-$};
						\node at (0,1.5){\small $+$};
					\end{scope}
				\end{scope}
			\end{tikzpicture}
		\end{center}
	\end{subfigure}
	\begin{subfigure}{0.6\linewidth}
	\begin{center}
		\includegraphics[width=\linewidth]{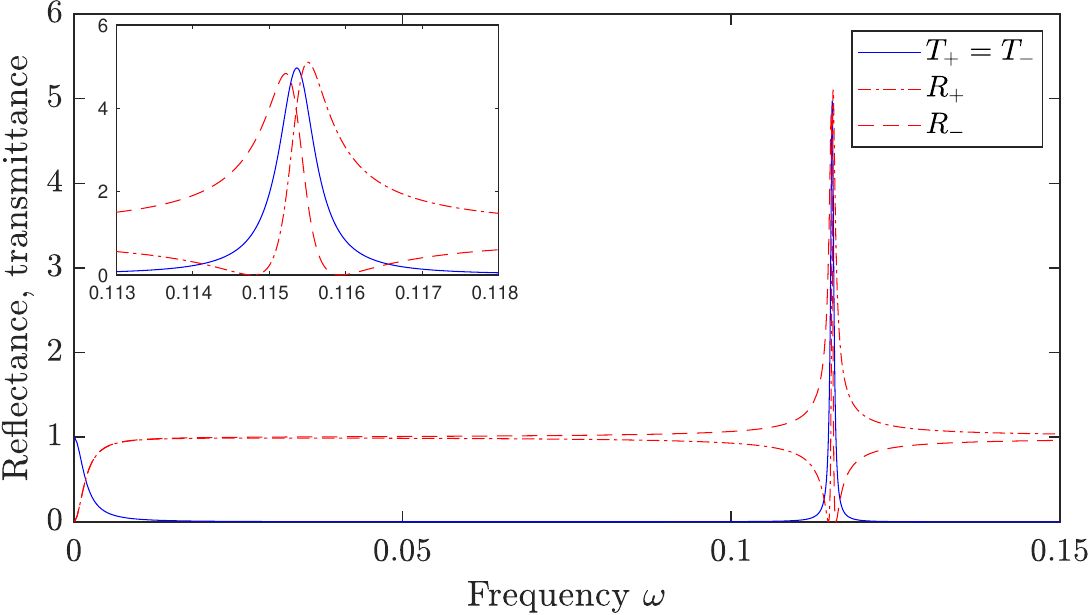}
	\end{center}
	\end{subfigure}
	\caption{Plot of the transmittance $T_\pm = |t_\pm|^2$ (blue) and reflectance $R_\pm = |r_\pm|^2$ (red) as functions of the frequency. Due to reciprocity we have $T_+ = T_-$. However, we observe that $R_+\neq R_-$. In particular, there are points when $R_\pm$ vanish while $R_\mp$ is non-zero, known as unidirectional reflection. Moreover, in this non-Hermitian case the transmission coefficients are not bounded by unity and can attain large peak values, known as extraordinary transmission.} \label{fig:unidir}
\end{figure}

\subsection{Time-modulated metamaterials}\label{sec:time}
We have now seen a variety of phenomena that can occur when the material parameters vary periodically in the spatial variable $x$. Mathematically, we can treat the time variable $t$ in a similar fashion. If the material parameters depend periodically on $t$, there can be conceptually similar phenomena which nevertheless have fundamentally different physical implications.

The Helmholtz equation we have studied so far originates from the scalar wave equation when posed in the frequency domain, and is valid only when the material parameters are constant in time. To study the time-dependent case, we return to the wave equation
\begin{equation}\label{eq:wave}
	\left(\frac{\p }{\p t } \frac{1}{\kappa(x,t)} \frac{\p}{\p t} - \nabla \cdot \frac{1}{\rho(x,t)} \nabla\right) u(x,t) = 0, \quad x\in \R^d, t\in \R.
\end{equation}
Here, $\kappa(x,t)$ and $\rho(x,t)$ are the material parameters. We consider the case of a  finite collection of resonators $D=D_1\cup \dots \cup D_N$ in $d=3$ spatial dimensions. We assume that the modulation is only performed inside the resonators, so that
 \begin{equation} \label{eq:resonatormod}
	\kappa(x,t) = \begin{cases}
		\kappa_0, & x \in \R^3 \setminus \overline{D}, \\  \kappa_i(t), & x\in D_i,
	\end{cases}, \qquad \rho(x,t) = \begin{cases}
		\rho_0, & x \in \R^3 \setminus \overline{D}, \\  \rho_i(t), & x\in D_i. \end{cases}
\end{equation}
We assume that $\kappa_i \in C^1(\R)$ for each $i=1,\dots,N$. Moreover, we assume that they are periodic in $t$ with period $T$ and frequency $\Omega = \frac{2\pi}{T}$ and consider only the case when the modulation of $1/\rho$ and $1/\kappa$ consist of a finite Fourier series with a large number of Fourier coefficients:
$$\frac{1}{\rho_i(t)} = \sum_{n = -M}^M r_{i,n} e^{\i n \Omega t}, \qquad \frac{1}{\kappa_i(t)} = \sum_{n = -M}^M k_{i,n} e^{\i n \Omega t},$$
for some $M\in \N$ satisfying $M = O\left(\delta^{-\gamma/2}\right)$ for some $0<\gamma<1$. 

As before we define the (time-dependent) contrast parameter and wave speeds as
$$\delta_i(t) = \frac{\rho_i(t)}{\rho_0}, \quad v_i(t) = \sqrt{\frac{\kappa_i(t)}{\rho_i(t)}}, \quad v_0 = \sqrt{\frac{\kappa_0}{\rho_0}}, $$
and assume that
$$\delta_i(t) = O(\delta), \quad v_i(t) = O(1) \quad v = O(1), \quad \text{for all } \ t\in \R,$$
for $i=1,\dots,N$, where $\delta \ll 1$.

The notion of frequency is slightly altered in this time-modulated setting. Since the wave equation \eqref{eq:wave} is periodic in $t$, we can apply the Floquet transform in $t$ and obtain the differential problem
\begin{equation} \label{eq:wave_transf}
	\begin{cases}\ \ds \left(\frac{\p }{\p t } \frac{1}{\kappa(x,t)} \frac{\p}{\p t} - \nabla \cdot \frac{1}{\rho(x,t)} \nabla\right) u(x,t) = 0,\\[0.3em]
		\	u(x,t)e^{-\i \omega t} \text{ is $T$-periodic in $t$}.
	\end{cases}
\end{equation}
It is apparent that $\omega$ plays the role of the quasiperiodicity in the Floquet-Bloch theory, and we will refer to $\omega$ as a quasifrequency. Since $\omega$ is defined modulo $\Omega$, we define the time-Brillouin zone $Y_t^* := \mathbb{C} / (\Omega \Z)$. Observe that we allow complex quasifrequencies; real $\omega$ correspond to solutions which are bounded in $t$, whereas non-real $\omega$ correspond to exponentially increasing or decaying solutions.

Due to the periodic nature of $Y_t^*$, the usual definition of subwavelength frequencies does not apply to quasifrequencies. For example, in the particular case when $\Omega = O(\delta^{1/2})$ (which will be of interest later on), the whole Brillouin zone scales as $O(\delta^{1/2})$, meaning that \emph{all} quasifrequencies tend to zero as $\delta \to 0$. In order to distinguish between these quasifrequencies, we introduce the following definition.
 	\begin{defn}[Subwavelength quasifrequency] \label{def:sub}
	A quasifrequency $\omega = \omega(\delta) \in Y^*_t$ of \eqref{eq:wave_transf} is said to be a subwavelength quasifrequency if there is a corresponding Bloch solution $u(x,t)$, depending continuously on $\delta$, which can be written as
	$$u(x,t)= e^{\i \omega t}\sum_{n = -\infty}^\infty v_n(x)e^{\i n\Omega t},$$
	where
	$$\omega \rightarrow 0 \ \text{and} \ M\Omega \rightarrow 0 \ \text{as} \ \delta \to 0,$$
	for some integer-valued function $M=M(\delta)$ such that, as $\delta \to 0$, we have
	$$\sum_{n = -\infty}^\infty \|v_n\|_{L^2(D)} = \sum_{n = -M}^M \|v_n\|_{L^2(D)} + o(1).$$
\end{defn}
The following theorem, proved in \cite{ammari2020time}, gives the capacitance matrix approximation to the subwavelength quasifrequencies as $\delta \to 0$.
\begin{thm}
	Assume that the material parameters are given by \eqref{eq:resonatormod}. Then, as $\delta \to 0$, the quasifrequencies $\omega \in Y^*_t$ to the wave equation \eqref{eq:wave} in the subwavelength regime are, to leading order, given by the quasifrequencies of the system of ordinary differential equations for $y_i(t)$,
	\begin{equation}\label{eq:C_ODE}
		\sum_{j=1}^N C_{ij} y_j(t) = -|D_i|\frac{\dx}{\dx t}\left(\frac{1}{\delta_iv_i^2}\frac{\dx y_i}{\dx t}\right),
	\end{equation}
	for $i=1,\dots,N$.
\end{thm}
The matrix appearing in the left-hand side of \eqref{eq:C_ODE} is the capacitance matrix. We can rewrite \eqref{eq:C_ODE} into the following system of Hill equations:
\begin{equation}	\label{eq:hill}
	\Psi''(t)+ M(t)\Psi(t)=0,
\end{equation}
where the vector $\Psi$ and the matrix $M$ are defined as
$$\Psi(t) = \left(\frac{y_i(t)}{\sqrt{\delta_i(t)v_i^2(t)}}\right)_{i=1}^N, \quad M(t) = W_1(t)C W_2(t) + W_3(t),$$
with $W_1, W_2$ and $W_3$ being the diagonal matrices with diagonal entries
$$\left(W_1\right)_{ii} = \frac{\sqrt{\delta_iv_i^2}}{|D_i|}, \qquad \left(W_2\right)_{ii} =\sqrt{\delta_iv_i^2}, \qquad \left(W_3\right)_{ii} = \sqrt{\delta_iv_i^2}\frac{\dx }{\dx t}\frac{1}{\delta_iv_i^2}\frac{\dx \sqrt{\delta_iv_i^2}}{\dx t},$$
for $i=1,\dots,N$.
\begin{figure}
	\centering
	\begin{subfigure}[b]{0.4\linewidth}
		\centering
		\begin{tikzpicture}
			\draw (-1.5,0) circle (1);
			\draw (1.5,0) circle (1);
			\node at (-1.5,0){$D_1$};
			\node at (1.5,0){$D_2$};
			\node at (-1.6,-1.5){\small $\kappa_1(t) =1+\epsilon\sin(\Omega t)$};
			\node at (1.6,-1.5){\small $\kappa_2(t) = 1-\epsilon\sin(\Omega t)$};
		\end{tikzpicture}
		\vspace{0.5cm}
		\caption{A $\mathcal{PT}$-symmetric pair of spherical resonators.}
	\end{subfigure}
	\hspace{0.1cm}
	\begin{subfigure}[b]{0.55\linewidth}
		\includegraphics[width=\linewidth]{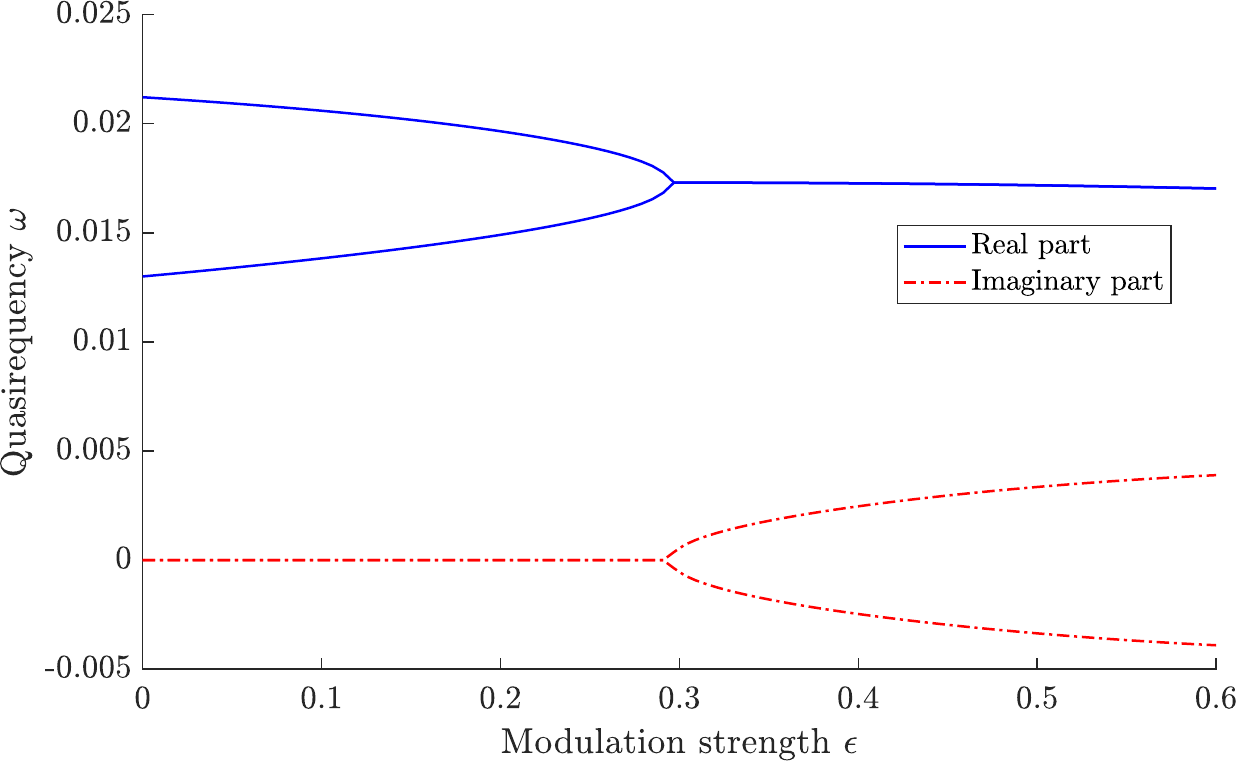}
		\caption{The two subwavelength resonant frequencies.}
	\end{subfigure}
	\caption{Similarly as in the non-Hermitian case demonstrated in \Cref{fig:EP2}, there can be exceptional points in the time-modulated case. The two subwavelength quasifrequencies of a pair of $\mathcal{PT}$-symmetric, time-modulated resonators can be approximated through the capacitance formulation, and an asymptotic exceptional point occurs at $\epsilon \approx0.3$.} \label{fig:EPtime}
\end{figure}
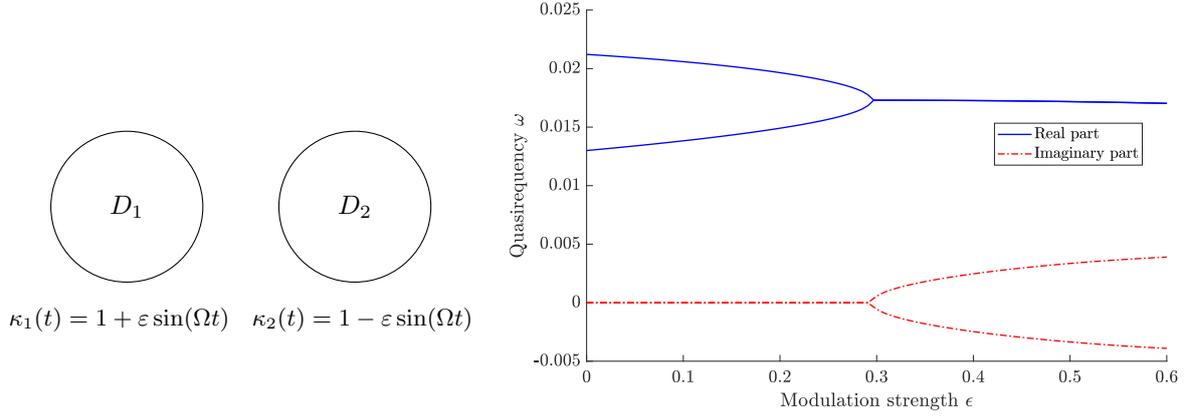
The time-modulated case shares many similarities with the non-Hermitian case (with complex parameters) studied before. Both these cases have energy input and output to the system, and not surprisingly we can find exceptional points in the time-modulated systems. \Cref{fig:EPtime} shows the emergence of an exceptional point in a dimer structure, similar to \Cref{fig:EP2} but instead due to the time-modulation.
\subsection{Near-zero metamaterials}
In \Cref{sec:doubleneg} we saw the emergence of exotic parameter values, namely negative effective parameters, due to the small-scale structure of the metamaterial. In the present section we will observe effective material parameters which are close to zero, in which case the Helmholtz equation reduces to the Laplace equation and wave propagation occurs without phase change (corresponding to ``infinite'' phase velocity).

\subsubsection{Near-zero refractive index in honeycomb crystals}
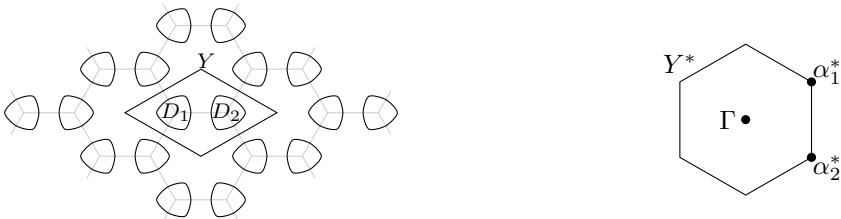
\begin{figure}[htb]
\begin{subfigure}[b]{0.5\linewidth}
	\centering
\begin{tikzpicture}[scale=1]
	\begin{scope}[xshift=-5cm,scale=1]
		\coordinate (a) at (1,{1/sqrt(3)});
		\coordinate (b) at (1,{-1/sqrt(3)});
		\pgfmathsetmacro{\rb}{0.25pt}
		\pgfmathsetmacro{\rs}{0.2pt}

		\draw (0,0) -- (1,{1/sqrt(3)}) node[xshift=2pt,yshift=3pt]{\footnotesize$Y$} -- (2,0) -- (1,{-1/sqrt(3)}) -- cycle;
		\begin{scope}[xshift = 1.33333cm]
			\draw plot [smooth cycle] coordinates {(0:\rb) (60:\rs) (120:\rb) (180:\rs) (240:\rb) (300:\rs) };
			\draw (0,0) node{\footnotesize $D_2$};
		\end{scope}
		\begin{scope}[xshift = 0.666667cm, rotate=60]
			\draw plot [smooth cycle] coordinates {(0:\rb) (60:\rs) (120:\rb) (180:\rs) (240:\rb) (300:\rs) };
			\draw (0,0) node{\footnotesize $D_1$};
		\end{scope}

		\draw[opacity=0.2] ({2/3},0) -- ({4/3},0)
		($0.5*(1,{1/sqrt(3)})$) -- ({2/3},0)
		($0.5*(1,{-1/sqrt(3)})$) -- ({2/3},0)
		($(1,{1/sqrt(3)})+0.5*(1,{-1/sqrt(3)})$) -- ({4/3},0)
		($0.5*(1,{1/sqrt(3)})+(1,{-1/sqrt(3)})$) -- ({4/3},0);

		\begin{scope}[shift = (a)]
			\begin{scope}[xshift = 1.33333cm]
				\draw plot [smooth cycle] coordinates {(0:\rb) (60:\rs) (120:\rb) (180:\rs) (240:\rb) (300:\rs) };
			\end{scope}
			\begin{scope}[xshift = 0.666667cm, rotate=60]
				\draw plot [smooth cycle] coordinates {(0:\rb) (60:\rs) (120:\rb) (180:\rs) (240:\rb) (300:\rs) };
			\end{scope}
			\draw[opacity=0.2] ({2/3},0) -- ({4/3},0)
			($0.5*(1,{1/sqrt(3)})$) -- ({2/3},0)
			($0.5*(1,{-1/sqrt(3)})$) -- ({2/3},0)
			($(1,{1/sqrt(3)})+0.5*(1,{-1/sqrt(3)})$) -- ({4/3},0)
			($0.5*(1,{1/sqrt(3)})+(1,{-1/sqrt(3)})$) -- ({4/3},0);
		\end{scope}
		\begin{scope}[shift = (b)]
			\begin{scope}[xshift = 1.33333cm]
				\draw plot [smooth cycle] coordinates {(0:\rb) (60:\rs) (120:\rb) (180:\rs) (240:\rb) (300:\rs) };
			\end{scope}
			\begin{scope}[xshift = 0.666667cm, rotate=60]
				\draw plot [smooth cycle] coordinates {(0:\rb) (60:\rs) (120:\rb) (180:\rs) (240:\rb) (300:\rs) };
			\end{scope}
			\draw[opacity=0.2] ({2/3},0) -- ({4/3},0)
			($0.5*(1,{1/sqrt(3)})$) -- ({2/3},0)
			($0.5*(1,{-1/sqrt(3)})$) -- ({2/3},0)
			($(1,{1/sqrt(3)})+0.5*(1,{-1/sqrt(3)})$) -- ({4/3},0)
			($0.5*(1,{1/sqrt(3)})+(1,{-1/sqrt(3)})$) -- ({4/3},0);
		\end{scope}
		\begin{scope}[shift = ($-1*(a)$)]
			\begin{scope}[xshift = 1.33333cm]
				\draw plot [smooth cycle] coordinates {(0:\rb) (60:\rs) (120:\rb) (180:\rs) (240:\rb) (300:\rs) };
			\end{scope}
			\begin{scope}[xshift = 0.666667cm, rotate=60]
				\draw plot [smooth cycle] coordinates {(0:\rb) (60:\rs) (120:\rb) (180:\rs) (240:\rb) (300:\rs) };
			\end{scope}
			\draw[opacity=0.2] ({2/3},0) -- ({4/3},0)
			($0.5*(1,{1/sqrt(3)})$) -- ({2/3},0)
			($0.5*(1,{-1/sqrt(3)})$) -- ({2/3},0)
			($(1,{1/sqrt(3)})+0.5*(1,{-1/sqrt(3)})$) -- ({4/3},0)
			($0.5*(1,{1/sqrt(3)})+(1,{-1/sqrt(3)})$) -- ({4/3},0);
		\end{scope}
		\begin{scope}[shift = ($-1*(b)$)]
			\begin{scope}[xshift = 1.33333cm]
				\draw plot [smooth cycle] coordinates {(0:\rb) (60:\rs) (120:\rb) (180:\rs) (240:\rb) (300:\rs) };
			\end{scope}
			\begin{scope}[xshift = 0.666667cm, rotate=60]
				\draw plot [smooth cycle] coordinates {(0:\rb) (60:\rs) (120:\rb) (180:\rs) (240:\rb) (300:\rs) };
			\end{scope}
			\draw[opacity=0.2] ({2/3},0) -- ({4/3},0)
			($0.5*(1,{1/sqrt(3)})$) -- ({2/3},0)
			($0.5*(1,{-1/sqrt(3)})$) -- ({2/3},0)
			($(1,{1/sqrt(3)})+0.5*(1,{-1/sqrt(3)})$) -- ({4/3},0)
			($0.5*(1,{1/sqrt(3)})+(1,{-1/sqrt(3)})$) -- ({4/3},0);
		\end{scope}
		\begin{scope}[shift = ($(a)+(b)$)]
			\begin{scope}[xshift = 1.33333cm]
				\draw plot [smooth cycle] coordinates {(0:\rb) (60:\rs) (120:\rb) (180:\rs) (240:\rb) (300:\rs) };
			\end{scope}
			\begin{scope}[xshift = 0.666667cm, rotate=60]
				\draw plot [smooth cycle] coordinates {(0:\rb) (60:\rs) (120:\rb) (180:\rs) (240:\rb) (300:\rs) };
			\end{scope}
			\draw[opacity=0.2] ({2/3},0) -- ({4/3},0)
			($0.5*(1,{1/sqrt(3)})$) -- ({2/3},0)
			($0.5*(1,{-1/sqrt(3)})$) -- ({2/3},0)
			($(1,{1/sqrt(3)})+0.5*(1,{-1/sqrt(3)})$) -- ({4/3},0)
			($0.5*(1,{1/sqrt(3)})+(1,{-1/sqrt(3)})$) -- ({4/3},0);
		\end{scope}
		\begin{scope}[shift = ($-1*(a)-(b)$)]
			\begin{scope}[xshift = 1.33333cm]
				\draw plot [smooth cycle] coordinates {(0:\rb) (60:\rs) (120:\rb) (180:\rs) (240:\rb) (300:\rs) };
			\end{scope}
			\begin{scope}[xshift = 0.666667cm, rotate=60]
				\draw plot [smooth cycle] coordinates {(0:\rb) (60:\rs) (120:\rb) (180:\rs) (240:\rb) (300:\rs) };
			\end{scope}
			\draw[opacity=0.2] ({2/3},0) -- ({4/3},0)
			($0.5*(1,{1/sqrt(3)})$) -- ({2/3},0)
			($0.5*(1,{-1/sqrt(3)})$) -- ({2/3},0)
			($(1,{1/sqrt(3)})+0.5*(1,{-1/sqrt(3)})$) -- ({4/3},0)
			($0.5*(1,{1/sqrt(3)})+(1,{-1/sqrt(3)})$) -- ({4/3},0);
		\end{scope}
		\begin{scope}[shift = ($(a)-(b)$)]
			\begin{scope}[xshift = 1.33333cm]
				\draw plot [smooth cycle] coordinates {(0:\rb) (60:\rs) (120:\rb) (180:\rs) (240:\rb) (300:\rs) };
			\end{scope}
			\begin{scope}[xshift = 0.666667cm, rotate=60]
				\draw plot [smooth cycle] coordinates {(0:\rb) (60:\rs) (120:\rb) (180:\rs) (240:\rb) (300:\rs) };
			\end{scope}
			\draw[opacity=0.2] ({2/3},0) -- ({4/3},0)
			($0.5*(1,{1/sqrt(3)})$) -- ({2/3},0)
			($0.5*(1,{-1/sqrt(3)})$) -- ({2/3},0)
			($(1,{1/sqrt(3)})+0.5*(1,{-1/sqrt(3)})$) -- ({4/3},0)
			($0.5*(1,{1/sqrt(3)})+(1,{-1/sqrt(3)})$) -- ({4/3},0);
		\end{scope}
		\begin{scope}[shift = ($-1*(a)+(b)$)]
			\begin{scope}[xshift = 1.33333cm]
				\draw plot [smooth cycle] coordinates {(0:\rb) (60:\rs) (120:\rb) (180:\rs) (240:\rb) (300:\rs) };
			\end{scope}
			\begin{scope}[xshift = 0.666667cm, rotate=60]
				\draw plot [smooth cycle] coordinates {(0:\rb) (60:\rs) (120:\rb) (180:\rs) (240:\rb) (300:\rs) };
			\end{scope}
			\draw[opacity=0.2] ({2/3},0) -- ({4/3},0)
			($0.5*(1,{1/sqrt(3)})$) -- ({2/3},0)
			($0.5*(1,{-1/sqrt(3)})$) -- ({2/3},0)
			($(1,{1/sqrt(3)})+0.5*(1,{-1/sqrt(3)})$) -- ({4/3},0)
			($0.5*(1,{1/sqrt(3)})+(1,{-1/sqrt(3)})$) -- ({4/3},0);
		\end{scope}
	\end{scope}
\end{tikzpicture}
\end{subfigure}
\begin{subfigure}[b]{0.4\linewidth}
	\centering
\begin{tikzpicture}[scale=1.5]
	\coordinate (a) at ({1/sqrt(3)},1);
	\coordinate (b) at ({1/sqrt(3)},-1);
	\coordinate (c) at ({2/sqrt(3)},0);
	\coordinate (K1) at ({1/sqrt(3)},{1/3});
	\coordinate (K2) at ({1/sqrt(3)},{-1/3});
	\coordinate (K3) at (0,{-2/3});
	\coordinate (K4) at ({-1/sqrt(3)},{-1/3});
	\coordinate (K5) at ({-1/sqrt(3)},{1/3});
	\coordinate (K6) at (0,{2/3});

	\draw[fill] (K1) circle(1pt) node[xshift=6pt,yshift=4pt]{$\alpha_1^*$};
	\draw[fill] (K2) circle(1pt) node[xshift=6pt,yshift=-4pt]{$\alpha_2^*$};
	\draw[fill] (0,0) circle(1pt) node[left]{$\Gamma$};

	\draw (K1) -- (K2) -- (K3) -- (K4) -- (K5) node[above]{$Y^*$} -- (K6) -- cycle;
\end{tikzpicture}
\vspace{10pt}
\end{subfigure}
\caption{Honeycomb crystal and corresponding Brillouin zone.}\label{fig:honeycomb}
\end{figure}

We consider a two-dimensional infinite honeycomb crystal in two dimensions depicted in Figure \ref{fig:honeycomb}. We let the lattice $\Lambda$ be generated by the lattice vectors
$$ l_1 = L\left( \frac{\sqrt{3}}{2}, \frac{1}{2} \right),~~l_2 = L\left( \frac{\sqrt{3}}{2}, -\frac{1}{2}\right).$$
We assume that each unit cell contains two resonators, $D=D_1\cup D_2$, such that each resonator is invariant under rotation by $2\pi/3$ and so that $D$ is invariant under rotation by $\pi$.

The dual lattice  $\Lambda^*$ is generated by $\alpha_1$ and $\alpha_2$ given by
$$ \alpha_1 = \frac{2\pi}{L}\left( \frac{1}{\sqrt{3}}, 1\right),~~\alpha_2 = \frac{2\pi}{L}\left(\frac{1}{\sqrt{3}}, -1 \right).$$
The points $$\alpha_1^*= \frac{2\alpha_1+\alpha_2}{3}, \quad \alpha^*_2 = \frac{\alpha_1+2\alpha_2}{3},$$ in the Brillouin zone are called \emph{Dirac points}. Next, we will study the band functions and Bloch modes around these points. For simplicity, we only consider the analysis around the Dirac point $\alpha_* := \alpha_1^*$, the other point having a similar behaviour.

\begin{figure}[tbh]
	\centering
	\begin{subfigure}[b]{0.48\linewidth}
		\centering
		\includegraphics[height=5.5cm]{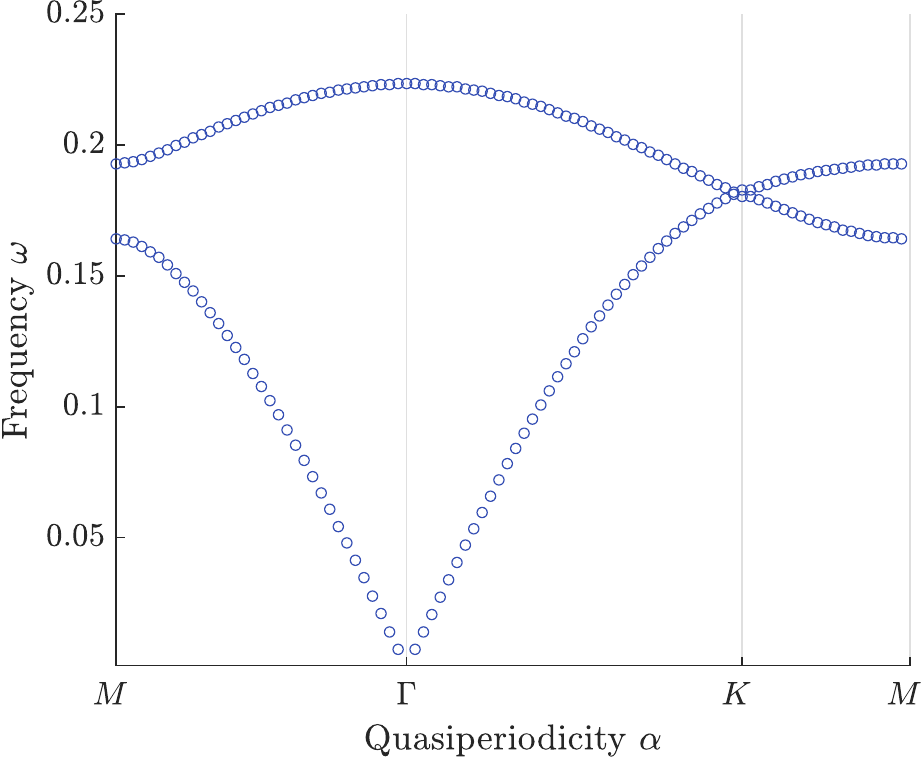}
		\caption{Subwavelength band structure of a honeycomb structure, exhibiting a Dirac cone at $\alpha = K$.}
		\label{fig:cone}
	\end{subfigure}\hfill
	\begin{subfigure}[b]{0.48\linewidth}
		\centering
		\includegraphics[height=6cm]{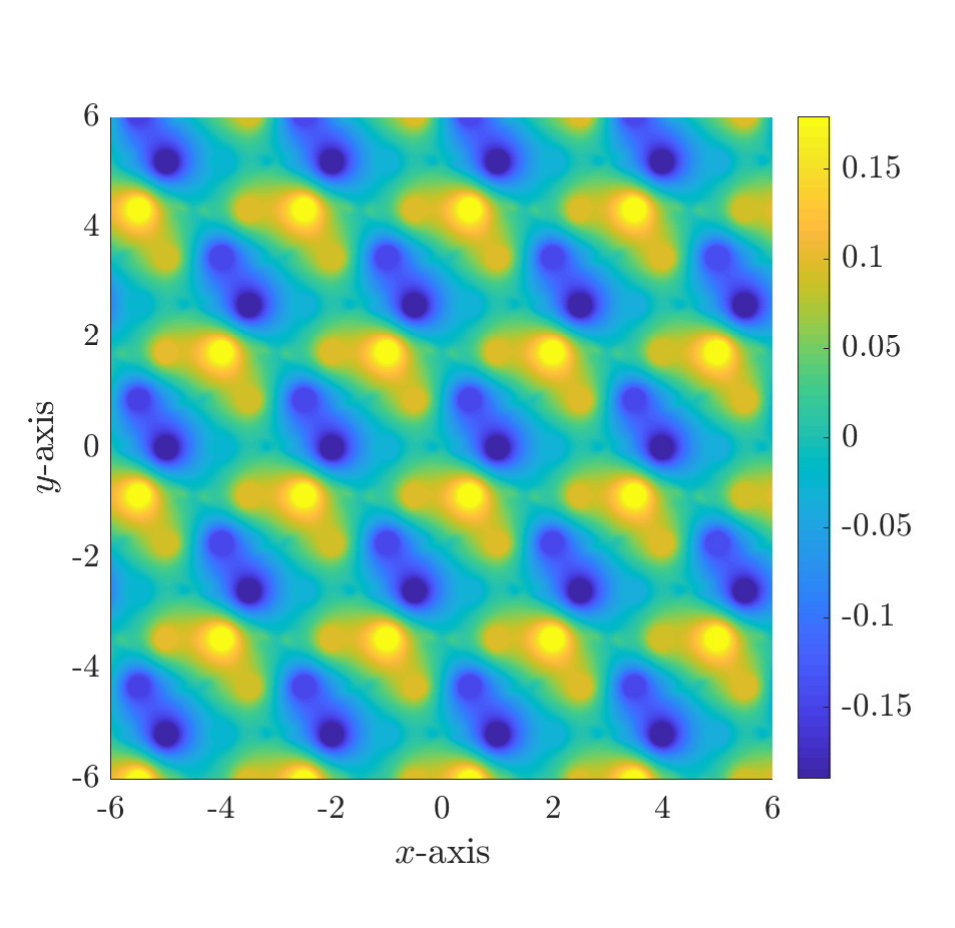}
		\caption{Small-scale behaviour of the eigenmodes at the Dirac point.}
		\label{fig:Smr}
	\end{subfigure}
	\caption{The honeycomb structure exhibits a Dirac cone at the corner of the Brillouin zone. Corresponding eigenmodes are rapidly oscillating, and is periodic across one hexagon in the honeycomb structure.} \label{fig:smallscale}
\end{figure}

\begin{figure}[tbh]
	\begin{subfigure}[b]{0.47\linewidth}
		\hspace{-10pt}
		\includegraphics[width=\linewidth]{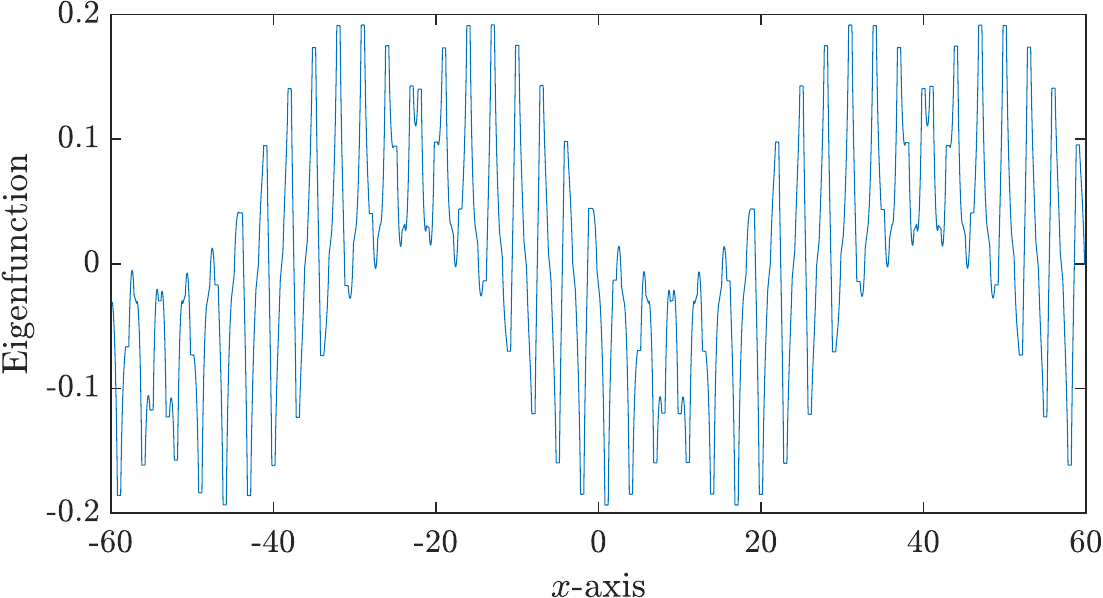}
		\caption{Large-scale behaviour of the eigenmodes close to the Dirac point. }
		\label{fig:h1D}
	\end{subfigure}\hfill
	\begin{subfigure}[b]{0.49\linewidth}
		\includegraphics[width=\linewidth]{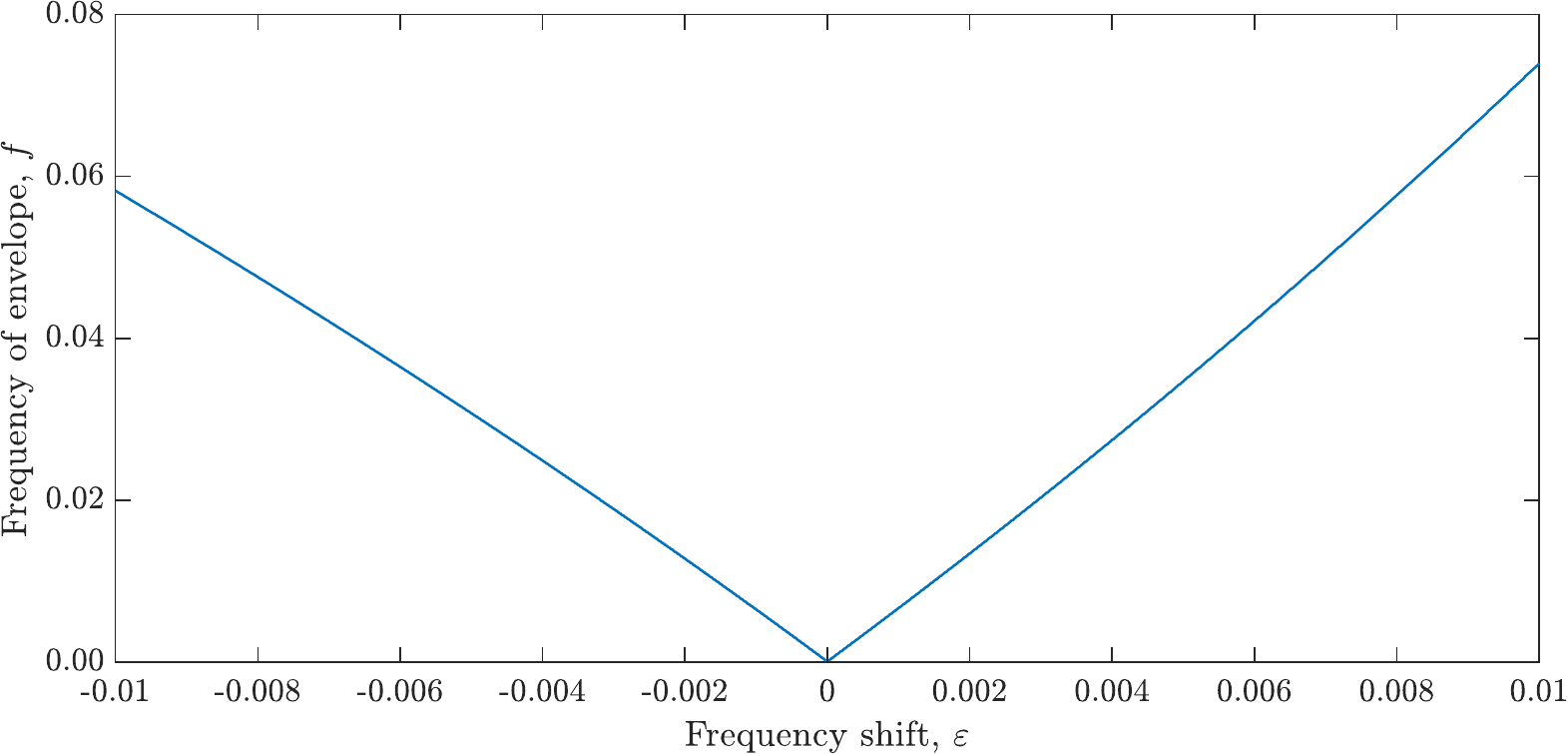}
		\caption{Spatial frequency $f$ of the envelopes as a function of the frequency shift $\varepsilon=\omega-\omega_*$.}
	\end{subfigure}
	\caption{For frequencies close to the Dirac frequency, the eigenfunctions oscillate on two distinct scales. The spatial frequency of the large-scale oscillations vanish at the centre of the Dirac cone, corresponding to zero refractive index.} \label{fig:largescale}
\end{figure}

At $\alpha = \alpha_*$, the generalised capacitance matrix $\C^\alpha$ has an eigenvalue of multiplicity 2: $\lambda_1^{\alpha_*} = \lambda_2^{\alpha_*}$. The next result shows that this asymptotic degeneracy is in fact an exact degeneracy, and moreover that the band functions intersect in a Dirac cone at this point (see \Cref{fig:cone}) \cite{ammari2020honeycomb}.
\begin{thm}\label{thm:honeycomb}
	For $\alpha$ close to $\alpha_*$ and $\delta$ small enough, the first two band functions form a Dirac cone, \ie{},
	\begin{equation} \label{eq:dirac}
		\begin{matrix}
			\ds \omega_1^\alpha = \omega_*- \mu|\alpha - \alpha_*| \big[ 1+ O(|\alpha-\alpha_*|) \big], \\[0.5em]
			\ds \omega_2^\alpha = \omega_*+ \mu|\alpha - \alpha_*| \big[ 1+ O(|\alpha-\alpha_*|) \big],
		\end{matrix}
	\end{equation}
	where $\omega_*$ and $\mu$ are independent of $\alpha$ and satisfy
	$$\omega_*= \sqrt{\lambda_1^{\alpha_*}} + O(\delta) \quad \text{and} \quad \mu =  |c|\sqrt\delta\mu_0 + O(\delta),  \quad \mu_0=\frac{1}{2}\sqrt{\frac{v_\mathrm{r}^2 }{|D_1|C_{11}^{\alpha_*}}}, \quad c = \left|\frac{\partial C_{12}^{\alpha}}{\partial \alpha_1}\Big|_{\alpha=\alpha_*}\right|,$$
	as $\delta \rightarrow 0$. Moreover, the error term $O(|\alpha-\alpha_*|)$ in \eqref{eq:dirac} is uniform in $\delta$.
\end{thm}
Next, we will see that the near-zero property follows as a direct consequence of the Dirac cone. We consider a homogenization setting close to the Dirac points. We rescale the unit cell by replacing $Y$ with $sY$ for some small $s>0$. To fix the order of the resonant frequencies, \ie{} $\omega_i^\alpha = O(1)$, we assume that $\delta = O(s^2)$ as $s\to 0$. We then have the following result \cite{ammari2020highfrequency}.
\begin{thm} \label{thm:nearzero}
	For frequencies $\omega$ close to the Dirac frequency $\omega_*$, namely, $\omega-\omega_* = \beta \sqrt\delta$, the following asymptotic behaviour of the Bloch eigenfunction $u^{\alpha_*/s + \tilde{\alpha}}_s$ holds:
	$$
	u_{s}^{\alpha_*/s+\tilde{\alpha}}(x) = \begin{bmatrix}A e^{\i\tilde{\alpha}\cdot x}\\ B e^{\i\tilde{\alpha}\cdot x}\end{bmatrix}\cdot \textbf{S}_D^{\alpha_*,k}\left(\frac{x}{s}\right) + O(s),
	$$
	where the macroscopic field $[\tilde{u}_{1}, \tilde{u}_{2}]^T:=[A e^{\i\tilde{\alpha}\cdot x}, B e^{\i\tilde{\alpha}\cdot x}]^T$ satisfies the two-dimensional Dirac equation
	$$
	\mu_0
	\begin{bmatrix}
		0 & (-c\i)( \p_1 - \i \p_2)
		\\
		(-\overline{c}\i)( \p_1 + \i \p_2) & 0
	\end{bmatrix}
	\begin{bmatrix}
		\tilde{u}_{1}
		\\
		\tilde{u}_{2}
	\end{bmatrix}
	= \frac{\omega-\omega_*}{\sqrt\delta}
	\begin{bmatrix}
		\tilde{u}_{1}
		\\
		\tilde{u}_{2}
	\end{bmatrix}.
	$$
\end{thm}
The system of Dirac equations can be considered as a homogenized equation for the honeycomb structure. Each $\tilde u_j$ satisfy the Helmholtz equation
\begin{equation} \label{eq:hom}
	\Delta \tilde{u}_j + \frac{(\omega-\omega_*)^2}{\mu^2} \tilde{u}_j = 0.
\end{equation}
In particular, at $\omega=\omega_*$ this effective equation reduces to the Laplace equation corresponding to effective zero refractive index. Equation \eqref{eq:hom} describes the large-scale behaviour of the eigenmodes, illustrated in \Cref{fig:largescale}. We emphasize that in addition to this large-scale behaviour, there will be small-scale oscillations described by the functions $\textbf{S}_D^{\alpha_*,k}$ as illustrated in \Cref{fig:smallscale}.

\begin{remark}
The results in \Cref{thm:honeycomb} and \Cref{thm:nearzero} 
are derived from first principles. They can be viewed as the classical wave analogues of those proved in \cite{feffermandirac} for
Schr\"odinger operator. The time-evolution of wave-packets in honeycomb systems of subwavelength resonators, which are spectrally concentrated near conical points, can be studied in the same way as in  \cite{feffermandirac,yi1,yi2}. 
\end{remark}

\begin{remark}
The results on the two-dimensional (double-degenerate) Dirac points in \Cref{thm:honeycomb} and \Cref{thm:nearzero} can be extended to three dimensions and higher-order degeneracies \cite{yi3,yi4}.
\end{remark}

\subsubsection{Double-near zero materials in time-modulated systems}

\begin{figure}[h]
	\begin{center}
		\includegraphics[width=0.75\linewidth]{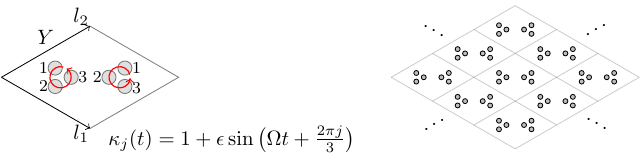}
	\end{center}
	\caption{Illustration of the trimer honeycomb lattice, with phase-shifted time-modulations inside the trimers.} \label{fig:trihoney}
\end{figure}

\begin{figure}[h]
	\begin{center}
		\begin{subfigure}[t]{0.4\linewidth}
			\vspace{0pt}
			\begin{center}
				\begin{tikzpicture}
					\draw (0.25,1.38) node{\includegraphics[width=1\linewidth]{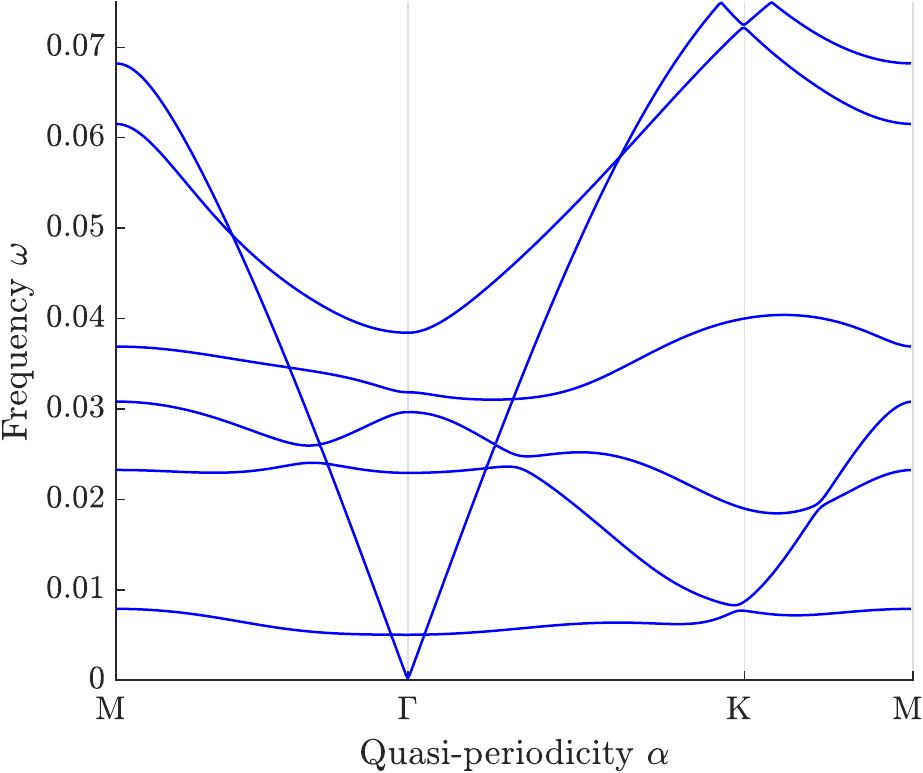}};
					%\draw[draw=none,fill,opacity=0.3] (0,0) -- (1.2,3.2) -- (-1.2,3.2) -- cycle;
					%\draw (0,3) node{\tiny first radiation};
					%\draw (0,2.8) node{\tiny continuum};
				\end{tikzpicture}
			\end{center}
			\caption{Static band structure}
		\end{subfigure}
		\hspace{20pt}
		\begin{subfigure}[t]{0.4\linewidth}
			\vspace{0pt}
			\begin{center}
				\begin{tikzpicture}
					\draw (0.25,1.38) node{\includegraphics[width=1\linewidth]{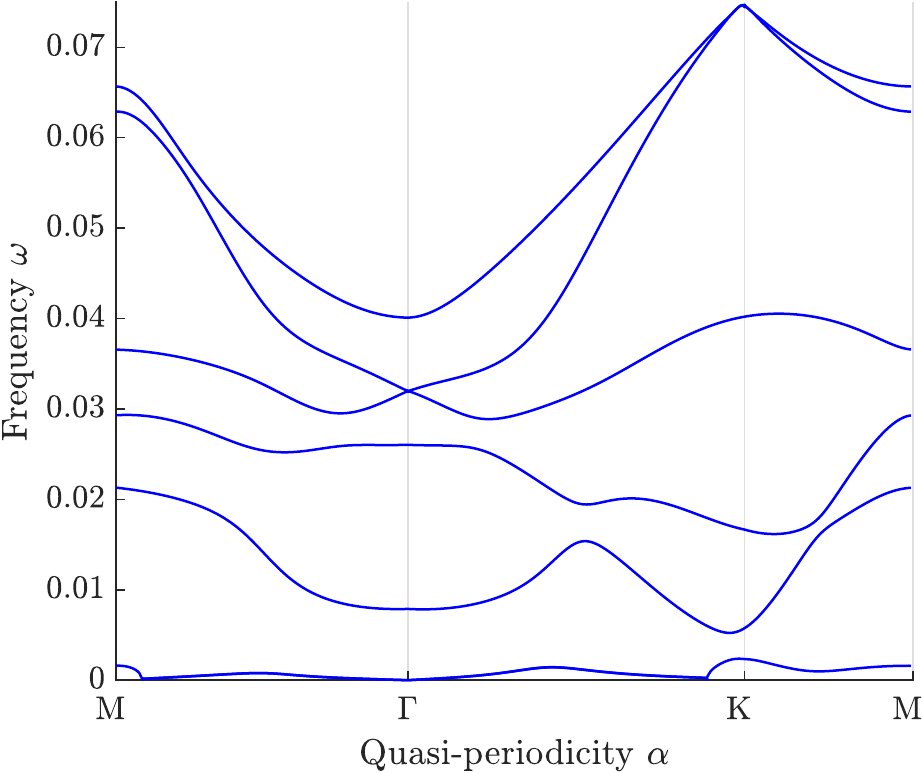}};
					%\draw[draw=none,fill,opacity=0.3] (0,0) -- (1.2,3.2) -- (-1.2,3.2) -- cycle;
					%\draw (0,3) node{\tiny first radiation};
					%\draw (0,2.8) node{\tiny continuum};
				\end{tikzpicture}
			\end{center}
			\caption{Modulated band structure}
		\end{subfigure}
	\end{center}
	\caption{In time-modulated structures, Dirac cones may appear at the origin of the Brillouin zone. Such points are associated to double near-zero materials, enabling wave transmission without phase changes and with strong interactions with incoming waves.}\label{fig:diracGamma}
\end{figure}

The Dirac cone observed in the previous section is located at the corner of the Brillouin zone, which means that corresponding modes do not lie in the radiation continuum. For physical structures of finite dimensions, these waves will be confined to the structure and will not interact with incoming waves. Physically, this can be seen as an impedance mismatch, leading to poor transmission, due to the fact that this structure corresponds to a \emph{single}-near zero material, where only one of the material parameters vanish.

In order to overcome this low transmission, it is desirable to create structures with Dirac cones at the origin of the Brillouin zone. Corresponding modes lie in the radiation continuum, and will therefore interact with incoming waves. In order to achieve this, we study a time-modulated honeycomb structure as illustrated in \Cref{fig:trihoney}. Using the theory from \Cref{sec:time} we can compute the (quasi-)band structure associated to this material. In \Cref{fig:diracGamma} we see the static (folded) band structure, with six bands in the subwavelength regime. For a particular modulation strength, we see that a Dirac cone degeneracy may appear at the origin of the Brillouin zone, enabling a double-near zero material around this point.

\section{Concluding remarks}

In this survey, we have reviewed several Helmholtz scattering problems posed in the subwavelength regime and repeatedly encountered the concept of capacitance. We have studied the mathematical properties of the generalised capacitance matrix, in both finite and infinite, periodic settings. We were then able to use the capacitance formulation to study several different interesting physical phenomena in the field of subwavelength metamaterials. This formulation emerged from a functional analytic approach, where the generalised capacitance matrix describes the perturbation of the kernel of a non-linear integral operator for asymptotically small parameter values. Similar approaches apply to a plethora of other subwavelength scattering problems, including high-contrast dielectric particles, plasmonic nanoparticles and Helmholtz resonators \cite{ammari2019dielectric, ammari2020mathematical, ammari2017plasmonicscalar, ammari2016plasmonicMaxwell, ammari2015superresolution}. Suitable capacitance formulations can thereby be used to characterize a wide range of subwavelength resonance phenomena.

Approximating classical wave systems in terms of generalised capacitance matrices shares similarities with the tight-binding approximation that is used widely in quantum theory, in the sense that both formulations provide a discrete approximation to a continuous differential problem. As we have seen, there are nevertheless fundamental differences, for instance due to the strong interactions between the subwavelength resonators. As observed in \Cref{rmk:tightbind}, the correspondence between the capacitance formulation and the tight-binding approximation holds only in the case of dilute resonators. Even in the dilute case, long-range interactions between subwavelength resonators are relatively strong and nearest-neighbour approximations are not generally appropriate. The strength of the capacitance formulation is that the capacitance matrix accounts for these strong interactions, thereby providing a unified mathematical model to study challenging problems in subwavelength physics. Recently, the capacitance matrix formulation has been used to reveal new insight into the mechanisms responsible for the fundamental features of Anderson localisation in systems of subwavelength resonators with randomly chosen material parameters. In \cite{anderson}, it is shown that the generalised capacitance matrix provides a characterisation of the localised modes and captures the long-range interactions of the wave-scattering system. This provides a rigorous framework to explain the exotic phenomena that are observed. 
On the other hand,  in \cite{skineffect} the capacitance matrix formulation is generalised to non-Hermitian problems with imaginary gauge potentials in order to prove the condensation of eigenmodes at one of the edges of a finite chain of subwavelength resonators with an imaginary gauge potential supported inside the resonators. 
Many of the approaches and results discussed in this survey have also recently found applications and experimental validation in the fields of applied physics and engineering; see, for instance, \cite{experimental1,experimental2,experimental3,experimental4,experimental5,experimental6,experimental7,experimental8,experimental9,experimental10,experimental11,experimental12,experimental13,experimental14,experimental15,experimental16,experimental17,experimental18}.

\appendix
\section{Abstract capacitance matrix}\label{app:abstract}
In this section, we describe how the capacitance formulation emerges from the structure of a general integral operator $\A$, which describes a subwavelength resonance problem. We work with the functional analytic approach described in \Cref{sec:approach}. At $\delta=0$, we assume that the operator $\A(\omega,0)$ has a characteristic value $\omega=0$ of multiplicity $2N$, admitting the following pole-pencil decomposition:
\begin{equation}\label{eq:polepencil}
	\A(\omega,0)^{-1} = \frac{K}{\omega^2} + \Rc(\omega), \quad \text{ for } \quad K = \sum_{i=1}^N \langle \Phi_i, \cdot\rangle \Psi_i,
\end{equation}
where $\ker(\A(0,0)) = \mathrm{span}\{\Psi_j\}$, $\ker(\A^*(0,0)) = \mathrm{span}\{\Phi_j\}$ and $\Rc$ is holomorphic for $\omega$ in a neighbourhood of $0$.
Moreover, we assume that $\A(\omega,\delta)$, for small but non-zero $\delta$, satisfies
$$\A(\omega,\delta) = \A(\omega,0) + \L(\omega,\delta),$$
for some operator $\L$ satisfying (in corresponding operator norm) $\|\L\| = O(\delta)$ uniformly for $\omega$ in a neighbourhood of $0$.

In this abstract setting, we can derive a capacitance formulation of the subwavelength resonances. We solve the equation $\A(\omega,\delta)\Phi = 0$. Multiplying with $\A(\omega,0)^{-1}$, we have
\begin{align*}
	0 &= \A(\omega,0)^{-1}\A(\omega,\delta)\Phi = \A(\omega,0)^{-1}\left(\A(\omega,0) + \L\right)\Phi = \left(I + \frac{K\L}{\omega^2} + \Rc \L\right)\Phi.
\end{align*}
Defining $\mathcal{B}(\omega,\delta) = \omega^2\Rc(\omega) \L(\omega,\delta)$ yields
\begin{equation}\label{eq:eigen}
	\left(\omega^2I + K\L + \mathcal{B}\right)\Phi = 0.
\end{equation}
The characteristic values are therefore determined by \eqref{eq:eigen}, which is in general a non-linear eigenvalue problem since $\L$ and $\mathcal{B}$ depend on $\omega$. If we consider the subwavelength resonances, we have $\|\mathcal{B}\| = O(\omega^2\delta)$ uniformly for $\omega$ and $\delta$ around $0$. Similarly, we have $\L = \L_0 + \hat\L$, where $\hat{\L}=O(\omega\delta)$. Therefore, the subwavelength resonances are approximated by the eigenvalues of the finite-rank operator $-K\L_0$ whose restriction to $\ker(\A(0,0))$ is given by the generalised capacitance matrix:
$$\C_{ij} =  -\langle \Phi_i, \L_0\Psi_j \rangle.$$
Then the characteristic values satisfy
$$\omega_n = \pm \sqrt{\lambda_n} + O(\delta),$$
where $\lambda_n$ are the eigenvalues of $\C$.

As an example of this formulation, for a finite collection of $N$ resonators in $d=3$ (as considered in \Cref{sec:finite}), the operator $\A$ is given by
\begin{equation*}
	\A(\omega,\delta)= \begin{pmatrix}
		\widetilde{\S}_D^\omega & -\S_D^k \\
		-\frac{1}{2}I+\widetilde{\K}_D^{\omega,*} & -\widetilde\delta \left(\frac{1}{2}I+\K_D^{k,*}\right)
	\end{pmatrix}.
\end{equation*}
With above notation, we then have that
\begin{equation*}
	\A(\omega,0)= \begin{pmatrix}
		\widetilde{\S}_D^\omega & -\S_D^k \\
		-\frac{1}{2}I+\widetilde{\K}_D^{\omega,*}  & 0
	\end{pmatrix}, \quad \L(\omega,\delta) = \begin{pmatrix}
	0 & 0 \\
	0 & -\widetilde\delta \left(\frac{1}{2}I+\K_D^{k,*}\right)
\end{pmatrix}, \quad \L_0(\delta) = \begin{pmatrix}
0 & 0 \\
0 & -\widetilde\delta \left(\frac{1}{2}I+\K_D^{0,*}\right)
\end{pmatrix}.
\end{equation*}
Moreover, it can be shown that $\A(\omega,0)^{-1}$ satisfies \eqref{eq:polepencil} where
$$\Phi_i = -\frac{v_i^2}{|D_i|}\begin{pmatrix} 0 \\ \chi_{\p D_i} \end{pmatrix},\qquad \Psi_j = \begin{pmatrix} \psi_i \\ \psi_i \end{pmatrix}, \quad \psi_i = (\S_D^0)^{-1}[\chi_{\p D_i}].$$
From this, it is straightforward to compute $\C_{ij}=\frac{\delta_iv_i^2}{|D_i|}\langle\chi_{\p D_i}, \psi_j\rangle$, as defined in \eqref{eq:GCM} in \Cref{defn:GCM}.

\bibliographystyle{abbrv}
\bibliography{capacitance}

\end{document}